\documentclass[11pt]{article}
\usepackage{amssymb,amsmath,amsthm,tikz,multirow,nccrules,float,colortbl,arydshln,multicol,ulem,graphicx,subfig}
\usetikzlibrary{arrows,calc}

\title{Tilings of the sphere by congruent quadrilaterals \uppercase\expandafter{\romannumeral2}: edge combination $a^3b$ with rational angles}
\author{Yixi Liao, Erxiao Wang\thanks{Corresponding author (wang.eric@zjnu.edu.cn).  Research was supported by Key projects of Zhejiang Natural Science Foundation No. LZ22A010003 and ZJNU Shuang-Long Distinguished Professorship Fund No. YS304319159.} \\
	Zhejiang Normal University}

\allowdisplaybreaks[4]  

\newcommand{\sub}{\subset}

\newcommand\aaa{\alpha}
\newcommand\bbb{\beta}
\newcommand\ccc{\gamma}
\newcommand\ddd{\delta}

\newcommand{\thin}{\hspace{0.1em}\rule{0.7pt}{0.8em}\hspace{0.1em}}
\newcommand{\thick}{\hspace{0.1em}\rule{1.5pt}{0.8em}\hspace{0.1em}}

\newtheorem{theorem}{Theorem}
\newtheorem{lemma}[theorem]{Lemma}
\newtheorem{remark}[theorem]{Remark}

\newtheorem{proposition}[theorem]{Proposition}
\newtheorem*{theorem*}{Theorem}

\theoremstyle{definition}
\newtheorem*{definition*}{Definition}
\newtheorem*{case*}{Case}
\newtheorem*{subcase*}{Subcase}

\numberwithin{equation}{section}

\begin{document}
	\date{}
	
	\maketitle	
\begin{abstract}
	Edge-to-edge tilings of the sphere by congruent quadrilaterals are completely classified in a series of three papers. This second one applies the powerful tool of trigonometric Diophantine equations to classify the case of $a^3b$-quadrilaterals with all angles being rational degrees. There are $12$ sporadic and $3$ infinite sequences of quadrilaterals admitting the $2$-layer earth map tilings together with their modifications,  and $3$ sporadic quadrilaterals admitting $4$ exceptional tilings. Among them only $3$ quadrilaterals are convex. New interesting  non-edge-to-edge triangular tilings are obtained as a byproduct. 
	
	
	{\it2000 MR Subject Classification} 52C20, 05B45
\end{abstract}

\section{Introduction}

In an edge-to-edge tiling of the sphere by congruent quadrilaterals, the tile can only have four edge arrangements \cite{ua2,lpwx1}:  $a^2bc,a^2b^2,a^3b,a^4$. Sakano and Akama \cite{sa} classified  tilings for $a^2b^2$ and $a^4$ via Ueno and Agaoka's \cite{ua} list of triangular tilings. Tilings for $a^2bc$ are classified in the first paper \cite{lpwx1} of this series via the methods in \cite{wy1,wy2,awy,wy3} developed for pentagonal tilings.  This second paper classifies tilings for $a^3b$ with all angles being rational multiples of $\pi$ (such quadrilaterals will be simply called \textit{rational} hereafter). We then classify tilings for $a^3b$ with some irrational angle in the third paper  \cite{lpwx2} to complete the classification. 

Recall that Coolsaet \cite{ck} classified convex rational quadrilaterals with three equal sides into $7$ infinite classes and $29$ sporadic examples. Akama and Cleemput \cite{ac} initiated some explorations of degree $3$ vertex types and certain forbidden cases for type $a^3b$, assuming also  convexity. 

\begin{figure}[htp]
	\centering
	\begin{tikzpicture}
		
		\fill[gray!50]  (-0.8+3,-0.8) -- (-0.8+3,0.8) -- (0.8+3,0.8) -- (0.8+3,-0.8);
		
		\foreach \a in {0,1}
		{
			\begin{scope}[xshift=3*\a cm]
				\draw
				(-0.8,-0.8) -- (-0.8,0.8) -- (0.8,0.8) -- (0.8,-0.8);

				\draw[line width=1.5]
				(-0.8,-0.8) -- (0.8,-0.8);	    	\end{scope}}   
		
		\draw  (4.5+0.5,0.1)--(5.5+0.5,0.1);

		\draw[line width=1.5] (4.5+0.5,-0.4)--(5.5+0.5,-0.4);
		
		\node at (-0.5,0.5) {\small $\bbb$};
		\node at (0.5,0.5) {\small $\ccc$};
		\node at (-0.5,-0.5) {\small $\aaa$};
		\node at (0.5,-0.5) {\small $\ddd$};
		
		\foreach \a in {0,1,2}
		\node at (90*\a:1) {\small $a$};
		\node at (0,-1.1) {\small $b$};
		
		\node at (3,1) {\small $a$};\node at (2,0) {\small $a$};\node at (4,0) {\small $a$};
		\node at (3,-1.1) {\small $b$};
		
		\node at (2+3.5,0.3) {\small $a$};
		
		\node at (2+3.5,-0.2) {\small $b$};

		\node at (-0.5+3,0.5) {\small $\ccc$};
		\node at (0.5+3,0.5) {\small $\bbb$};
		\node at (-0.5+3,-0.5) {\small $\ddd$};
		\node at (0.5+3,-0.5) {\small $\aaa$};
		
	\end{tikzpicture} 
	\caption{Quadrilaterals with the edge combination $a^3b$.}
	\label{quad}
\end{figure}
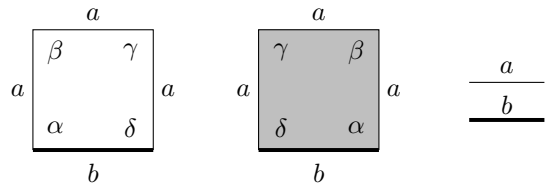

An $a^3b$-quadrilateral is given by Fig.\,\ref{quad}, with normal edge $a$, thick edge $b$ and angles $\alpha,\beta,\gamma,\delta$ as indicated. The second picture is the mirror image or flip of the first. The angles determine the orientation. Conversely, the edge lengths and the orientation also determine the angles. So we may present the tiling by shading instead of indicating all angles. Throughout this paper, an $a^3b$-tiling is always an edge-to-edge tiling of the sphere by congruent simple quadrilaterals in Fig.\,\ref{quad}, such that all vertices have degree $\ge 3$.

	\begin{theorem*}
There are $15$ sporadic and $3$ infinite sequences of rational quadrilaterals which admit $a^3b$-tilings (Table \ref{Tab-1.1} and \ref{Tab-1.2}). Except the last $3$ sporadic cases, they are all $2$-layer earth map tilings $T(f\aaa\bbb\ddd,2\ccc^{\frac{f}{2}})$ for some even integers $f\ge6$, together with their modifications when $\beta$ is an integer multiple of $\gamma$. The total number $\mathcal{Q}(f)$ of quadrilaterals in Table \ref{Tab-1.1} and \ref{Tab-1.2} and their total number $\mathcal{T}(f)$ of tilings are: 
\begin{table*}[htp]               
	\centering     
	\resizebox{\textwidth}{10mm}{\begin{tabular}{c|cccccccccccccc}	 
			
			$f$ &$6,30$&$8$&$12$&$16$&$18$&$20$&$36$&$12k$&$12k+2$ &$12k+4$&$12k+6$&$12k+8$&$12k+10$ \\
			\hline
			$k$&&&&&&&&$2,\ge4$&$\ge1$&$\ge2$&$\ge3$&$\ge2$&$\ge0$\\
			\hline
			$\mathcal{Q}(f)$&$4$&$1$&$8$&$4$&$4$&$5$&$5$&$3$&$3$&$3$&$3$&$3$&$3$ \\
			\hline			
			$\mathcal{T}(f)$&$4$&$1$&$12$&$14$&$6$&$13$&$8$&$6$&$k+6$&$k+11$&$3$&$k+10$&$k+8$ \\		
			\hline	
	\end{tabular}}
\end{table*}
	\end{theorem*}

	\begin{table*}[htp]                        
		\centering      
		\resizebox{\textwidth}{48mm}{\begin{tabular}{c|c|c|c}	 
				
				$f$ & $ (\aaa,\bbb,\ccc,\ddd),a,b$&Page & all vertices and tilings \\
				\hline 
				\multirow{3}{*}{$6$} & $(6,3,4,3)/6,1/2,1/6$&\pageref{discrete-1}&\multirow{3}{*}{$6\aaa\bbb\ddd,2\ccc^3$ }  \\
				&$(1,8,4,3)/6,0.391,1$&\pageref{discrete-2'},\pageref{discrete-2},\pageref{discrete-2''}&\\
				&$(12,4,6,2)/9,0.567,0.174$&\pageref{discrete-3}&\\
				\cline{1-4}
				\multirow{5}{*}{$12$}
				& $(2,10,3,6)/9,0.339,0.532$ &\pageref{discrete-4}& \multirow{5}{*}{$12\aaa\bbb\ddd,2\ccc^6$}\\				
				&$(1,21,5,8)/15,0.424,0.741$&\pageref{discrete-5}&\\
				&$(4,9,5,17)/15,0.424,0.165$&\pageref{discrete-6}&\\ 
				&$(9,28,10,23)/30,0.335,0.415$&\pageref{discrete-7}&\\ 
				&$(3,16,10,41)/30,0.469,0.146$&\pageref{discrete-8}&\\
				\hline   		
				\multirow{2}{*}{$20$} & $(5,32,6,23)/30,0.335,0.415$ &\pageref{discrete-9}& \multirow{2}{*}{$20\aaa\bbb\ddd,2\ccc^{10}$}\\
				&$(1,16,6,43)/30,0.469,0.273$&\pageref{discrete-10}&\\
				\cline{1-4}		
				$30$&$(1,42,4,17)/30,0.424,0.549$&\pageref{discrete-11}&$30\aaa\bbb\ddd,2\ccc^{15}$\\
				\hline 
				\hline 
				\multirow{3}{*}{$18$}&\multirow{3}{*}{$(3,20,4,13)/18,0.339,0.452$}&\pageref{discrete-12}&$18\aaa\bbb\ddd,2\ccc^9$\\
				&&\pageref{discrete-13}&$16\aaa\bbb\ddd,2\bbb\ccc^4,2\aaa\ccc^5\ddd$\\
				&&\pageref{discrete-14}&$14\aaa\bbb\ddd,2\aaa^2\ccc\ddd^2,4\bbb\ccc^4$\\
				\hline
				\hline 
				$16$&$(1,4,2,2)/4,1/4,1/2$&\pageref{discrete-15}&$8\bbb\ddd^2,8\aaa^2\bbb\ccc,2\ccc^4$: $2$ tilings\\
				\hline
				$36$&$(5,4,7,3)/9,0.174,0.258$&\pageref{discrete-16}&$18\bbb\ccc^2,6\aaa^3\ddd,6\aaa^2\bbb^2,6\aaa\bbb\ddd^3,2\ddd^6$\\
				\hline 												
				$36$&$(15,6,10,7)/18,0.225,0.118$&\pageref{discrete-17}&$14\aaa^2\bbb,8\aaa\ddd^3,10\bbb\ccc^3,6\bbb^2\ccc\ddd^2$\\
				\hline 
				
		\end{tabular}}
	\caption{Fifteen sporadic quadrilaterals and their tilings.}\label{Tab-1.1}        
	\end{table*}

	\begin{table*}[htp]

		\centering     		
		\resizebox{\textwidth}{35mm}{\begin{tabular}{c|c|c}	 
			
			$(\aaa,\bbb,\ccc,\ddd)$&all vertices and tilings & Page\\
			\hline 			 
			\multirow{4}{*}{ $(\frac 4f,1-\frac{4}{f},\frac 4f,1)$}&$\forall$ even $f\ge10: f\aaa\bbb\ddd,2\ccc^{\frac{f}{2}}$ &\pageref{a=1} \\
            &$f=4k(k\ge3)$: $(f-2)\aaa\bbb\ddd,2\aaa\ccc^{\frac f4-1}\ddd,2\bbb\ccc^{\frac f4+1}$& \pageref{a=1}\\
			&$(f-4)\aaa\bbb\ddd,2\bbb^2\ccc^2,4\aaa\ccc^{\frac f4-1}\ddd$: $2$ tilings&\pageref{a=1}\\
			&$f=12$: $6\aaa\bbb\ddd,2\bbb^3,6\aaa\ccc^2\ddd$&\pageref{a=1}\\
			\hline 			 
			\multirow{4}{*}{ $(\frac 2f,\frac{4f-4}{3f},\frac{4}{f},\frac{2f-2}{3f})$}&$\forall$ even $f\ge 6: f\aaa\bbb\ddd,2\ccc^{\frac{f}{2}}$ & \pageref{b=2d-1}\\
			&$f=6k+4(k\ge1)$: $(f-2)\aaa\bbb\ddd,2\bbb\ccc^{\frac{f+2}{6}},2\aaa\ccc^{\frac{f-1}{3}}\ddd$&\pageref{b=2d-3}\\
			&$(f-4)\aaa\bbb\ddd,2\aaa^2\ccc^{\frac{f-4}{6}}\ddd^2,4\bbb\ccc^{\frac{f+2}{6}}$: $\lfloor \frac{k+2}{2} \rfloor$ tilings&\pageref{b=2d-4}\\
			&$(f-6)\aaa\bbb\ddd,2\aaa\ddd^3,2\aaa^2\bbb\ccc^{\frac{f-4}{6}},4\bbb\ccc^{\frac{f+2}{6}}$: $3$ tilings&\pageref{b=2d-2}\\
			\hline 
			\multirow{4}{*}{ $(\frac2f,\frac{2f-4}{3f},\frac4f,\frac{4f-2}{3f})$}&$\forall$ even $f\ge10: f\aaa\bbb\ddd,2\ccc^{\frac{f}{2}}$ &\pageref{a>1 b=2d 1} \\
			&$f=6k+2(k\ge2)$: $(f-2) \aaa\bbb\ddd,2\aaa\ccc^{\frac{f-2}{6}}\ddd,2\bbb\ccc^{\frac{f+1}{3}}$&\pageref{a>1 b=2d 2}\\
			&$(f-4)\aaa\bbb\ddd,4\aaa\ccc^{\frac{f-2}{6}}\ddd,2\bbb^2\ccc^{\frac{f+4}{6}}$: $\lfloor \frac{k+3}{2} \rfloor$ tilings&\pageref{a>1 b=2d 3}\\
			&$(f-6)\aaa\bbb\ddd,2\bbb^3\ccc,6\aaa\ccc^{\frac{f-2}{6}}\ddd$&\pageref{a>1 b=2d 4}\\
			\hline
		\end{tabular}}
	\caption{Three infinite sequences of quadrilaterals and their tilings.}\label{Tab-1.2}        
	\end{table*}

 In Table \ref{Tab-1.1} and \ref{Tab-1.2}  the angles and edge lengths are expressed in units of $\pi$,  and the last column counts all vertices and also all tilings when they are not uniquely determined by the vertices. All exact and numerical geometric data are provided in the appendix. A rational fraction, such as $\alpha=\frac{2}{9}$, means the precise value $\frac{2\pi}{9} $. A decimal expression, such as $a\approx0.3918$, means an approximate value  $0.3918\pi \le a < 0.3919\pi$. We put $\pi$ back in any trigonometric functions to avoid confusion.  
 
\begin{figure}[htp]
 	\centering
 	\includegraphics[scale=0.17]{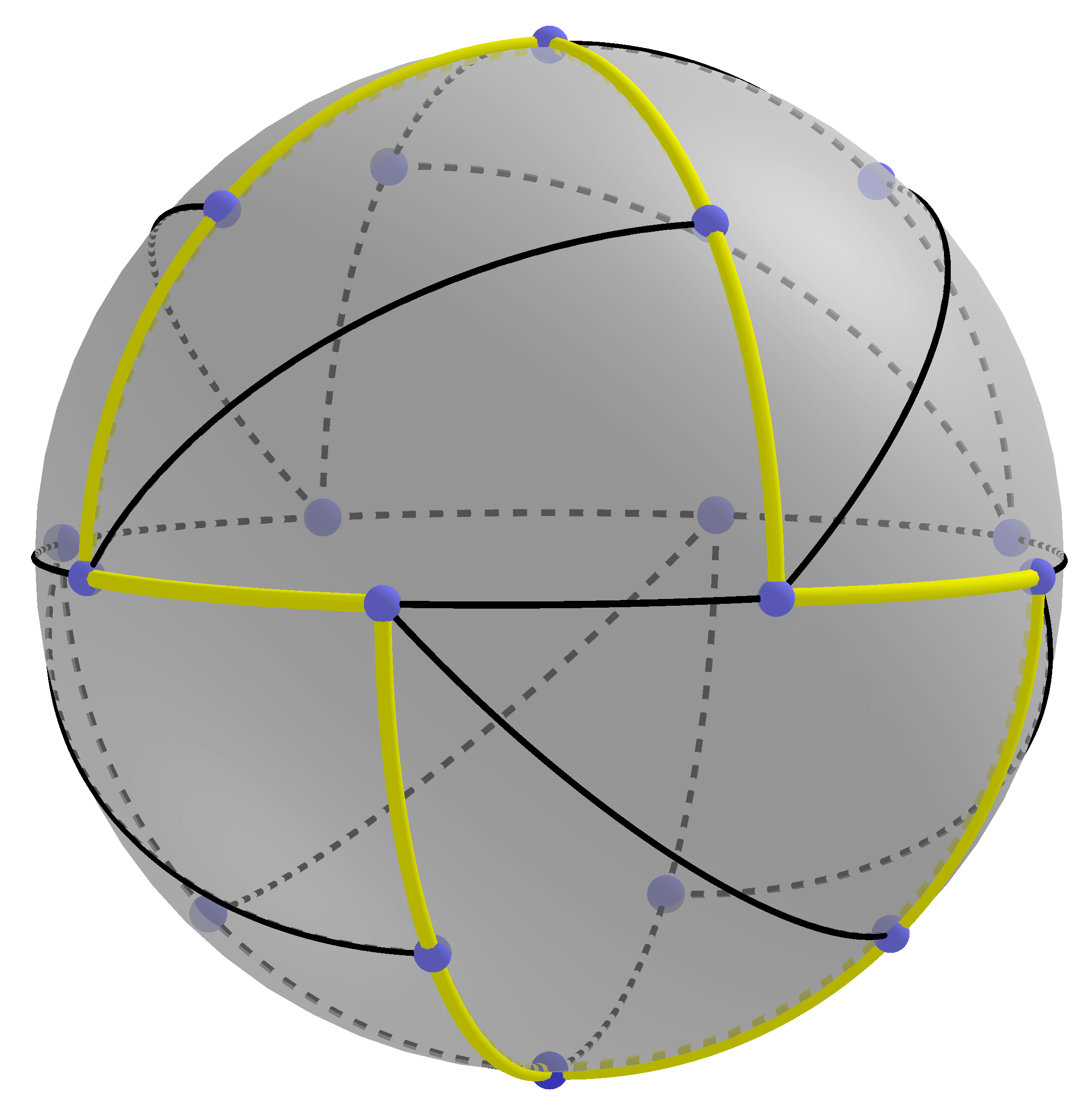}				
 	\includegraphics[scale=0.175]{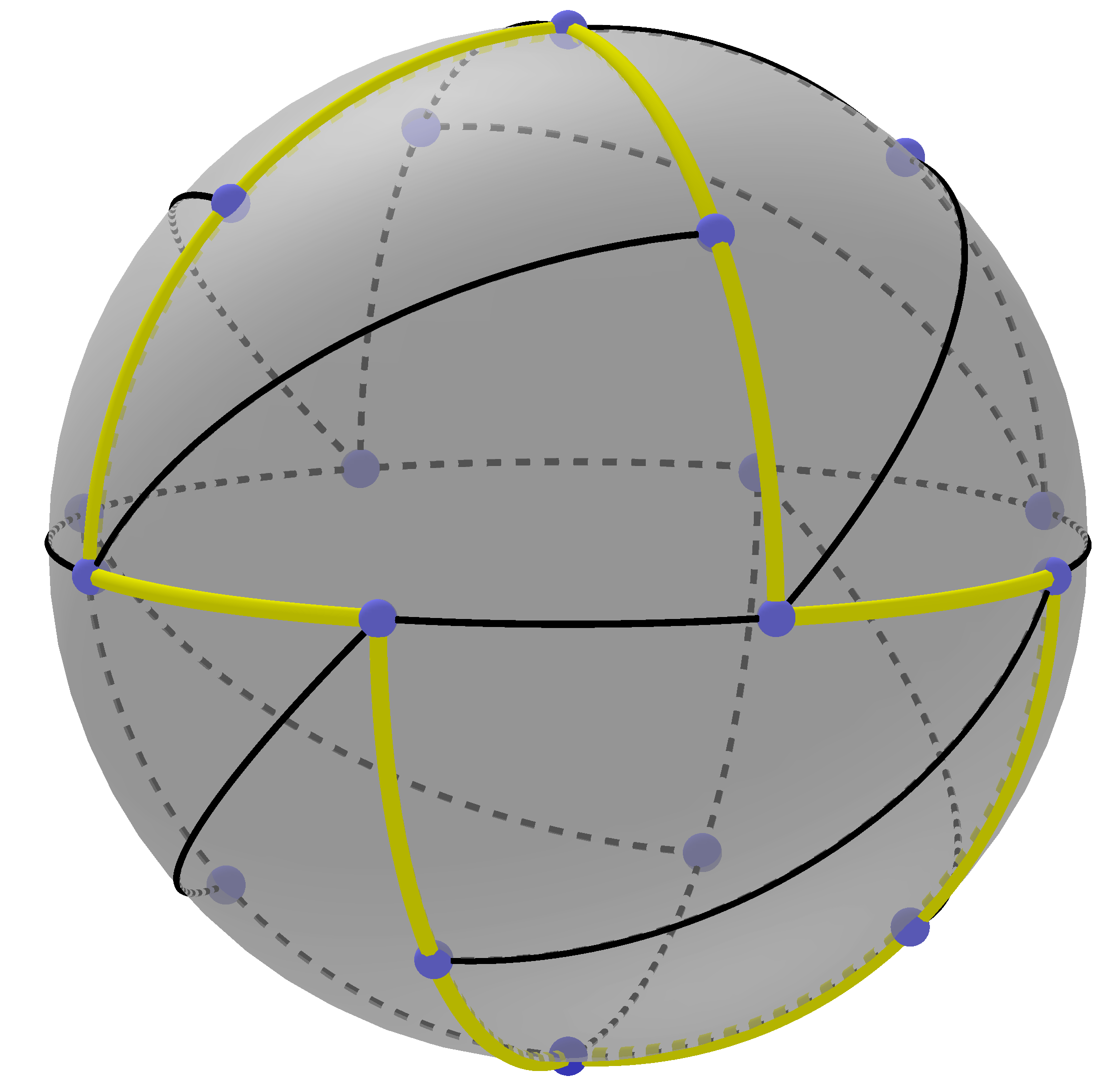}  
 	\includegraphics[scale=0.29]{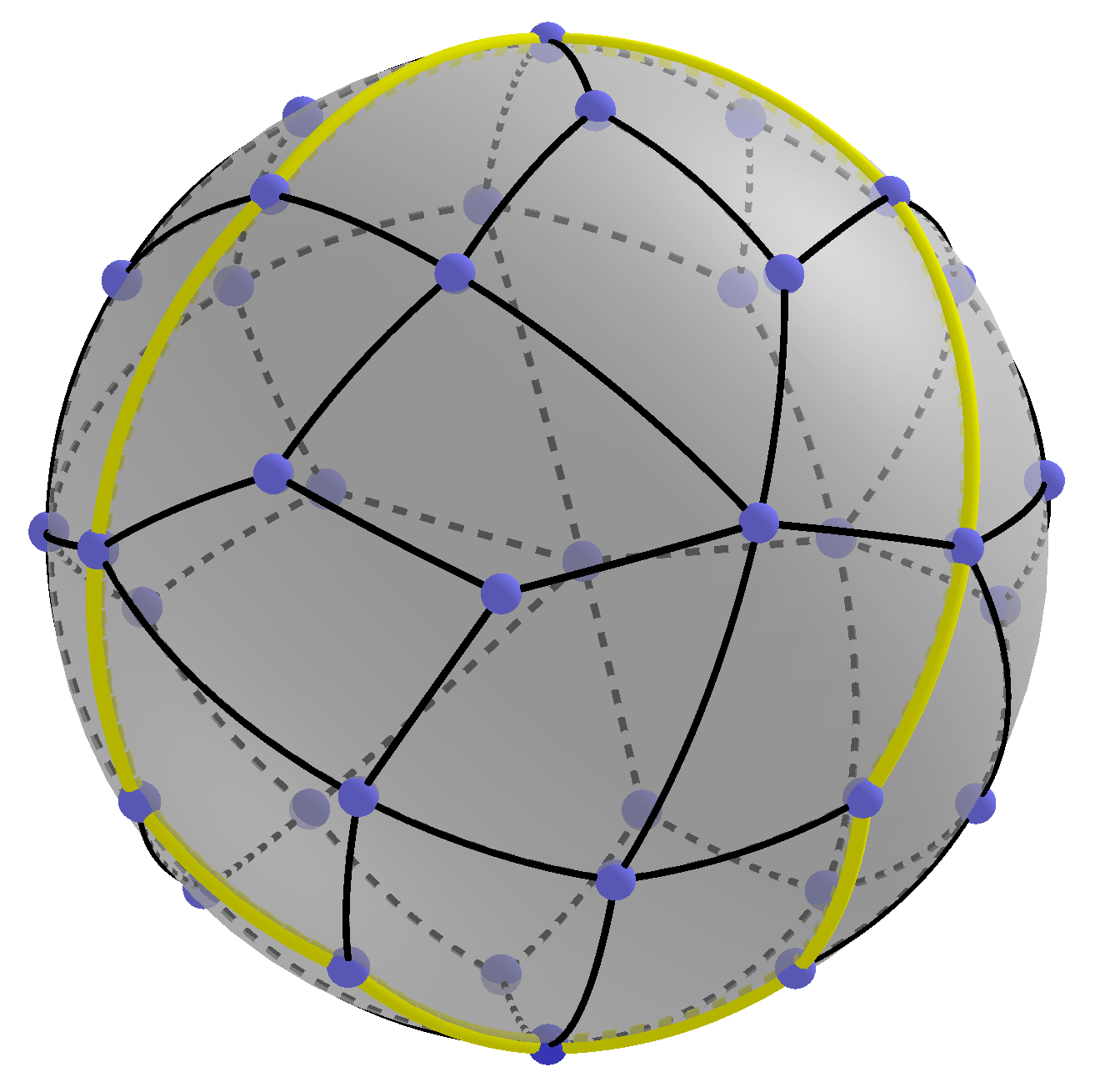}
 	\includegraphics[scale=0.263]{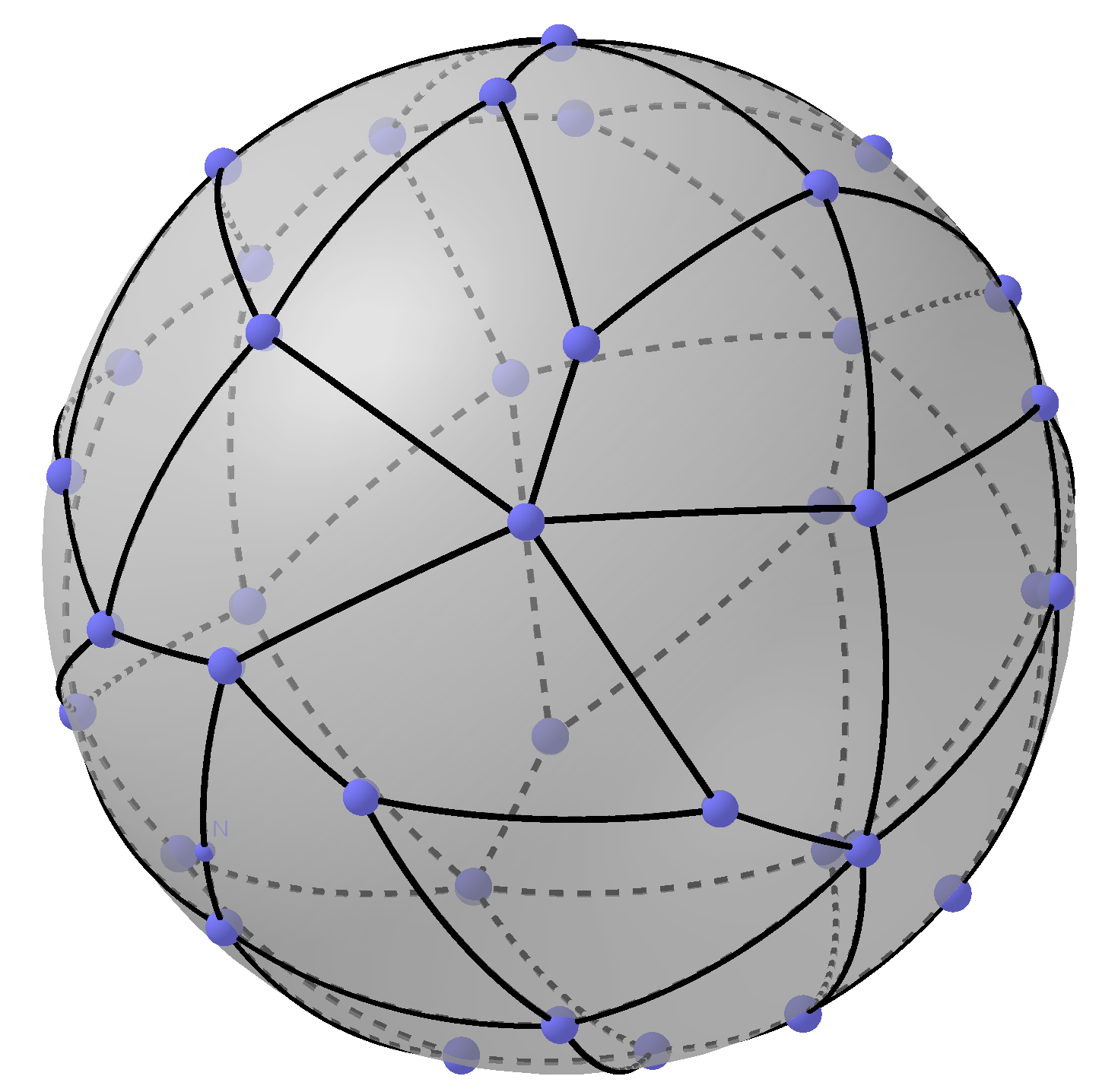}

 	\caption{Four exceptional tilings with $f=16,16,36,36$.} 
 	\label{fig 1-1-1}	
 \end{figure}

 Four exceptional tilings for the last $3$ sporadic quadrilaterals in Table \ref{Tab-1.1} ($f=16,16,36,36$)  are shown in Fig.\,\ref{fig 1-1-1}. The first $3$ tilings have repeated time zones which could be generalized combinatorially. But the quadrilaterals only exist for some particular $f$ due to geometric constraint. We remark that the last tiling ($f=36$) is the only tiling, among all edge-to-edge triangular, quadrilateral, pentagonal tilings of the sphere,  which has no apparent relation with any platonic solids or earth map tilings.

\subsection*{Modifications of special $2$-layer earth map $a^3b$-tilings}
Once all angles are fixed, there are only finitely many combinations of them summing to $2$ or form a vertex in the tiling. Then one may apply brute-force trial-and-error to find all tilings. However the  following hindsight can help us to understand most tilings in a constructive way. 

\begin{figure}[htp]
	\centering
	\includegraphics[scale=0.188]{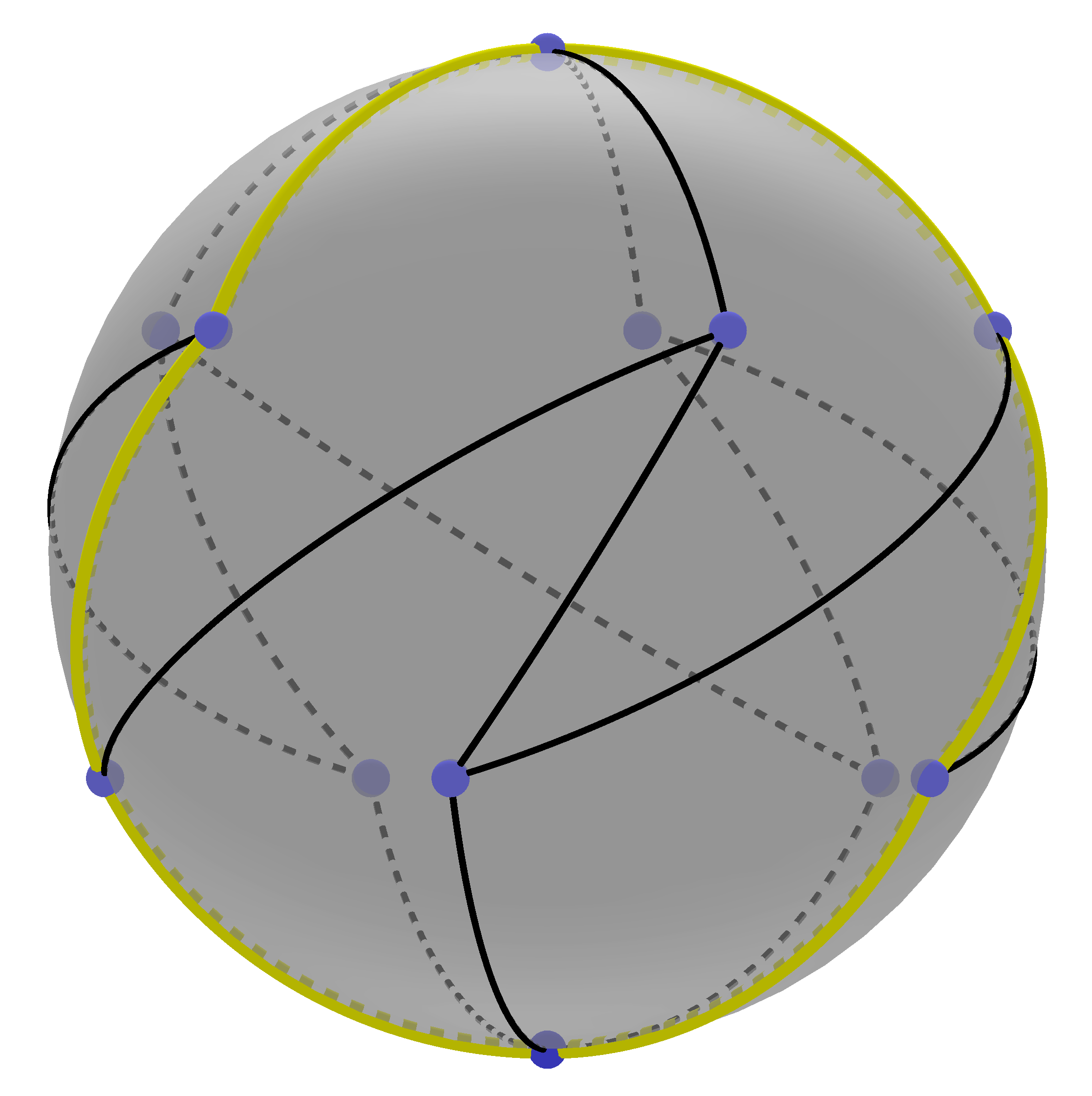}	\hspace{30pt}		
	\includegraphics[scale=0.22]{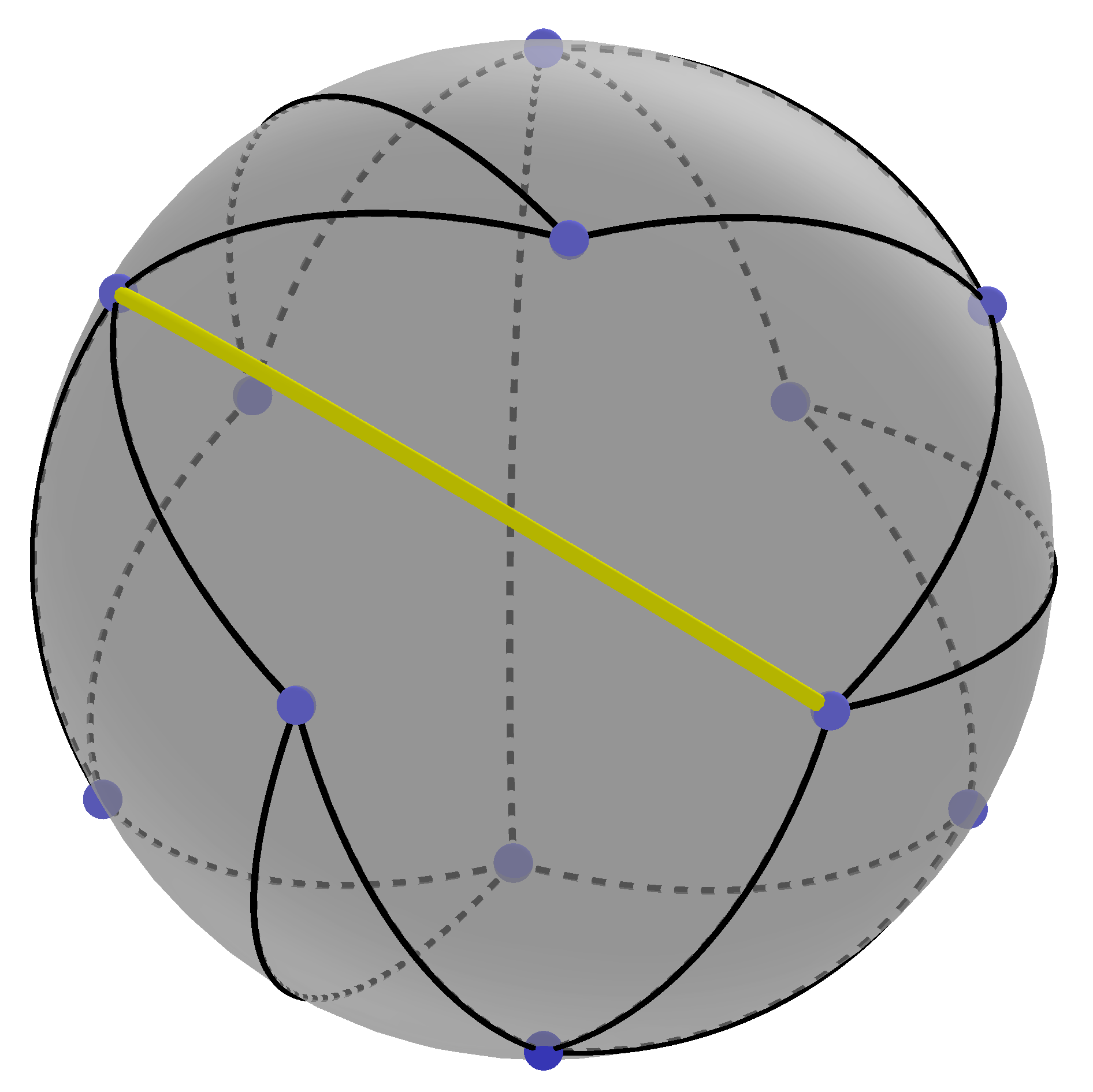}  
	
	\caption{Two very different tilings of Case $(1,6,2,3)/5$ in Table \ref{Tab-1.2}.} 
	\label{fig 1-1-3}	
\end{figure}

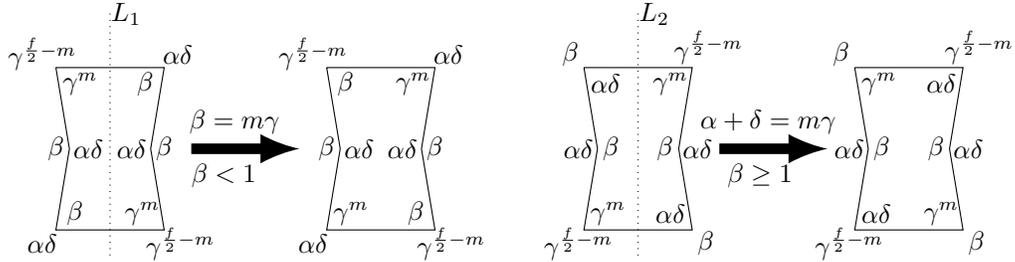
\begin{figure}[htp]
	\centering		
	\begin{tikzpicture}[>=latex,scale=0.36]
		%
		%
		%
		%
		%
		%
		%
		%
		%
		%
		%

		\begin{scope}[xshift=0 cm]
			\draw(-2,3)--(2,3)--(1.5,0)--(2,-3)--(-2,-3)--(-1.5,0)--(-2,3);
			
			\draw[line width=4pt, ->](3,0)--(7,0);
			
			\node at (4.6,1){\small $\bbb=m\ccc$};
			\node at (4.2,-1){\small $\bbb<1$};
			
			\draw[dotted] (0,5)--(0,-4);
			\node at (0.6,5){\small $L_1$};
			
			\node at (-1.1,2.4){\small $\ccc^{m}$};
			\node at (1.2,-2.4){\small $\ccc^{m}$};
			
			\node at (1.3,2.4){\small $\bbb$};
			\node at (-1.3,-2.4){\small $\bbb$};
			
			\node at (-0.8,0){\small $\aaa\ddd$};
			\node at (0.8,0){\small $\aaa\ddd$};
			
			\node at (-2,0){\small $\bbb$};
			\node at (2,0){\small $\bbb$};
			
			\node at (2.5,3.5){\small $\aaa\ddd$};
			\node at (-2.5,-3.5){\small $\aaa\ddd$};
			\node at (-2.5,3.5){\small $\ccc^{\frac{f}{2}-m}$};
			\node at (2.6,-3.5){\small $\ccc^{\frac{f}{2}-m}$};	    			
		\end{scope}	
		
		\begin{scope}[xshift=10 cm]
			\draw(-2,3)--(2,3)--(1.5,0)--(2,-3)--(-2,-3)--(-1.5,0)--(-2,3);
			
			\node at (-1.3,2.4){\small $\bbb$};
			\node at (1.3,-2.4){\small $\bbb$};
			
			\node at (1.2,2.4){\small $\ccc^{m}$};
			\node at (-1.1,-2.4){\small $\ccc^{m}$};
			
			\node at (-0.8,0){\small $\aaa\ddd$};
			\node at (0.8,0){\small $\aaa\ddd$};
			
			\node at (-2,0){\small $\bbb$};
			\node at (2,0){\small $\bbb$};
			
			\node at (2.5,3.5){\small $\aaa\ddd$};
			\node at (-2.5,-3.5){\small $\aaa\ddd$};
			\node at (-2.5,3.5){\small $\ccc^{\frac{f}{2}-m}$};
			\node at (2.6,-3.5){\small $\ccc^{\frac{f}{2}-m}$};
		\end{scope}	
		
		\begin{scope}[xshift=19.5 cm]
			\draw(-2,3)--(2,3)--(1.5,0)--(2,-3)--(-2,-3)--(-1.5,0)--(-2,3);
			
			\draw[dotted] (0,5)--(0,-4);
			\node at (0.6,5){\small $L_2$};
			
			\node at (1.2,2.4){\small $\ccc^{m}$};
			\node at (-1.1,-2.4){\small$\ccc^{m}$};
			
			\node at (-1.2,2.4){\small $\aaa\ddd$};
			\node at (1.2,-2.4){\small $\aaa\ddd$};
			
			\node at (-1,0){\small $\bbb$};
			\node at (1,0){\small $\bbb$};
			
			\node at (-2.2,0){\small $\aaa\ddd$};
			\node at (2.2,0){\small $\aaa\ddd$};
			
			\node at (-2.5,3.5){\small $\bbb$};
			\node at (2.5,-3.5){\small $\bbb$};
			\node at (2.6,3.75){\small$\ccc^{\frac{f}{2}-m}$};
			\node at (-2.2,-3.5){\small $\ccc^{\frac{f}{2}-m}$};
			\draw[line width=4pt, ->](3,0)--(7,0);
			
			\node at (4.8,1){\small $\aaa+\ddd=m\ccc$};
			\node at (4.5,-1){\small$\bbb\ge1$};
		\end{scope}	
		
		\begin{scope}[xshift=29.5 cm]
			\draw(-2,3)--(2,3)--(1.5,0)--(2,-3)--(-2,-3)--(-1.5,0)--(-2,3);
			
			\node at (1.2,2.4){\small $\aaa\ddd$};
			\node at (-1.2,-2.4){\small $\aaa\ddd$};
			
			\node at (-1.1,2.4){\small $\ccc^{m}$};
			\node at (1.2,-2.4){\small $\ccc^{m}$};
			
			\node at (-1,0){\small $\bbb$};
			\node at (1,0){\small $\bbb$};
			
			\node at (-2.2,0){\small $\aaa\ddd$};
			\node at (2.2,0){\small $\aaa\ddd$};
			
			\node at (-2.5,3.5){\small $\bbb$};
			\node at (2.5,-3.5){\small $\bbb$};
			\node at (2.6,3.75){\small $\ccc^{\frac{f}{2}-m}$};
			\node at (-2.2,-3.5){\small $\ccc^{\frac{f}{2}-m}$};
		\end{scope}

		%
		%
		%
		%
		%
		%

		
	\end{tikzpicture}
	
	\caption{Two basic flip modifications for certain $2$-layer earth map tilings.} \label{flip1}
\end{figure}

It turns out that a $2$-layer earth map $a^3b$-tiling $T(f\aaa\bbb\ddd,2\ccc^{ \frac{f}{2}})$ admits some modification
if and only if $\beta$ is an integer multiple of $\gamma$.
An authentic 3D picture for a $2$-layer earth map tiling is shown in the left of Fig.\,\ref{fig 1-1-3}. The structure of any $2$-layer earth map tiling  is shown in Fig.\,\ref{1-1}. When $\bbb=m\ccc<1$, $m$ continuous time zones ($2m$ tiles) form a dumb-bell like hexagon enclosed by $6$ $a$-edges in the first picture of Fig.\,\ref{flip1}.  Simply flip along the middle vertical line $L_1$ (or equivalently along the middle horizontal line), and one gets a new tiling of the sphere with different vertices. This is called the $1$st basic flip modification. When $\aaa+\ddd=m\ccc\le1$, we get the $2$nd basic flip modification in the right of Fig.\,\ref{flip1}.

A closer look at the inner and outer sides of this hexagon reveals that these two flips are essentially the same: $\aaa+\ddd=m\gamma$ is equivalent to $\beta=(\frac{f}{2}-m)\gamma$, and the sphere is divided by the $6$ $a$-edges into two complementary hexagons, either of which may be flipped. However it is more convenient to flip the smaller one so that we can flip several separated regions to get more tilings. So we still use both basic flips in Fig.\,\ref{flip1} but assuming afterwards that  $m\le\frac{f}{4}$. Case $(3,20,4,13)/18$ of Table \ref{Tab-1.1} and some sub-sequence of each infinite sequence of Table \ref{Tab-1.2} admit $2$ or $3$ basic flips. 

\begin{figure}[htp]
	\centering
	\includegraphics[scale=0.16]{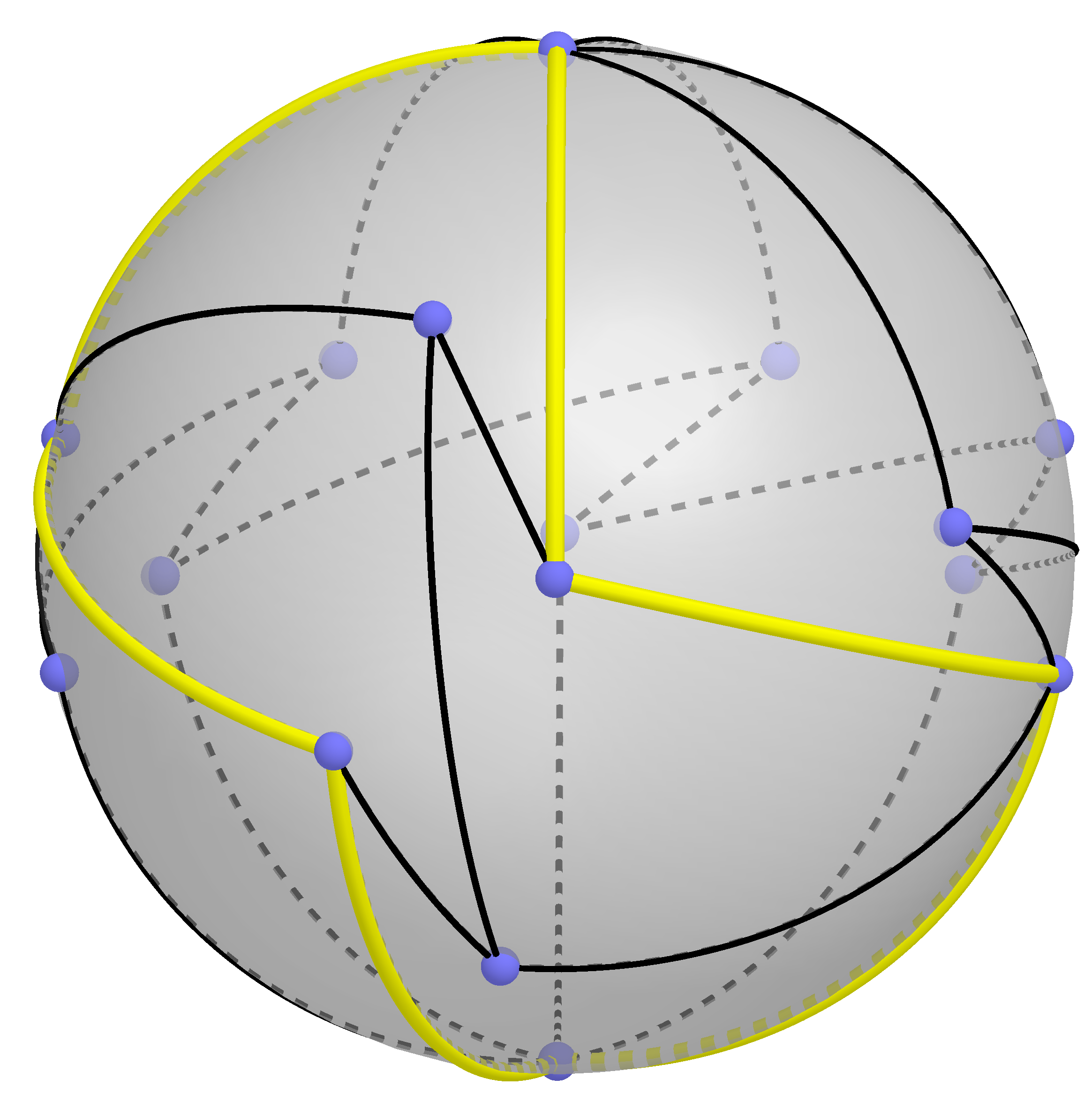}				
	\includegraphics[scale=0.155]{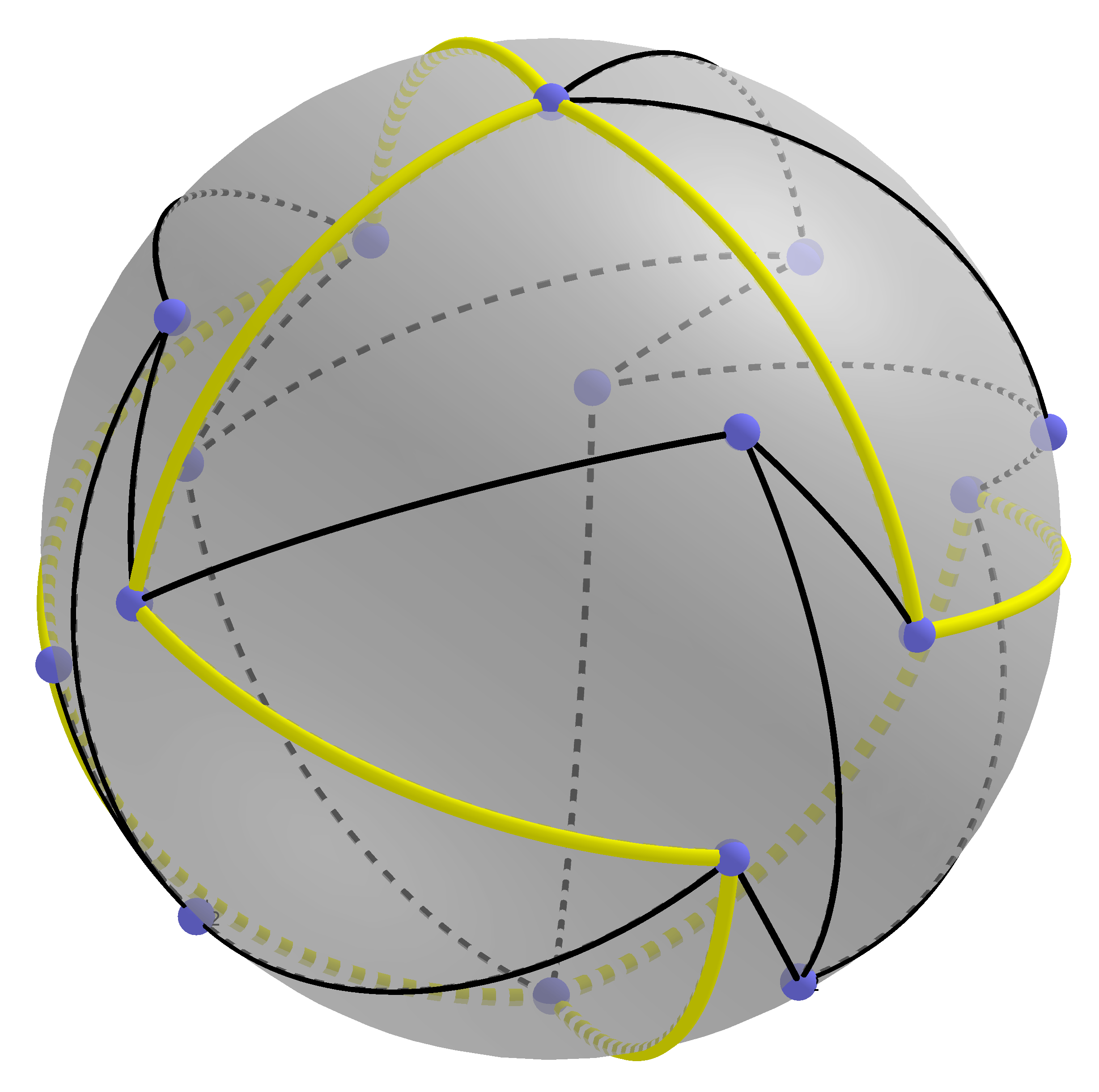}  
	\includegraphics[scale=0.16]{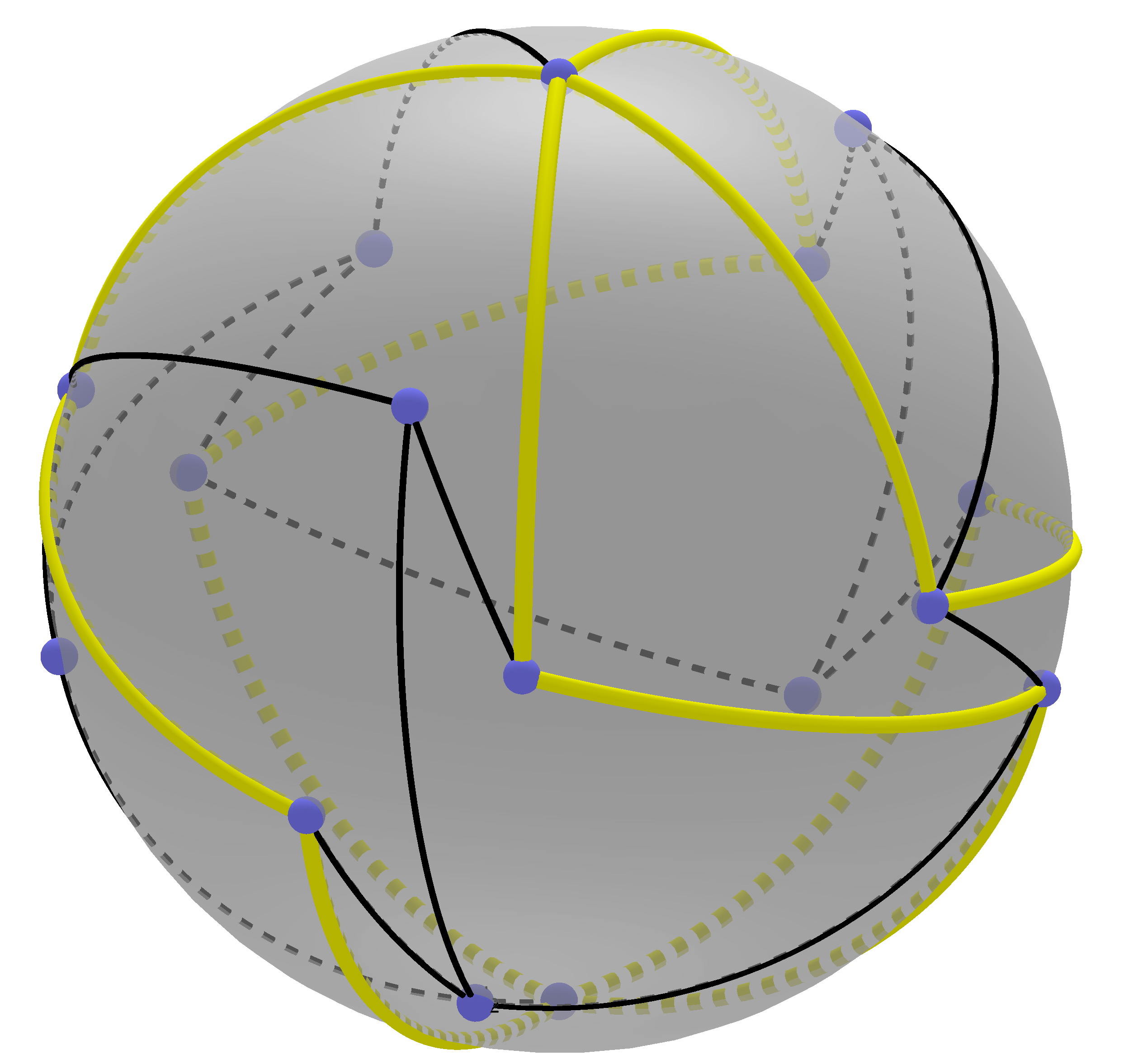}
	\includegraphics[scale=0.16]{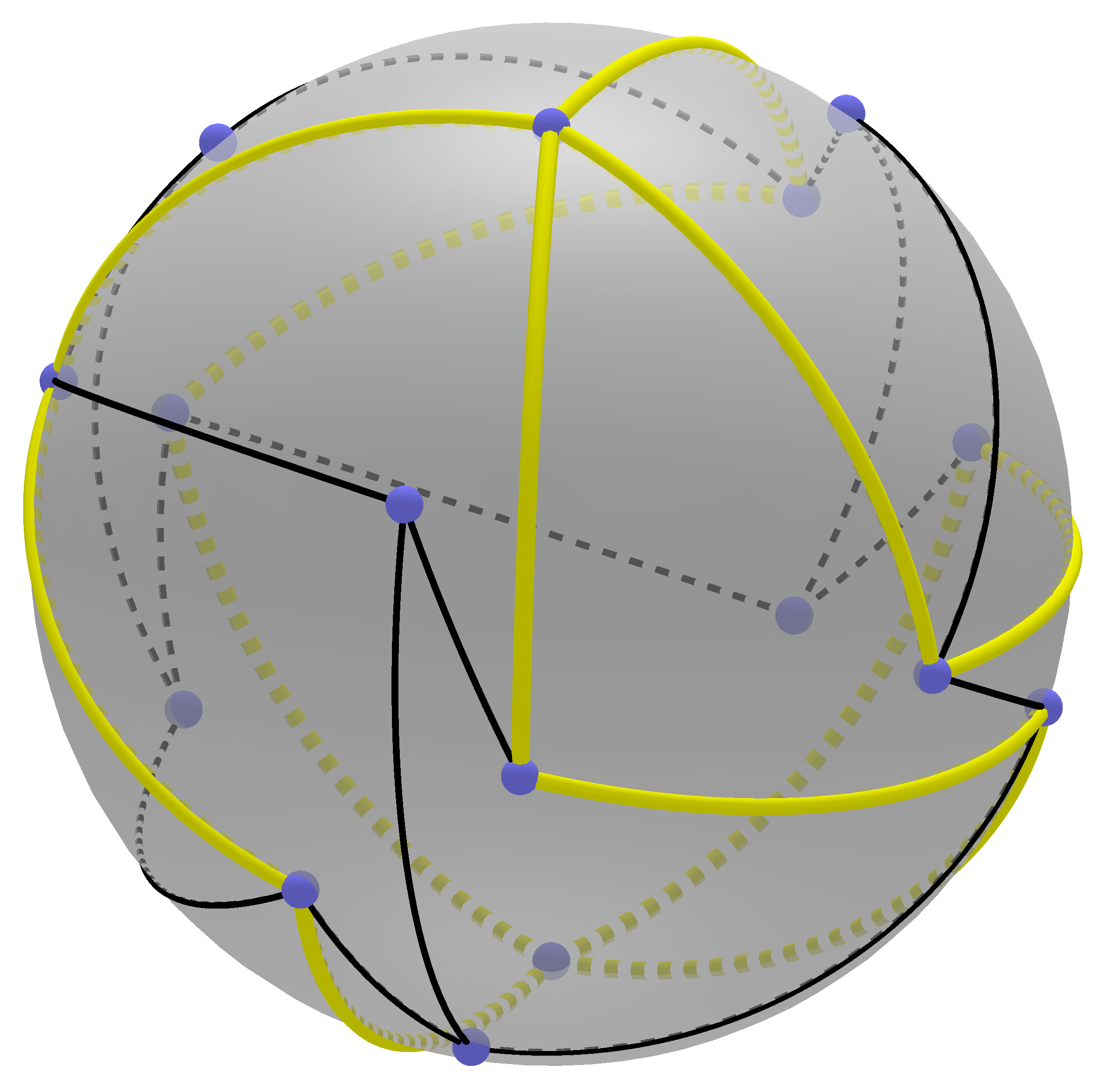}
	
	\caption{Four flips for Case $(1,4,2,9)/7$ in Table \ref{Tab-1.2}.} 
	\label{fig 1-1-2}	
\end{figure}

Fig.\,\ref{fig 1-1-2} shows $4$ different flips of the $2$-layer earth map tiling in the third case of Table \ref{Tab-1.2} with $f=14$ tiles. Flipping once, we get the $1$st picture. Flipping twice, we get the $2$nd and $3$rd pictures when the space between two flips is $0$ or $1$ time zone. Flipping three times, we get the $4$th picture. 

For Case $(\frac 2f,\frac{4f-4}{3f},\frac{4}{f},\frac{2f-2}{3f})$ with $f=6k+4\,(k\ge1)$,  there is another kind of modification giving $3$ more tilings, and we will explain it later using Fig.\,\ref{b=2d 4} and \ref{flip2}. An authentic $3$D picture for such a new tiling with $f=10$ is shown in the right of Fig.\,\ref{fig 1-1-3}.

\subsection*{Non-edge-to-edge triangular tilings}
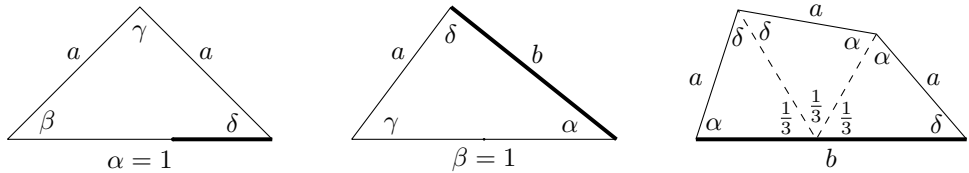
\begin{figure}[htp]
  	\centering
  	\begin{tikzpicture}[>=latex,scale=0.88]	
  		\begin{scope}[xshift=5.2 cm] 		
  			\draw (0,0)--(1.4*5/11,4.3*5/11)--(6*5/11,3.5*5/11)--(9*5/11,0);
  			\draw[line width=1.5] (0,0)--(9*5/11,0);
  			\draw[dashed] (6*5/11,3.5*5/11)--(4*5/11,0)--(1.4*5/11,4.3*5/11);
  			\node at (0.6*5/11,0.5*5/11){\small $\aaa$};\node at (5.2*5/11,3.2*5/11){\small $\aaa$};\node at (6.2*5/11,2.7*5/11){\small $\aaa$}; 
  			\node at (1.4*5/11,3.4*5/11){\small $\ddd$};\node at (2.3*5/11,3.6*5/11){\small $\ddd$};\node at (8*5/11,0.5*5/11){\small $\ddd$};
  			\node at (3*5/11,0.7*5/11){\small $\frac13$};\node at (4*5/11,1.1*5/11){\small $\frac13$}; \node at (5*5/11,0.7*5/11){\small $\frac13$};
  			\node at (0*5/11.2,2*5/11){\small $a$};  \node at (4*5/11,4.3*5/11){\small $a$}; \node at (7.9*5/11,2*5/11){\small $a$}; \node at (4.5*5/11,-0.6*5/11){\small $b$}; 			
  		\end{scope}	
  		\begin{scope}[xshift=-5.2 cm] 		
  			\draw (0,0)--(2,2)--(4,0)
  			(0,0)--(2.5,0);
  			\draw[line width=1.5] (2.5,0)--(4,0);
  			\node at (0.6,0.25){\small $\bbb$};
  			\node at (2,1.6){\small $\ccc$};
  			\node at (3.4,0.25){\small $\ddd$};
  			\fill (2.5,0) circle (0.03);
  			\node at (1,1.3){\small $a$};  	\node at (3,1.3){\small $a$};		\node at (2,-0.3){\small $\aaa=1$};  				
  		\end{scope}	 		
  		\begin{scope}[xshift=0 cm] 		 		
  		\draw (0,0)--(4,0)
  		(0,0)--(1.5,2);
  		\draw[line width=1.5] (1.5,2)--(4,0);
  		\node at (0.6,0.2){\small $\ccc$};
  		\node at (1.5,1.6){\small $\ddd$};
  		\node at (3.3,0.2){\small $\aaa$};
  		\fill (2,0) circle (0.02);
  		\node at (0.7,1.3){\small $a$};	  	\node at (2.8,1.3){\small $b$};		\node at (2,-0.3){\small $\bbb=1$};   		
  	   \end{scope}
  
  	\end{tikzpicture}
  	\caption{The degenerate and subdivision ways to get triangular tilings.}
  	\label{degradation}
  \end{figure}

When any angle of the quadrilateral is $1$, it degenerates to a triangle as shown in the first two pictures of Fig.\,\ref{degradation}. Then the first infinite sequence and two sporadic cases with $f=6,16$ produce many new examples of non-edge-to-edge triangular tilings. 

The second infinite sequence of quadrilaterals satisfy $\ccc=2\aaa,\bbb=2\ddd$ and can be subdivided into $3$ congruent triangles (observed first in \cite{ck}) as shown in the third  picture of Fig.\,\ref{degradation}, which also induce new non-edge-to-edge triangular tilings. Note that the sporadic case with $f=16$ admits such subdivision too, but inducing only some edge-to-edge triangular tiling. 
 
These are new examples, comparing to early explorations of non-edge-to-edge triangular tilings in  \cite{D1,D2,DD1,DD2,DD3}. 

  
  

\subsection*{Outline of the paper}
The classification for $a^2bc$ in \cite{lpwx1} is mainly the analysis around a special tile. However $a^3b$ is handled by a new efficient method, different from all methods developed for triangular and  pentagonal tilings. While the cost is to solve some  trigonometric Diophantine equations, the idea behind this new method is very simple: too many linearly independent vertex types in a tiling would force all angles to be rational, or the vertex types must be very limited. This paper will identify all rational $a^3b$-quadrilaterals suitable for tiling. Then the third of our series \cite{lpwx2} handles the irrational angle case in a fast way due to strong constraints on vertex types. 

This paper is organized as follows. Section \ref{basic_facts} includes general results from \cite{lpwx1} and some technical results specific to $a^3b$. Section \ref{sec-convex} looks for all possible tilings from Coolsaet's list of convex rational $a^3b$-quadrilaterals. Sections \ref{sec-concave-a} and \ref{sec-concave-b} solve some trigonometric Diophantine equations to identify all \textit{concave} rational $a^3b$-quadrilaterals suitable for tiling, and then find all of their tilings. Sections \ref{sec-degenerate-a} and \ref{sec-degenerate-b} handle two degenerate cases when the quadrilateral becomes some triangle, and thus complete the classification.

	\section{Basic Facts}
	\label{basic_facts}
	We will always express angles in $\pi$ radians for simplicity. So the sum of all angles at a vertex is $2$. We present some basic facts and techniques in this section. 

    Let $v,e,f$ be the numbers of vertices, edges, and tiles in a quadrilateral tiling. Let $v_k$ be the number of vertices of degree $k$. Euler's formula $v-e+f=2$ implies (see \cite{lpwx1})
\begin{align}
	f&=6+ \sum_{k=4}^{\infty}(k-3)v_k
	=6+v_4+2v_5+3v_6+\cdots, \label{vcountf} \\
	v_3 &=8+\sum_{k=5}^{\infty}(k-4)v_k=8+v_5+2v_6+3v_7+\cdots. \label{vcountv}
\end{align}
    So $f\ge 6$ and $v_3 \ge 8$. 

\begin{lemma} [{\cite[Lemma 2]{lpwx1}}] \label{anglesum} 
	If all tiles in a tiling of the sphere by $f$ quadrilaterals have the same four angles $\aaa,\bbb,\ccc,\ddd$, then 
	\[
	\aaa+\bbb+\ccc+\ddd = 2+\frac{4}{f} , 
	\]
	ranging in $(2,\frac83]$. In particular no vertex contains all four angles.
\end{lemma}

\begin{lemma} [{\cite[Lemma 3]{wy2}}] \label{geometry1}
	If the  quadrilateral in Fig.\,\ref{quad} is simple, then $\bbb<\ccc$ is equivalent to $\aaa>\ddd$.
\end{lemma}

\begin{lemma}\label{geometry3}
	If the  quadrilateral in Fig.\,\ref{quad} is simple, then $\bbb=\ddd$ if and only if $\aaa=1$. Furthermore, if it is convex with all angles $<1$, then $\bbb>\ddd$ is equivalent to $\aaa<\ccc$, and  $\bbb<\ddd$ is equivalent to $\aaa>\ccc$.    
\end{lemma}

\begin{proof}
	
	If $\aaa=1$, we get an isosceles triangle in the first picture of Fig.\,\ref{g3proof}, thus $\bbb=\ddd$. If $\bbb=\ddd$ and $\aaa\neq1$, then $\angle CBD = \angle BDC$ implies $\angle ABD = \angle ADB$. So we get $a=b$, a contradiction. 
	When the quadrilateral is convex with all angles $<1$, the line $AC$ in the second of Fig.\,\ref{g3proof} is inside the quadrilateral, and divides $\aaa$ and $\ccc$ as $\ccc=\theta+\ccc'$ and $\aaa=\theta+\aaa'$. Then 
	\[
	\aaa<\ccc
	\iff 
	\aaa'<\ccc'
	\iff
	a<b.
	\]
	By the same reason, we have $\bbb>\delta \iff a<b$. Therefore, $\bbb>\delta$ is equivalent to $\aaa<\ccc$. Similarly $\bbb<\ddd$ is equivalent to $\aaa>\ccc$.
\end{proof}

\begin{figure}[htp]
	\centering
	\begin{tikzpicture}[>=latex,scale=0.8]		
		\begin{scope}[xshift=-6 cm,scale=1] 		
			\draw (0,0)--(2,2)--(4,0)
			(0,0)--(2.5,0);
			\draw[line width=1.5] (2.5,0)--(4,0);
			\node at (0.6,0.25){\small $\bbb$};
			\node at (2,1.6){\small $\ccc$};
			\node at (3.4,0.25){\small $\ddd$};
			\fill (2.5,0) circle (0.03	);
			\node at (1,1.3){\small $a$};  	\node at (3,1.3){\small $a$};		\node at (2,-0.3){\small $a+b$};  				
		\end{scope}
		
		\foreach \a in {0,1}
		{
			\begin{scope}[xshift=6*\a cm]
				
				\draw[dashed]
			    (0,0) -- (2.65,1.53);
				
				\node at (0.5,2) {\small $B$};
				\node at (2.7,1.75) {\small $C$};
				\node at (-0.2,0) {\small $A$};
				
			\end{scope}
		}
		
		\draw
		(0,0) -- (70:2) -- ++(-10:2) -- ++(-50:2);
		
		\draw[dashed]
		(0,0) -- (2.65,1.53);
		
		\draw[line width=1.5]
		(0,0) -- (3.95,0);	
		
		\node at (0.8,1.5) {\small $\bbb$};
		\node at (2.05,1.4) {\small $\theta$};
		\node at (0.35,0.45) {\small $\theta$};
		\node at (2.55,1.2) {\small $\ccc'$};
		\node at (0.8,0.25) {\small $\aaa'$};
		\node at (3.4,0.25) {\small $\ddd$};
		
		\node at (4.1,0) {\small $D$};
		
		\begin{scope}[xshift=6 cm]
			
			\draw
			(0,0) -- (70:2) -- ++(-10:2); 
			
			\draw
			(2.65,1.53) -- ++(12:1.25);
			
			\draw[line width=1.5]
			(0,0) to[out=0, in=-110] (3.9,1.8);	
			
			\node at (4,2) {\small $D$};
			
		\end{scope}

	\end{tikzpicture}
	\caption{Proof of Lemmas \ref{geometry3} and \ref{geometry4}.}
	\label{g3proof}
\end{figure}

\begin{lemma}\label{geometry4}
	If the  quadrilateral in Fig.\,\ref{quad} is simple, and $\ddd\le 1$, then $2\aaa+\bbb>1$ and $\bbb+2\ccc>1$.
\end{lemma}

\begin{proof}
	If all angles are $< 1$, then the quadrilateral is convex and the line $AC$ is inside the quadrilateral  in the second picture of Fig.\,\ref{g3proof}.  Thus $\theta<\aaa,\ccc$. This implies  $2\aaa+\bbb>2\theta+\bbb>1$ and $\bbb+2\ccc>\bbb+2\theta>1$.
	
	If $\bbb\ge 1$ or both $\aaa,\ccc\ge 1$, then both inequalities certainly hold. If $\ddd=1$, then $\aaa=\ccc$ by Lemma \ref{geometry3}, and $2\aaa+\bbb>1$ as the angle sum of a triangle. So we only need to consider the following two cases: 
	\begin{enumerate}
		\item $\aaa,\bbb,\ddd< 1$ and $\ccc\ge 1$. 
		\item $\bbb,\ccc,\ddd< 1$ and $\aaa\ge 1$.
	\end{enumerate}
	Case $1$ is shown in the third picture of Fig.\,\ref{g3proof}, and it suffices to show $2\aaa+\bbb>1$. By $\aaa,\ddd<1$ and $AB=CD=a<1$, both $B$ and $C$ lie in the interior of the same hemisphere bounded by the great circle $\bigcirc AD$. By $\bbb< 1\le\ccc$ and $\aaa>\ddd$ (Lemma \ref{geometry1}), the line $AC$ is inside the quadrilateral. Then $\aaa\ge \theta$ as in the second picture of Fig.\,\ref{g3proof}, and  $2\aaa+\bbb>2\theta+\bbb>1$. Case $2$ can be proved similarly.
\end{proof}

\begin{lemma} [{Parity Lemma, \cite[Lemma 10]{wy2}}] \label{beven}
	In an $a^3b$-tiling, the total number of $ab$-angles $\aaa$ and $\ddd$ at any vertex is even. 
\end{lemma}

In a tiling of the sphere by $f$ congruent tiles, each angle of the tile appears $f$ times in total. If one vertex has more $\aaa$ than $\bbb$, there must exist another vertex with more $\bbb$ than $\aaa$.
Such global counting induces many interesting and useful results. 

\begin{lemma}[{Balance Lemma, \cite[Lemma 11]{wy2}}] \label{balance}
	 If either $\aaa^2\cdots$ or $\delta^2\cdots$ is not a vertex, then any vertex either has no $\aaa,\delta$, or is of the form $\aaa\delta\cdots$ with no more $\aaa,\delta$ in the remainder.
\end{lemma}

The very useful tool \textit{adjacent angle deduction} (abbreviated as AAD) has been introduced in \cite[Section 2.5]{wy1}. The following is \cite[Lemma 10]{wy1}.

\begin{lemma}\label{aadlemma}
	Suppose $\lambda$ and $\mu$ are the two angles adjacent to $\theta$ in a quadrilateral.
	\begin{itemize}
		\item If $\lambda\thin\lambda\cdots$ is not a vertex, then $\theta^n$ has the unique AAD $\thin^{\lambda}\theta^{\mu}\thin^{\lambda}\theta^{\mu}\thin^{\lambda}\theta^{\mu}\thin\cdots$.
		\item If $n$ is odd, then we have the AAD $\thin^{\lambda}\theta^{\mu}\thin^{\lambda}\theta^{\mu}\thin$ at $\theta^n$.
	\end{itemize}
\end{lemma}

\begin{lemma}\label{lemma-1}
	There is no tiling of the sphere by congruent quadrilaterals with two angles $\ge1$. 
\end{lemma}	
\begin{proof}	
	If any two angles, say $\aaa,\bbb$, are greater than or equal to $1$, then  $\aaa\cdots=\aaa\ccc^x\ddd^y(x+y\ge 2)$, $\bbb\cdots=\bbb\ccc^p\ddd^q(p+q\ge 2)$. Given that $\#\aaa=\#\bbb=f$, we deduce that $\#\ccc+\#\ddd\ge 4f$, which contradicts $\#\ccc+\#\ddd=2f$. 	
\end{proof}	

\begin{proposition}\label{symmetric}	
	There is no tiling of the sphere by congruent symmetric $a^3b$-quadrilaterals ($\aaa=\ddd$ and $\bbb=\ccc$).
\end{proposition}
\begin{proof}
	The convex case with all angles $<1$ has been proved by Akama and Cleemput in \cite{ac}. If any angle is $\ge1$, we get two angles $\ge1$ by symmetry, then Lemma \ref{lemma-1} applies. 
\end{proof}

\begin{lemma}\label{proposition-7}
	Assume $\ccc^{\frac{f}{2}}$ is a vertex in an $a^3b$-tiling. If $\bbb^2\cdots$ or $\ddd^2\cdots$ is not a vertex, and $\bbb\ddd\cdots=\aaa\bbb\ddd$, then the tiling must be a $2$-layer earth map tiling $T(f\aaa\bbb\ddd,2\ccc^{\frac{f}{2}})$ in Fig.\,\ref{1-1}. In particular if all $\bbb$-vertices are $\aaa\bbb\ddd$, then the tiling must be a $2$-layer earth map tiling. 
\end{lemma}

\begin{proof}
By Lemma \ref{aadlemma}, when $\bbb^2\cdots$ or $\ddd^2\cdots$  is not a vertex, we have the unique AAD $\ccc^{\frac{f}{2}}=\thin^{\bbb}\ccc^{\ddd}\thin^{\bbb}\ccc^{\ddd}\thin\cdots$. In Fig.\,\ref{1-1}, $\ccc_1\ccc_2\ccc_3\cdots$ determines $T_1,T_2,T_3$. Then  $\bbb_2\ddd_1\cdots=\aaa_4\bbb_2\ddd_1$ determines $T_4$;  $\bbb_3\ddd_2\cdots=\aaa_5\bbb_3\ddd_2,\aaa_2\bbb_4\ddd_5\cdots=\aaa_2\bbb_4\ddd_5$ determines $T_5$. The argument started at $\aaa_4\bbb_2\ddd_1$ can be repeated at $\aaa_5\bbb_3\ddd_2$. More repetitions give the 
	unique tiling of $f$ tiles with $2\ccc^{\frac{f}{2}}$ and $f\aaa\bbb\ddd$.
	\begin{figure}[htp]
		\centering
		\begin{tikzpicture}[>=latex,scale=0.6] 
			\foreach \a in {0,1,2}
			{
				\begin{scope}[xshift=2*\a cm] 
					\draw (0,0)--(0,-2)
					(2,0)--(2,-2)--(3.5,-3)--(3.5,-5)
					(1.5,-3)--(1.5,-5)
					(0,-2)--(1.5,-3);
					\draw[line width=1.5] (1.5,-3)--(2,-2);
					\node at (1.4,-2.55){\small $\aaa$};
					\node at (1.8,-3.2){\small $\ddd$};
					\node at (3.1,-3.2){\small $\bbb$};
					\node at (1,0){\small $\ccc$};
					\node at (2.5,-5){\small $\ccc$};
					
					\node at (2.1,-2.45){\small $\aaa$};
					\node at (1.7,-1.8){\small $\ddd$};
					\node at (0.3,-1.7){\small $\bbb$};
				\end{scope}
			}
			
			\fill (8,-2) circle (0.05); \fill (7.7,-2) circle (0.05);
			\fill (8.3,-2) circle (0.05);
			
			\node[draw,shape=circle, inner sep=0.5] at (1,-1) {\small $1$};
			\node[draw,shape=circle, inner sep=0.5] at (3,-1) {\small $2$};
			\node[draw,shape=circle, inner sep=0.5] at (5,-1) {\small $3$};
			\node[draw,shape=circle, inner sep=0.5] at (2.5,-4) {\small $4$};
			\node[draw,shape=circle, inner sep=0.5] at (4.5,-4) {\small $5$};
			\node[draw,shape=circle, inner sep=0.5] at (6.5,-4) {\small $6$};
			
		\end{tikzpicture}
		\caption{A $2$-layer earth map tiling $T(f\aaa\bbb\ddd,2\ccc^{\frac{f}{2}})$.} \label{1-1}
	\end{figure}
\end{proof}

\begin{lemma}\label{proposition-5'}
	In an $a^3b$-tiling, if $\aaa\ge1$, then either $\aaa\bbb\ddd$ or $\aaa\ccc\ddd$ is a vertex, and the only other possible vertex with $\aaa$ or $\ddd$ must be $\aaa\ccc^l\ddd$ or $\aaa\bbb^l\ddd$ respectively for some $l\ge2$. 
\end{lemma}
\begin{proof}
$\aaa\ge1$ implies $\aaa^2\cdots$ is not a vertex. Then Balance Lemma \ref{balance} and Lemma \ref{anglesum} imply that any vertex with $\aaa$ or $\ddd$ must be of two types $\aaa\bbb^l\ddd$ or $\aaa\ccc^m\ddd$. If there exists only one type, say $\aaa\bbb^l\ddd$, then $l=1$ by Balance Lemma \ref{balance}. If there exist both types with $l,m\ge2$, then the only solution  satisfying Balance Lemma \ref{balance} is: $\{\frac f2 \aaa\bbb^2\ddd,\frac f2 \aaa\ccc^2\ddd\}$.  This implies $\bbb=\ccc$, contradicting Proposition \ref{symmetric}. Therefore, one of $l,m$ must be $1$, and the other must be $\ge2$ since $\bbb\neq\ccc$. 
\end{proof}

\begin{lemma}\label{geometry5}
  In an $a^3b$-tiling, the $a$-edge and two diagonals are always $<1$. If both $\bbb,\ccc<1$, then $b<1$. 
\end{lemma}
\begin{proof}
	By Lemma \ref{lemma-1}, there are just three types of simple quadrilaterals suitable for tiling: convex with all angles $<1$, $\aaa\ge1$, or $\bbb\ge1$, as shown in Fig.\,\ref{quadrilateral}. It is clear that $a<1$, otherwise $BC$ and $CD$ would intersect at the antipodal of $C$, contradicting the simpleness. 
	
	For the first two types in Fig.\,\ref{quadrilateral}, both $\bbb<1$ and $\ccc<1$, then $a<1$ implies that both $A$ and $D$ lie in the interior of the same hemisphere bounded by the great circle $\bigcirc BC$. Therefore, two diagonals and $b$-edge are all $<1$. 
	
	For the last type in Fig.\,\ref{quadrilateral}, both $\aaa<1$ and $\ddd<1$, then $a<1$ implies that both $B$ and $C$ lie in the interior of the same hemisphere bounded by the great circle $\bigcirc AD$.  Therefore, both diagonals are $<1$. 	
\end{proof}

\begin{lemma}\label{geometry6}
	For $a^3b$-quadrilaterals, the following equations hold: 
	\begin{equation}
		\begin{split}
			\cos b=& \cos^3 a(1-\cos \bbb)(1-\cos \ccc)-\cos^2 a\sin \bbb\sin \ccc+\\ &\cos a(\cos \bbb+\cos \ccc-\cos \bbb\cos \ccc) +\sin\bbb\sin\ccc;  \label{4}
		\end{split}
	\end{equation} 
	\begin{equation}
		\cos a=\frac{\sin \aaa+\cos\ddd\sin\ccc}{2\sin\ddd\sin^2\frac{\ccc}{2}}=\frac{\sin \ddd+\cos\aaa\sin\bbb}{2\sin\aaa\sin^2\frac{\bbb}{2}} \quad (\aaa,\ddd\neq1);   \label{4-2}
	\end{equation}
	\begin{align}
		\sin(\aaa-\frac{\ccc}{2})\sin\frac{\bbb}{2}=\sin\frac{\ccc}{2}\sin(\ddd-\frac{\bbb}{2}), \label{4-7}
	\end{align}
	\begin{align}
		\text{or}\,\,\,	\sin(\aaa+\frac{\ccc}{2})\sin\frac{\bbb}{2}=-\sin\frac{\ccc}{2}\sin(\ddd+\frac{\bbb}{2}).  \label{4-8}
	\end{align} 
\end{lemma}
\begin{proof}	
	The equation \eqref{4} always holds by the extended cosine law in \cite[Lemma 11]{wy1}. By Lemma \ref{lemma-1}, there are just three types of simple quadrilaterals suitable for tiling: convex with all angles $<1$, $\aaa\ge1$, or $\bbb\ge1$, as shown in Fig.\,\ref{quadrilateral}. When $\aaa,\bbb,\ccc,\ddd<1$, Lemma \ref{geometry5} implies that all edges and diagonals are $<1$, and the equations \eqref{4-2},\eqref{4-7} were proved in \cite[Theorem 2.1]{ck} while the equation \eqref{4-8} was dismissed. 
	
	The same proof works for the other two concave cases in Fig.\,\ref{quadrilateral}. If $\aaa>1$, the sine law $\frac{\sin(2-\aaa)}{\sin y}=\frac{\sin(\psi-\ddd)}{\sin a}$ is equivalent to $\frac{\sin\aaa}{\sin y}=\frac{\sin(\ddd-\psi)}{\sin a}$. If $\bbb>1$, the sine law $\frac{\sin(2-\bbb)}{\sin x}=\frac{\sin(\phi)}{\sin a}$ is equivalent to $\frac{\sin\bbb}{\sin x}=\frac{\sin(-\phi)}{\sin a}$. Then every step to derive the equation \eqref{4-2} is exactly the same as the convex case, and one of the equations \eqref{4-7} or \eqref{4-8} must hold. But the equation \eqref{4-8} can no longer be dismissed easily as the convex case. 
	\begin{figure}[htp]
		\centering
		\begin{tikzpicture}[>=latex,scale=0.45] 
			\begin{scope} 
				\draw (-1.1*1.4,2.4*1.4-0.5)--(0,0-0.5)--(2.5*1.4,0-0.5)--(4*1.4,1.9*1.4-0.5);
				\draw[line width=1.5] (-1.1*1.4,2.4*1.4-0.5)--(4*1.4,1.9*1.4-0.5);
				\draw[dotted] (0,0-0.5)--(4*1.4,1.9*1.4-0.5)
				(-1.1*1.4,2.4*1.4-0.5)--(2.5*1.4,0-0.5);
				
				\node at (-1.35*1.4,2.5*1.4-0.5){\small $A$};\node at (-0.05*1.4,-0.3*1.4-0.5){\small $B$};\node at (2.5*1.4,-0.3*1.4-0.5){\small $C$};
				\node at (4.25*1.4,2*1.4-0.5){\small $D$};
				\node at (0.85*1.4,1.4*1.4-0.5){\small $x$};\node at (2.2*1.4,1.3*1.4-0.5){\small $y$};
				\node at (-0.6*1.4,1.8*1.4-0.5){\small $\phi$};  \node at (1.8*1.4,0.2*1.4-0.5){\small $\phi$};  \node at (0.8*1.4,0.2*1.4-0.5){\small $\psi$};  \node at (3.3*1.4,1.3*1.4-0.5){\small $\psi$}; 
				
			\end{scope}
			\begin{scope}	[xshift=12 cm] 
				\draw (5*0.9,-4*0.9+3)--(0,0+3)--(-5*0.9,-4*0.9+3)--(1*0.9,-3*0.9+3);

				\draw[line width=1.5] (5*0.9,-4*0.9+3)--(1*0.9,-3*0.9+3);
				
				\draw[dotted] (0,0+3)--(1*0.9,-3*0.9+3)
				(5*0.9,-4*0.9+3)--(-5*0.9,-4*0.9+3);
				
				\node at (0,0.5*0.9+3){\small $C$};  \node at (0.8*0.9,-3.5*0.9+3){\small $A$};  \node at (-5.05*0.9,-4.45*0.9+3){\small $B$};
				\node at (4.9*0.9,-4.45*0.9+3){\small $D$};   \node at (0,-1.5*0.9+3){\small $x$}; \node at (0,-4.5*0.9+3){\small $y$};

				\node at (-3.3*0.9,-3.2*0.9+3){\small $\bbb$};  \node at (3.3*0.9,-3.2*0.9+3){\small $\ddd$};
				
				\node at (-0.2*0.9,-0.5*0.9+3){\small $\phi$};  \node at (0.5*0.9,-2.7*0.9+3){\small $\phi$};
				
				\node at (-4.1*0.9,-3.7*0.9+3){\small $\psi$}; \node at (4.1*0.9,-3.7*0.9+3){\small $\psi$};

				\node at (-2.5*0.9,-1.5*0.9+3){\small $a$}; \node at (2.3*0.9,-2.9*0.9+3){\small $b$}; \node at (2.5*0.9,-1.5*0.9+3){\small $a$};\node at (-1.8*0.9,-3.1*0.9+3){\small $a$};

			\end{scope}
			\begin{scope}	[xshift=23 cm] 
				\draw (-4*0.65,5*0.65)--(-6*0.65,-1*0.65)--(0,0)--(6*0.65,-1*0.65);

				\draw[line width=1.5] (-4*0.65,5*0.65)--(6*0.65,-1*0.65);
				
				\draw[dotted] (-6*0.65,-1*0.65)--(6*0.65,-1*0.65)
				(0,0)--(-4*0.65,5*0.65);
				
				\node at (-6*0.65,-1.6*0.65){\small $C$};  \node at (5.8*0.65,-1.6*0.65){\small $A$};  \node at (0,-0.55*0.65){\small $B$};
				\node at (-4.4*0.65,5.5*0.65){\small $D$};   \node at (0,-1.7*0.65){\small $x$};    \node at (-2.5*0.65,2.2*0.65){\small $y$};

				\node at (-0.9*0.65,0.35*0.65){\small $\psi$}; \node at (4.2*0.65,-0.35*0.65){\small $\aaa$};
				\node at (-5.3*0.65,-0.4*0.65){\small $\ccc$}; \node at (-3.8*0.65,3.9*0.65){\small $\psi$};
				
				\node at (3*0.65,-0.8*0.65){\small $\phi$};  \node at (-3*0.65,-0.8*0.65){\small $\phi$};

				\node at (-5.5*0.65,2*0.65){\small $a$}; \node at (2.1*0.65,2*0.65){\small $b$}; \node at (-2*0.65,0.1*0.65){\small $a$};\node at (2*0.65,0.1*0.65){\small $a$};
				
			\end{scope}
		\end{tikzpicture}
		\caption{Proof of Lemmas \ref{geometry5} and \ref{geometry6}.} \label{quadrilateral}
	\end{figure}
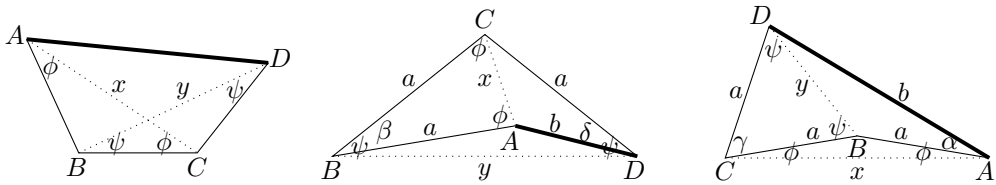	
\end{proof}
	It is amazing that all rational solutions (rational multiples of $\pi$) to \eqref{4-7} or \eqref{4-8} can be found via the algebra of cyclotomic fields in Conway-Jones \cite{cj}, as Coolsaet \cite{ck} did for convex $a^3b$-quadrilaterals using Myerson's Theorem \cite{mg}. We summarize the algorithm as the following easy-to-use proposition. 
	\begin{proposition}\label{proposition-6}
		All solutions of $\sin x_1 \sin x_2=\sin x_3 \sin x_4$ with rational angles $0\le x_1,x_2,x_3,x_4\le\frac12$ fall into the following four cases:
		\begin{enumerate}	
		\item $x_1  x_2=x_3 x_4=0$;
		
		\item $\{x_1, x_2\}=\{ x_3, x_4\}$;
		
		\item $\{x_1, x_2\}=\{ \frac{1}{6},\theta\}$ and $\{x_3,x_4\}=\{\frac{\theta}{2},\frac{1}{2}-\frac{\theta}{2}\}$, or $\{x_3,x_4\}=\{\frac{1}{6},\theta\}$ and $\{x_1,x_2\}=\{\frac{\theta}{2},\frac{1}{2}-\frac{\theta}{2}\}$, for some $0<\theta\le\frac{1}{2}$;
					
		\item Up to reordering, all other solutions $x_1,x_2,x_3,x_4$ satisfying $0<x_1<x_3\le x_4<x_2\le\frac12$ are in Table \ref{Tab-3.1}.  						
		\end{enumerate}
	\begin{table}[htp]  
		\centering
		\begin{minipage}{0.45\textwidth}
			\begin{tabular}{cccc}
				$x_1$ &$x_2$& $x_3$& $x_4$ \\
				\hline
				1/21&8/21&1/14&3/14\\
				1/14&5/14&2/21&5/21\\
				4/21&10/21&3/14&5/14\\
				1/20&9/20&1/15&4/15\\
				2/15&7/15&3/20&7/20\\
				1/30&3/10&1/15&2/15\\
				1/15&7/15&1/10&7/30\\
				1/10&13/30&2/15&4/15\\
			\end{tabular}
		\end{minipage}
		\raisebox{0.59em}{\begin{minipage}{0.45\textwidth}
				\begin{tabular}{cccc}
					$x_1$&$x_2$&$x_3$&$x_4$\\
					\hline
					4/15&7/15&3/10&11/30\\
					1/30&11/30&1/10&1/10\\
					7/30&13/30&3/10&3/10\\
					1/15&4/15&1/10&1/6\\
					\cline{2-2}
					2/15&\multicolumn{1}{|c|}{7/15}&1/6&3/10\\
					\cline{2-2}
					1/12&5/12&1/10&3/10\\
					1/10&3/10&1/6&1/6\\
				\end{tabular}
		\end{minipage}}
	\caption{Nongeneric solutions in Proposition \ref{proposition-6}. }\label{Tab-3.1}
	\end{table}	   
	\end{proposition}
	\begin{remark}
	  We always have $\sin \frac16 \,  \sin \theta=\sin\frac{\theta}{2} \, \sin (\frac12-\frac{\theta}{2})$. But it is a lengthy check when we transform all angles in this formula to the range $(0,1/2]$ for different ranges of $\theta$ to assure that the double angle relation and the `summing to $1/2$' relation still hold. So we get the 3rd case.
	\end{remark}
	
	\begin{remark}\label{remark23}
		After Case $1$ we can assume all $x_i>0$. Case $2$ and Case $3$ have a common solution $\{ x_1, x_2\}=\{ x_3, x_4\}=\{\frac16,\frac13\}$. 
		
	\end{remark}

    \begin{remark}
	The $7/15$ highlighted in a box in Table \ref{Tab-3.1} was $8/15$ in Myerson's original table, which is an obvious typo since $8/15 > 1/2$. This typo remained in \cite{ck} but the results there were still all right luckily. 
    \end{remark}

\noindent We will use Lemma/Proposition $n'$ to denote the use of Lemma/Proposition $n$ after interchanging $\aaa\leftrightarrow\ddd$ and $\bbb\leftrightarrow\ccc$. 

By Lemma \ref{lemma-1}, the quadrilateral in our tiling can have at most one angle $\ge1$.  Up to the symmetry of interchanging $\aaa\leftrightarrow\ddd$ and $\bbb\leftrightarrow\ccc$, we need only to consider five possibilities: convex (all angles $<1$), concave ($\aaa>1$ or $\bbb>1$), or degenerate ($\aaa=1$ or $\bbb=1$), which will be discussed in the following sections respectively. 
	
	\section{Convex case $\aaa,\bbb,\ccc,\ddd<1$}
	\label{sec-convex}
		
Coolsaet \cite[Theorem 3.2]{ck} classified simple convex \textit{rational} quadrilateral with three equal sides  into $29$ sporadic examples in the first column of Table \ref{Tab-2.2} and $7$ infinite classes (up to interchanging $\aaa\leftrightarrow\ddd$ and $\bbb\leftrightarrow\ccc$): 
\begin{enumerate}
	\item $\aaa=\ccc$ and $\bbb=\ddd$ (and all four sides are equal); 
	\item $\aaa=\ddd$ and $\bbb=\ccc$;
	\item $\aaa=\frac{\ccc}{2}$ and $\ddd=\frac{\bbb}{2}$, $\aaa,\ddd<\frac{1}{2}$;
	\item $\aaa=\frac{3\ccc}{2},\bbb=\frac{1}{3}$ and $\ddd=\frac{2}{3}-\frac{\ccc}{2}$, with $\frac{1}{2}<\ccc<\frac{2}{3}$;
	\item $\aaa=\frac{1}{6}+\frac{\ccc}{2},\bbb=2\ccc$ and $\ddd=\frac{1}{2}+\frac{\ccc}{2}$, with $\frac{1}{3}<\ccc<\frac{1}{2}$; 
	\item $\aaa=\frac{1}{6}+\frac{\ccc}{2},\bbb=2\ccc$ and $\ddd=\frac{1}{2}+\frac{3\ccc}{2}=3\aaa$, with $\frac{4}{15}<\ccc<\frac{1}{3}$; 
	\item $\aaa=\frac{1}{6}+\frac{\ccc}{2},\bbb=2-2\ccc$ and $\ddd=\frac{3}{2}-\frac{3\ccc}{2}$, with $\frac{1}{2}<\ccc<\frac{5}{6}$. 
\end{enumerate}	
In fact Coolsaet assumed additionally that all edges and diagonals are $<1$, and our Lemma \ref{geometry5} shows that such assumptions are satisfied automatically for $a^3b$-tilings. Case $1$ and $2$ are immediately dismissed due to $a\neq b$ and Proposition \ref{symmetric}.  In this section, we will find  all possible tilings from Table \ref{Tab-2.2} and from $5$ remaining cases. 
	
	\subsection*{Sporadic cases in Table \ref{Tab-2.2}} \label{discrete-7}
A quadrilateral is qualified to tile the sphere only if its angle sum is $2+\frac{4}{f}$ for some even integer $f\ge6$, every angle can be extended to a vertex, and there must also exist degree $3$ vertices by the equation \eqref{vcountv}. These basic criteria dismiss most sporadic examples in Table \ref{Tab-2.2}, as the details showing in the second and third columns. There are only three subcases left. But  $(21,8,26,7)/30$ implies $\aaa\cdots=\aaa\bbb^4\ddd$ and  $(17,16,26,11)/30$ implies $\ccc\cdots=\aaa^2\ccc$, both contradicting Balance Lemma \ref{balance}. So only  $(23,10,28,9)/30$ admits a $2$-layer earth map tiling 
$T(12\aaa\ccc\ddd,2\bbb^6)$. In fact the only other possible vertex is $\aaa\bbb\ddd^3$, but  Lemma $\ref{proposition-7}'$ shows that there is no other tilings. This is $f=12$,  $(9, 28, 10, 23)/30$ in Table \ref{Tab-1.1} after	interchanging $\aaa\leftrightarrow\ddd$ and $\bbb\leftrightarrow\ccc$. 
	\begin{table}[htp]  
		\centering     
		\begin{minipage}{0.49\textwidth}
			\scriptsize{\begin{tabular}{ccc}
					$(\aaa,\bbb,\ccc,\ddd)$&$f$& \\
					\hline 
					(29,16,18,23)/42&84&no degree 3 vertex \\
					(31,16,18,23)/42&42&no degree 3 vertex \\
					(35,16,30,17)/42&12&no degree 3 vertex \\
					(37,16,30,17)/42&21/2&$f$ is not even\\
					(35,18,40,17)/42&84/13&$f$ is not even\\
					(11,30,40,7)/42&42&no degree 3 vertex \\
					(29,30,40,23)/42&84/19&$f$ is not even \\
					(49,16,42,17)/60&60&no degree 3 vertex \\
					(53,16,42,17)/60&30&no degree 3 vertex \\
					\hline 
					\multicolumn{1}{|c}{(21,8,26,7)/30}&60&\multicolumn{1}{c|}{$\aaa\cdots=\aaa\bbb^4\ddd$} \\
					\hline 
					(49,18,56,17)/60&12&no degree 3 vertex \\
					\hline 
					\multicolumn{1}{|c}{(23,10,28,9)/30}&12&\multicolumn{1}{c|}{} \\
					\hline 
					(11,7,9,8)/15&12&no degree 3 vertex \\
					(13,7,9,8)/15&60/7&$f$ is not even \\
					(17,14,28,9)/30&15&$f$ is not even\\				
			\end{tabular}}
		\end{minipage}
		\raisebox{0.45em}{\begin{minipage}{0.49\textwidth}
				\scriptsize{\begin{tabular}{ccc}
						$(\aaa,\bbb,\ccc,\ddd)$&$f$& \\
						\hline					
						(25,16,18,19)/30&20/3&$f$ is not even \\
						(23,16,18,19)/30&15/2&$f$ is not even \\
						(25,16,22,17)/30&6&no $\aaa\cdots$ \\
						(27,16,22,17)/30&60/11&$f$ is not even \\
						(23,32,54,13)/60&120&no degree 3 vertex \\
						(31,32,54,19)/60&15&$f$ is not even\\
						\hline 
						\multicolumn{1}{|c}{(17,16,26,11)/30}&12&\multicolumn{1}{c|}{$\ccc\cdots=\aaa^2\ccc$} \\
						\hline 
						(31,36,50,23)/60&12&no degree 3 vertex \\
						(11,9,13,8)/15&60/11&$f$ is not even \\
						(19,18,28,13)/30&20/3&$f$ is not even \\
						(25,18,28,17)/30&30/7&$f$ is not even \\
						(19,42,56,13)/60&24&no degree 3 vertex \\
						(37,42,56,29)/60&60/11&$f$ is not even \\
						(23,22,28,19)/30&15/4&$f$ is not even\\
				\end{tabular}}
		\end{minipage}}
	\caption{29 sporadic convex rational $a^3b$-quadrilaterals.}\label{Tab-2.2}        
	\end{table}

    \subsection*{Case 3. $\aaa=\frac{\ccc}{2},\ddd=\frac{\bbb}{2}$, $\aaa,\ddd<\frac{1}{2}$}

    By Lemma \ref{anglesum}, we get $\frac{2}{3}<\aaa+\ddd\le \frac{8}{9}$ and $\frac{4}{3}<\bbb+\ccc\le \frac{16}{9}$. Without loss of generality, let $\bbb>\ccc$. So we get $\bbb>\frac{2}{3}$, and $\ddd>\frac{1}{3}$ by Lemma \ref{geometry3}. 
    By $\bbb<1$, we get $\ccc>\frac{1}{3},\aaa>\frac{1}{6}$ and $\ddd<\frac{1}{2}$.
    By $R(\bbb^2\cdots)<\bbb=2\ddd,\ccc=2\aaa$ and Parity Lemma, there is no $\bbb^2\cdots$ vertex. Similarly, there is no $\bbb\ddd^2\cdots$ vertex. By $\aaa<R(\bbb\ddd\cdots)<3\aaa$ and Parity Lemma, there is no $\bbb\ddd\cdots$ vertex. By $R(\bbb\cdots)<3\ccc$, $\ccc=2\aaa$ and Parity Lemma, we get $\bbb\cdots=\bbb\ccc^2,\aaa^2\bbb\ccc$ or $\aaa^4\bbb$. They all satisfy   $\#\aaa+\#\ccc\ge2\#\bbb$.   
    If $\aaa^2\bbb\ccc$ or $\aaa^4\bbb$ is a vertex, then $\#\aaa+\#\ccc>2\#\bbb$, contradicting Balance Lemma \ref{balance}. If $\bbb\cdots=\bbb\ccc^2$, then $\#\ccc>\#\bbb$, again a contradiction. We conclude that there is no tiling in this case. 
   
    \subsection*{Case 4. $\aaa=\frac{3\ccc}{2},\bbb=\frac{1}{3},\ddd=\frac{2}{3}-\frac{\ccc}{2},\frac{1}{2}<\ccc<\frac{2}{3}$} 
    
    We have $\frac{3}{4}<\aaa<1$ and $\frac{1}{3}<\ddd<\frac{5}{12}$. By $R(\ccc\cdots)<2\aaa$, $0<R(\aaa\ccc\ddd\cdots)<\bbb,\ccc,2\ddd$ and Parity Lemma, there is no $\aaa\ccc\dots$ vertex.
    By  $0<R(\ccc^3\cdots)<2\bbb,\ccc,2\ddd$ and Parity Lemma, we get  $\ccc^3\cdots=\bbb\ccc^3$. By $2\bbb<R(\ccc^2\cdots)<3\bbb,3\ddd$, $0<R(\ccc^2\ddd^2\cdots)<\bbb$ and Parity Lemma, we get  $\ccc^2\cdots=\bbb\ccc^3$.  By $
    R(\ccc\ddd^2\cdots)=2\bbb<2\ddd, \, 4\bbb<R(\ccc\cdots)<5\bbb$ and Parity Lemma, we get $\ccc\cdots=\bbb\ccc^3$ or $\bbb^2\ccc\ddd^2$.    
    By Balance Lemma, $\bbb\ccc^3$ is a vertex.
    Therefore, $\aaa=\frac{5}{6},\bbb=\frac{1}{3},\ccc=\frac{5}{9}$ and $\ddd=\frac{7}{18}$. Then we get $f=36$. By Parity Lemma, we get all \textit{anglewise vertex combinations} abbreviated as 
    $\text{AVC}\subset\{\aaa^2\bbb,\aaa\ddd^3,\bbb\ccc^3,\bbb^2\ccc\ddd^2,\bbb^6\}$.
    
    If $\bbb^6$ is a vertex, we have the AAD $\bbb^6=\thin^{\ccc}\bbb^{\aaa}\thin^{\aaa}\bbb^{\ccc}\thin\cdots$ or $\thin^{\ccc}\bbb^{\aaa}\thin^{\ccc}\bbb^{\aaa}\thin\cdots$. This gives a vertex $\aaa\thin\aaa\cdots$ or $\aaa\ccc\cdots$, contradicting the AVC. 	Then 
    there is only one solution satisfying Balance Lemma \ref{balance}: $\{14\aaa^2\bbb,8\aaa\ddd^3,10\bbb\ccc^3,6\bbb^2\ccc\ddd^2\}$. 
    
    In Fig.\,\ref{f=36 special}, we have the unique AAD $\bbb^2\ccc\ddd^2=\thin^{\ccc}\ddd_1^{\aaa}\thick^{\aaa}\ddd_2^{\ccc}\thin^{\ccc}\bbb_3^{\aaa}\thin^{\ddd}\ccc_4^{\bbb}\thin^{\aaa}\bbb_5^{\ccc}\thin$ which determines $T_1,T_2,T_3,T_4,T_5$. Then $\aaa_5\bbb_4\cdots=\aaa^2\bbb$ determines $T_6$; $\aaa_4\bbb_6\cdots=\aaa^2\bbb$ determines $T_7$; $\aaa_1\ddd_4\ddd_7\cdots=\aaa\ddd^3$ determines $T_8$;  $\aaa_8\ddd_3\cdots=\aaa\ddd^3$ determines $T_9,T_{10}$. We have $\ccc_2\ccc_3\ccc_{10}\cdots=\thin\ccc_2\thin\ccc_3\thin\ccc_{10}\thin^{\ccc}\bbb_{11}^{\aaa}\thin$ or $\thin\ccc_2\thin\ccc_3\thin\ccc_{10}\thin^{\aaa}\bbb_{11}^{\ccc}\thin$. We might as well take $\ccc_2\ccc_3\ccc_{10}\cdots=\thin\ccc_2\thin\ccc_3\thin\ccc_{10}\thin^{\ccc}\bbb_{11}^{\aaa}\thin$ which determines $T_{11}$. Similarly, we can determine  $T_{12},T_{13},T_{14},T_{15},T_{16},T_{17}$ and $T_{18}$. We have $\bbb_{10}\ccc_{11}\ccc_{13}\cdots=\thin\ccc_{11}\thin\bbb_{10}\thin\ccc_{13}\thin^{\bbb}\ccc_{19}^{\ddd}\thin$ or $\thin\ccc_{11}\thin\bbb_{10}\thin\ccc_{13}\thin^{\ddd}\ccc_{19}^{\bbb}\thin$. We might as well take $\bbb_{10}\ccc_{11}\ccc_{13}\cdots=\thin\ccc_{11}\thin\bbb_{10}\thin\ccc_{13}\thin^{\bbb}\ccc_{19}^{\ddd}\thin$ which determines $T_{19}$. Similarly, we can determine $T_{20},T_{21}$, $\dots,T_{36}$. For other choices of $\ccc_2\ccc_3\ccc_{10}\cdots$ and $\bbb_{10}\ccc_{11}\ccc_{13}\cdots$, we still get this tiling or its equivalent opposite. This is Case $(15, 6, 10, 7)/18$ in Table \ref{Tab-1.1}. \label{discrete-17}
    
    \begin{figure}[htp]
    	\centering
    	\begin{tikzpicture}[>=latex,scale=0.43] 
    		
    		\foreach \a in {-1}
    		{
    			\fill[gray!50]
    			(\a*-9,-4) -- (\a*-9,3) -- (\a*-1,8)-- (\a*-2,6)--(\a*-1,4)--(\a*1,2)--(\a*-1,2)--(\a*-2,0)--(\a*-1,-2)--(\a*0,-4)--(\a*-2,-4)--(\a*-3,-5)--(\a*-1,-8)--(\a*-7,-4)--(\a*-5,-4)--(\a*-5,-2)--(\a*-4,0)--(\a*-3,2)--(\a*-4,3)--(\a*-7,2)--(\a*-9,-4);
    			
    			\fill[gray!50]
    			(\a*8,2)--(\a*8,-5)--(\a*0,-10)--(\a*1,-8)--(\a*0,-6)--(\a*-2,-4)--(\a*0,-4)--(\a*1,-2)--(\a*0,0)--(\a*-1,2)--(\a*1,2)--(\a*2,3)--(\a*0,6)--(\a*6,2)--(\a*4,2)--(\a*4,0)--(\a*3,-2)--(\a*2,-4)--(\a*3,-5)--(\a*6,-4)--(\a*8,2);

    			\draw (\a*-11,3)--(\a*-9,3)--(\a*-1,8)--(\a*-7,2)--(\a*-9,-4)--(\a*0,-10)--(\a*-1,-8)--(\a*-7,-4)--(\a*-6,0)--(\a*-7,2)--(\a*-4,3)--(\a*-2,6)--(\a*-1,4)--(\a*0,6)--(\a*-1,8)--(\a*8,2)
    			(\a*-1,4)--(\a*-3,2)--(\a*-4,0)--(\a*-2,0)--(\a*-1,-2)--(\a*-5,-2)--(\a*-6,0)--(\a*-3,2)--(\a*-1,2)
    			(\a*-5,-2)--(\a*-5,-4)--(\a*-3,-5)--(\a*-1,-8)--(\a*0,-6)--(\a*-2,-4)--(\a*-5,-2)
    			(\a*1,-8)--(\a*0,-6)--(\a*2,-4)--(\a*0,-4)--(\a*-1,-2)--(\a*0,0)--(\a*-1,2)
    			(\a*0,0)--(\a*1,-2)--(\a*3,-2)--(\a*2,-4)--(\a*5,-2)--(\a*6,-4)--(\a*3,-5)--(\a*1,-8)
    			(\a*-1,4)--(\a*1,2)--(\a*4,0)--(\a*4,2)--(\a*2,3)--(\a*0,6)--(\a*6,2)--(\a*5,-2)--(\a*4,0)--(\a*0,0)
    			(\a*8,2)--(\a*6,-4)--(\a*0,-10)--(\a*8,-5)--(\a*10,-5);
    			
    			\draw[line width=1.5] (\a*-9,3)--(\a*-9,-4)--(\a*-5,-4)
    			(\a*-5,-2)--(\a*-4,0)
    			(\a*-4,3)--(\a*-3,2)
    			(\a*-2,0)--(\a*-1,2)--(\a*1,2)--(\a*2,3)
    			(\a*-2,6)--(\a*-1,8)
    			(\a*4,2)--(\a*8,2)--(\a*8,-5)
    			(\a*-3,-5)--(\a*-2,-4)--(\a*0,-4)--(\a*1,-2)
    			(\a*4,0)--(\a*3,-2)
    			(\a*2,-4)--(\a*3,-5)
    			(\a*0,-10)--(\a*1,-8);

    			\node[draw,shape=circle, inner sep=0.5] at (\a*-5,-0.2) {\small $1$};
    			\node[draw,shape=circle, inner sep=0.5] at (\a*-3,-1) {\small $2$};
    			\node[draw,shape=circle, inner sep=0.5] at (\a*-2,-3) {\small $3$};
    			\node[draw,shape=circle, inner sep=0.5] at (\a*-4,-3.6) {\small $4$};
    			\node[draw,shape=circle, inner sep=0.5] at (\a*-5.8,-3) {\small $5$};
    			\node[draw,shape=circle, inner sep=0.5] at (\a*-4.5,-5) {\small $6$};
    			\node[draw,shape=circle, inner sep=0.5] at (\a*-1.5,-6) {\small $7$};
    			\node[draw,shape=circle, inner sep=0.5] at (\a*0,-4.8) {\small $8$};
    			\node[draw,shape=circle, inner sep=0.5] at (\a*1.5,-3) {\small $9$};
    			\node[draw,shape=circle, inner sep=0.5] at (\a*0,-2) {\footnotesize $10$};
    			\node[draw,shape=circle, inner sep=0.5] at (\a*-1,0) {\footnotesize $11$};
    			\node[draw,shape=circle, inner sep=0.5] at (\a*-2.5,1) {\footnotesize $12$};
    			\node[draw,shape=circle, inner sep=0.5] at (\a*2,-1) {\footnotesize $13$};
    			\node[draw,shape=circle, inner sep=0.5] at (\a*3.9,-2) {\footnotesize $14$};
    			\node[draw,shape=circle, inner sep=0.5] at (\a*4,-3.7) {\footnotesize $15$};
    			\node[draw,shape=circle, inner sep=0.5] at (\a*1.5,-6) {\footnotesize $16$};
    			\node[draw,shape=circle, inner sep=0.5] at (\a*0,-8) {\footnotesize $17$};
    			\node[draw,shape=circle, inner sep=0.5] at (\a*3.1,-6) {\footnotesize $18$};
    			\node[draw,shape=circle, inner sep=0.5] at (\a*1,1) {\footnotesize $19$};
    			\node[draw,shape=circle, inner sep=0.5] at (\a*-1,2.8) {\footnotesize $20$};
    			\node[draw,shape=circle, inner sep=0.5] at (\a*4.8,0.8) {\footnotesize $21$};
    			\node[draw,shape=circle, inner sep=0.5] at (\a*2.8,1.6) {\footnotesize $22$};
    			\node[draw,shape=circle, inner sep=0.5] at (\a*0.2,4) {\footnotesize $23$};
    			\node[draw,shape=circle, inner sep=0.5] at (\a*3,3.2) {\footnotesize $24$};
    			\node[draw,shape=circle, inner sep=0.5] at (\a*-2.2,4) {\footnotesize $25$};
    			\node[draw,shape=circle, inner sep=0.5] at (\a*-5,1.5) {\footnotesize $26$};
    			\node[draw,shape=circle, inner sep=0.5] at (\a*-1,6) {\footnotesize $27$};
    			\node[draw,shape=circle, inner sep=0.5] at (\a*-4,4) {\footnotesize $28$};
    			\node[draw,shape=circle, inner sep=0.5] at (\a*-7.5,-2) {\footnotesize $29$};
    			\node[draw,shape=circle, inner sep=0.5] at (\a*-4,-6.7) {\footnotesize $30$};
    			\node[draw,shape=circle, inner sep=0.5] at (\a*6.8,-4.5) {\footnotesize $31$};
    			\node[draw,shape=circle, inner sep=0.5] at (\a*-5,-8) {\footnotesize $32$};
    			\node[draw,shape=circle, inner sep=0.5] at (\a*-7.8,2.7) {\footnotesize $33$};
    			\node[draw,shape=circle, inner sep=0.5] at (\a*6.5,0) {\footnotesize $34$};
    			\node[draw,shape=circle, inner sep=0.5] at (\a*3,4.7) {\footnotesize $35$};
    			\node[draw,shape=circle, inner sep=0.5] at (\a*5,6) {\footnotesize $36$};
    		}

    	\end{tikzpicture}
    	\caption{ $T(14\aaa^2\bbb,8\aaa\ddd^3,10\bbb\ccc^3,6\bbb^2\ccc\ddd^2)$.} \label{f=36 special}
    \end{figure}
	
		\subsection*{Case 5. $\aaa=\frac{1}{6}+\frac{\ccc}{2},\bbb=2\ccc,\ddd=\frac{1}{2}+\frac{\ccc}{2},\frac{1}{3}<\ccc<\frac{1}{2}$}
	
	We have $\frac{1}{3}<\aaa<\ccc<\frac{1}{2}$ and $\frac{2}{3}<\ddd<\bbb<1$. 
	By $R(\bbb^2\cdots)<2\aaa,\bbb,2\ccc,\ddd$ and Parity Lemma, we get $\bbb^2\cdots=\bbb^2\ccc$. By $R(\aaa\bbb\ddd\cdots)< \text{all angles}$,  $0<R(\aaa^2\bbb\cdots)<2\aaa,2\ccc$, $2\ccc<R(\bbb\cdots)<4\ccc,2\ddd$ and Parity Lemma, we get $\bbb\cdots=\aaa\bbb\ddd,\bbb^2\ccc,\aaa^2\bbb\ccc$ or $\bbb\ccc^3$. However  $\aaa\bbb\ddd,\bbb^2\ccc$ or $\bbb\ccc^3$ implies $f=9$ or $15$, contradicting the fact that $f$ is even. So we have only $\bbb\cdots=\aaa^2\bbb\ccc$ with $f=12$. But this again contradicts Balance Lemma \ref{balance}. We conclude that there is no tiling in this case.

	\subsection*{Case 6. $\aaa=\frac{1}{6}+\frac{\ccc}{2},\bbb=2\ccc,\ddd=\frac{1}{2}+\frac{3\ccc}{2}=3\aaa,\frac{4}{15}<\ccc<\frac{1}{3}$}
	
	We have $\frac{3}{10}<\aaa<\frac{1}{3},\frac{8}{15}<\bbb<\frac{2}{3}$ and $\frac{9}{10}<\ddd<1$. 
	By $R(\ddd^2\cdots)< \text{all angles}$, there is no $\ddd^2\cdots$ vertex. By  $0<R(\aaa^3\ddd\cdots)<\text{all angles}$, $0<R(\aaa\bbb\ddd\cdots)< \text{all angles}$, $2\ccc<R(\aaa\ddd\cdots)<3\ccc$ and Parity Lemma, there is no $\ddd\cdots$ vertex, a contradiction. We conclude that there is no tiling in this case.

	\subsection*{Case 7. $\aaa=\frac{1}{6}+\frac{\ccc}{2},\bbb=2-2\ccc,\ddd=\frac{3}{2}-\frac{3\ccc}{2},\frac{1}{2}<\ccc<\frac{5}{6}$} 
	
	By Lemma \ref{anglesum}, we have $\aaa=\frac{7}{12}-\frac{1}{f},\bbb=\frac{1}{3}+\frac{4}{f},\ccc=\frac{5}{6}-\frac{2}{f}$ and $\ddd=\frac{1}{4}+\frac{3}{f}$. So we have $\frac{5}{12}<\aaa<\ccc<\frac{5}{6}$ and $\frac{1}{4}<\ddd<\bbb<1$. 
	
	If $\bbb>\ccc$, then we get $6<f<12$. So we have $\frac{5}{12}<\aaa<\frac12<\ddd<\frac{3}{4}$ and $\frac{1}{2}<\ccc<\frac{2}{3}<\bbb<1$. By $R(\bbb^2\cdots)<2\aaa,\bbb,\ccc,2\ddd$ and Parity Lemma, there is no $\bbb^2\cdots$ vertex. By $0<R(\aaa^2\bbb\cdots)<\text{all angles}$, $R(\aaa\bbb\ddd\cdots)<\text{all angles}$, $R(\bbb\ddd^2\cdots)<\text{all angles}$, $R(\bbb\cdots)=2\ccc$ and Parity Lemma, we get $\bbb\cdots=\aaa\bbb\ddd,\bbb\ccc^2$ or $\bbb\ddd^2$. Suppose $\aaa\bbb\ddd$ or $\bbb\ddd^2$ is a vertex. Then we get $f=\frac{36}{5}$ or $\frac{60}{7}$, a contradiction. So we have $\bbb\cdots=\bbb\ccc^2$. But this again contradicts Balance Lemma \ref{balance}.
	
	Therefore, $\bbb<\ccc$, then we get $f>12$. So we have $\frac{1}{4}<\ddd<\frac12<\aaa<\frac{7}{12}$ and $\frac{1}{3}<\bbb<\frac{2}{3}<\ccc<\frac56$.
	
	If $\aaa^k\bbb^l\ccc^m\ddd^n$ is a vertex, then we have 
	\[
	(\tfrac{7}{12}-\tfrac{1}{f})k+(\tfrac13+\tfrac{4}{f})l+(\tfrac{5}{6}-\tfrac{2}{f})m+(\tfrac{1}{4}+\tfrac{3}{f})n=2.
	\] 
	We also have $\aaa > \frac{1}{2},\bbb>\frac{1}{3},\ccc>\frac{2}{3},\ddd>\frac{1}{4}$.  This implies $k \le 3,l\le 5,m \le 2,n\le 7$. We substitute
	the finitely many combinations of exponents satisfying the bounds into the equation above and solve for $f$. By the angle values and Parity Lemma, we get all possible AVC in Table \ref{Tab-2.1}. Its first row  ``$f=\text{all}$'' means that the angle combinations can be vertices for any $f$; all other rows are mutually exclusive.  All possible tilings based on the AVC of Table \ref{Tab-2.1} are deduced as follows. 
		
	\begin{table}[htp]  
		\centering     
			\begin{tabular}{c|c}
				\hline
				$f$&vertex\\
				\hline
				\hline
				all&$\bbb\ccc^2$, $\aaa^3\ddd$\\
				$20$&$\aaa\bbb^2\ddd$ \\
				$24$&$\bbb^4$, $\bbb\ccc\ddd^2$, \sout{$\bbb\ddd^4$} \\
				$36$&$\aaa^2\bbb^2$, $\aaa\bbb\ddd^3$, \sout{$\bbb^3\ddd^2$}, $\ddd^6$ \\
				$60$&$\bbb^3\ccc$, \sout{$\bbb^5$}, $\bbb^2\ddd^4$, $\ccc\ddd^4$\\
				$84$&$\aaa\bbb^3\ddd$, $\aaa\ddd^5$\\
				$132$& $\bbb^4\ddd^2$, \sout{$\bbb\ddd^6$}\\
			\end{tabular}
		\caption{AVC for  $\aaa=\frac{7}{12}-\frac{1}{f},\bbb=\frac{1}{3}+\frac{4}{f},\ccc=\frac{5}{6}-\frac{2}{f},\ddd=\frac{1}{4}+\frac{3}{f}$.}\label{Tab-2.1} 
	\end{table}

	By Table \ref{Tab-2.1}, there is no $\aaa\ccc\cdots$ vertex. Then $\bbb\ddd^4$ cannot be a vertex, since its AAD gives $\aaa\ccc\cdots$. Similarly, we deduce that each of  $\bbb^3\ddd^2,\bbb^5$ and $\bbb\ddd^6$ cannot be a vertex. If $\text{AVC}\sub\{\bbb\ccc^2,\aaa^3\ddd\}$, there is no solution satisfying Balance Lemma \ref{balance}. All other subcases are discussed as follows. 
	\subsubsection*{Table \ref{Tab-2.1}, $f=20,132$}
	
    There is no solution satisfying Balance Lemma \ref{balance}. 
	
	\subsubsection*{Table \ref{Tab-2.1}, $f=24$}
	
	By Table \ref{Tab-2.1}, we get $\text{AVC}\sub \{\bbb\ccc^2,\aaa^3\ddd,\bbb^4,\bbb\ccc\ddd^2 \}$. There is only one solution satisfying Balance Lemma \ref{balance}: $\{8\bbb\ccc^2,8\aaa^3\ddd,2\bbb^4,8\bbb\ccc\ddd^2\}$.

	In Fig.\,\ref{case f=24}, by AVC, we know $\aaa\ccc\cdots$ is not a vertex. So we have the AAD $\bbb^4=\thin^{\aaa}\bbb_1^{\ccc}\thin^{\ccc}\bbb_2^{\aaa}\thin^{\aaa}\bbb_3^{\ccc}\thin^{\ccc}\bbb_4^{\aaa}\thin$. This determines $T_1,T_2,T_3,T_4$. Then $\ccc_1\ccc_2\cdots=\bbb_5\ccc_1\ccc_2$ determines $T_5$.  Then $\aaa_5\ddd_2\cdots=\aaa_5\aaa_6\aaa_7\ddd_2$ determines $T_6,T_7$. Then $\ccc_5\ddd_1\cdots=\bbb_9\ccc_5\ddd_1\ddd_8$ determines $T_8$. By $\bbb_9$, we have $\aaa_9\ccc_8\cdots$ or $\aaa_9\ddd_5\ddd_6\cdots$, contradicting the AVC.

	\begin{figure}[htp]
		\centering
		\begin{tikzpicture}[>=latex,scale=0.53] 
			\draw (-4,0)--(4,0)
			(-4,2)--(4,2)
			(0,2)--(0,-2)
			(-2,-2)--(2,-2)
			(-2,2)--(0,4)--(4,4)--(4,0)
			(-4,2)--(-4,0);

			\draw[line width=1.5] (-2,2)--(-2,-2)
			(2,2)--(2,-2)
			(0,4)--(2,2);
			
			\draw[dotted] (-4,2)--(-4,4)--(0,4);
								
			\node[draw,shape=circle, inner sep=0.5] at (-1,1) {\small $1$};
			\node[draw,shape=circle, inner sep=0.5] at (1,1) {\small $2$};
			\node[draw,shape=circle, inner sep=0.5] at (-1,-1) {\small $4$};
			\node[draw,shape=circle, inner sep=0.5] at (1,-1) {\small $3$};
			\node[draw,shape=circle, inner sep=0.5] at (0,2.95) {\small $5$};
			\node[draw,shape=circle, inner sep=0.5] at (2.9,2.95) {\small $6$};
			\node[draw,shape=circle, inner sep=0.5] at (3,1) {\small $7$};
			\node[draw,shape=circle, inner sep=0.5] at (-3,1) {\small $8$};
			\node[draw,shape=circle, inner sep=0.5] at (-3,2.95) {\small $9$};

			\node at (0.35,0.4){\small $\bbb$};\node at (-0.35,0.4){\small $\bbb$};\node at (0.35,-0.4){\small $\bbb$};\node at (-0.35,-0.4){\small $\bbb$};
			\node at (0.35,1.6){\small $\ccc$};\node at (-0.35,1.6){\small $\ccc$};
			\node at (0.35,-1.65){\small $\ccc$};\node at (-0.35,-1.65){\small $\ccc$};
			\node at (1.6,0.4){\small $\aaa$};\node at (1.6,-0.4){\small $\aaa$};
			\node at (-1.6,0.4){\small $\aaa$};\node at (-1.6,-0.4){\small $\aaa$};
			\node at (1.6,1.6){\small $\ddd$};\node at (1.6,-1.6){\small $\ddd$};
			\node at (-1.6,1.6){\small $\ddd$};\node at (-1.6,-1.6){\small $\ddd$};
			
			\node at (2.4,0.35){\small $\ddd$};\node at (2.4,1.65){\small $\aaa$};
			\node at (3.6,0.35){\small $\ccc$};\node at (3.6,1.6){\small $\bbb$};
			\node at (3.6,2.4){\small $\bbb$};\node at (3.65,3.65){\small $\ccc$};
			\node at (2.2,2.35){\small $\aaa$};\node at (1,3.65){\small $\ddd$};
			\node at (0,2.2){\small $\bbb$};\node at (-0.1,3.6){\small $\ddd$};
			\node at (1.1,2.3){\small $\aaa$}; \node at (-1.1,2.3){\small $\ccc$};
			
			\node at (-2.4,0.35){\small $\aaa$};\node at (-2.4,1.6){\small $\ddd$};
			\node at (-3.6,0.4){\small $\bbb$};\node at (-3.6,1.6){\small $\ccc$};
			
			\node at (-2.2,2.35){\small $\bbb$};

		\end{tikzpicture}
		\caption{$f=24$, $\text{AVC}=\{\bbb\ccc^2,\aaa^3\ddd,\bbb^4,\bbb\ccc\ddd^2\}$ admit no tiling.} \label{case f=24}
	\end{figure}
	
	\subsubsection*{Table \ref{Tab-2.1}, $f=60$}
	
	By Table \ref{Tab-2.1}, we get $\text{AVC}\sub \{\bbb\ccc^2,\aaa^3\ddd,\bbb^3\ccc,\ccc\ddd^4,\bbb^2\ddd^4 \}$. 
	By AVC, $\aaa\ccc\cdots$ is not a vertex. By Lemma \ref{aadlemma}, the AAD of $\bbb^2\ddd^4$ must be $\thick^{\aaa}\ddd^{\ccc}\thin^{\ccc}\ddd^{\aaa}\thick\cdots$. This gives a vertex $\thin^{\bbb}\ccc^{\ddd}\thin^{\ddd}\ccc^{\bbb}\thin\cdots$. By AVC, we have $\thin^{\bbb}\ccc^{\ddd}\thin^{\ddd}\ccc^{\bbb}\thin\cdots=\thin^{\bbb}\ccc^{\ddd}\thin^{\ddd}\ccc^{\bbb}\thin^{\aaa}\bbb^{\ccc}\thin$. This gives a vertex $\aaa\bbb\cdots$, contradicting the AVC. Therefore, $\bbb^2\ddd^4$ is not a vertex. Similarly, $\ccc\ddd^4$ is not a vertex. Then $\ddd\cdots=\aaa^3\ddd$, contradicting Balance Lemma \ref{balance}.

	\subsubsection*{Table \ref{Tab-2.1}, $f=84$}
	
	By Table \ref{Tab-2.1}, we get $\text{AVC}\sub \{\bbb\ccc^2,\aaa^3\ddd,\aaa\bbb^3\ddd,\aaa\ddd^5 \}$. There is only one solution satisfying Balance Lemma \ref{balance}: $\{42\bbb\ccc^2,20\aaa^3\ddd,14\aaa\bbb^3\ddd,10\aaa\ddd^5\}$.
	
	By AVC, we know $\aaa\ccc\cdots$ is not a vertex. So we have the AAD $\aaa\bbb^3\ddd=\thin^{\bbb}\aaa^{\ddd}\thick^{\aaa}\ddd^{\ccc}\thin^{\ccc}\bbb^{\aaa}\thin^{\aaa}\bbb^{\ccc}\thin^{\ccc}\bbb^{\aaa}\thin$. This gives a vertex $\thin^{\ddd}\ccc^{\bbb}\thin^{\bbb}\ccc^{\ddd}\thin\cdots$. By AVC, we have $\thin^{\ddd}\ccc^{\bbb}\thin^{\bbb}\ccc^{\ddd}\thin\cdots=\thin^{\ddd}\ccc^{\bbb}\thin^{\bbb}\ccc^{\ddd}\thin^{\ccc}\bbb^{\aaa}\thin$. This gives a vertex $\ccc\ddd\cdots$ \textit{not} in the AVC.

	\subsubsection*{Table \ref{Tab-2.1}, $f=36$}
	
	By Table \ref{Tab-2.1}, we get $\text{AVC}\sub \{\bbb\ccc^2,\aaa^3\ddd,\aaa^2\bbb^2,\aaa\bbb\ddd^3,\ddd^6 \}$. 
	We get $\aaa=\frac59,\bbb=\frac{4}{9},\ccc=\frac79,\ddd=\frac13$. 
	
	If $\ddd^6$ is not a vertex, there is only one solution satisfying Balance Lemma \ref{balance}: $\{18\bbb\ccc^2,6\aaa^3\ddd,4\aaa^2\bbb^2,10\aaa\bbb\ddd^3\}$.
	
	In Fig.\,\ref{f=36-2}, we have the unique AAD $\aaa^3\ddd=\thick\aaa_1\thin\aaa_2\thick\aaa_3\thin\ddd_4\thick$ which determines $T_1,T_2,T_3,T_4$. Then $\bbb_3\ccc_4\cdots=\bbb_3\ccc_4\ccc_5$. By $\ccc_5$, $\ccc_3\cdots=\bbb_5\ccc_3\cdots$ determines $T_5$. Then $\bbb_1\bbb_2\cdots=\aaa_6\aaa_7\bbb_1\bbb_2$ determines $T_6,T_7$; $\bbb_7\ccc_1\cdots=\bbb_7\ccc_1\ccc_8$. By $\ccc_8$, $\ccc_7\cdots=\bbb_8\ccc_7\cdots$ determines $T_8$. Then $\aaa_4\ddd_1\ddd_8\cdots=\aaa_4\bbb_{10}\ddd_1\ddd_8\ddd_9$ determines $T_9$. By $\bbb_{10}$, $\ccc_9\cdots=\bbb_{11}\ccc_9\ccc_{10}$ determines $T_{10}$. By $\bbb_{11}$, $\ddd_{10}\cdots=\aaa_{11}\ddd_{10}\cdots$ determines $T_{11}$. Then $\aaa_{10}\bbb_4\ddd_5\cdots=\aaa_{10}\bbb_4\ddd_5\ddd_{12}\ddd_{13}$ determines $T_{12},T_{13}$; $\aaa_{11}\aaa_{13}\ddd_{10}\\\cdots=\aaa_{11}\aaa_{13}\aaa_{14}\ddd_{10}$ determines $T_{14}$; $\bbb_6\ccc_2\cdots=\bbb_6\ccc_2\ccc_{15}$. By $\ccc_{15}$, $\ccc_6\cdots=\bbb_{15}\ccc_6\cdots$ determines $T_{15}$. Then $\ddd_2\ddd_3\ddd_{15}\cdots=\aaa_{16}\bbb_{17}\ddd_2\ddd_3\ddd_{15}$ determines $T_{16}$; $\bbb_5\ccc_3\cdots=\bbb_5\ccc_3\ccc_{17}$ determines $T_{17}$; $\aaa_5\aaa_{12}\ddd_{17}\cdots=\aaa_5\aaa_{12}\aaa_{18}\ddd_{17}$ determines $T_{18}$; $\bbb_{12}\bbb_{18}\cdots=\aaa_{19}\aaa_{20}\bbb_{12}\bbb_{18}$ determines $T_{19}$; $\bbb_{20}\ccc_{12}\ccc_{13}\cdots=\bbb_{20}\ccc_{12}\ccc_{13}$ determines $T_{20}$. Then we get $\bbb_{13}\bbb_{14}\ccc_{20}\cdots$, contradicting the AVC.

	\begin{figure}[htp]
		\centering
		\begin{tikzpicture}[>=latex,scale=0.5] 
			\draw (0,0)--(8,0)--(11,-5)--(7,-10)--(7,-12)--(0,-12)--(0,-10)--(-2,-7)--(-2,-3)--(0,0)
			(6,0)--(9,-5)--(7,-8)--(7,-10)--(2,-10)--(0,-12)
			(9,-5)--(7,-5)--(4,-3)--(-2,-3)
			(4,0)--(2,-3)
			(7,-8)--(4,-7)--(4,-3)
			(0,-3)--(0,-7)--(4,-7)
			(0,-5)--(4,-5)
			(0,-7)--(-2,-7)
			(2,-7)--(2,-10);
			\draw[line width=1.5] (0,0)--(2,-3)--(2,-7)--(0,-10)
			(0,-5)--(-2,-5)
			(4,0)--(7,-5)--(4,-7)--(4,-12)
			(9,-5)--(11,-5);
			\node at (0,-0.6){\small $\aaa$};\node at (0.7,-0.4){\small $\ddd$}; \node at (2,-0.4){\small $\ccc$};\node at (3.3,-0.4){\small $\bbb$};\node at (4,-0.6){\small $\aaa$};\node at (4.7,-0.4){\small $\ddd$};\node at (5.7,-0.4){\small $\ccc$};\node at (6.7,-0.4){\small $\bbb$};\node at (7.7,-0.4){\small $\ccc$};
			\node at (-1.4,-2.6){\small $\bbb$};\node at (0,-2.6){\small $\ccc$};\node at (1.3,-2.6){\small $\ddd$}; \node at (2,-2.4){\small $\aaa$};\node at (2.6,-2.6){\small $\bbb$};
			\node at (4,-2.6){\small $\ccc$};\node at (5.8,-3.8){\small $\ddd$};
			\node at (-1.7,-3.4){\small $\ccc$};\node at (-0.3,-3.4){\small $\bbb$};\node at (0.3,-3.4){\small $\ccc$};\node at (1.7,-3.4){\small $\ddd$};\node at (2.3,-3.4){\small $\ddd$};\node at (3.7,-3.4){\small $\ccc$};\node at (4.3,-3.7){\small $\bbb$};
			\node at (-1.7,-4.6){\small $\ddd$};\node at (-0.3,-4.6){\small $\aaa$};\node at (0.3,-4.6){\small $\bbb$};
			\node at (1.7,-4.6){\small $\aaa$};\node at (2.3,-4.6){\small $\aaa$};\node at (3.7,-4.6){\small $\bbb$};\node at (4.3,-5){\small $\ccc$};\node at (6.3,-5){\small $\aaa$};\node at (7.2,-4.6){\small $\aaa$};\node at (8.3,-4.6){\small $\bbb$};\node at (9.2,-4.6){\small $\aaa$};\node at (10.3,-4.6){\small $\ddd$};
			\node at (-1.7,-5.4){\small $\ddd$};\node at (-0.3,-5.4){\small $\aaa$};\node at (0.3,-5.4){\small $\bbb$};\node at (1.7,-5.4){\small $\aaa$};\node at (2.3,-5.4){\small $\ddd$};\node at (3.7,-5.4){\small $\ccc$};\node at (7.1,-5.4){\small $\aaa$};\node at (8.4,-5.4){\small $\bbb$};\node at (9.1,-5.4){\small $\aaa$};\node at (10.3,-5.4){\small $\ddd$};
			\node at (-1.7,-6.6){\small $\ccc$};\node at (-0.3,-6.6){\small $\bbb$};\node at (0.3,-6.6){\small $\ccc$};\node at (1.7,-6.6){\small $\ddd$};\node at (2.3,-6.6){\small $\aaa$};\node at (3.7,-6.6){\small $\bbb$};\node at (4.3,-6.4){\small $\ddd$};\node at (4.9,-6.9){\small $\ddd$};\node at (6.9,-7.6){\small $\ccc$};\node at (7.3,-8.2){\small $\bbb$};
			\node at (-1.2,-7.4){\small $\bbb$};\node at (0,-7.4){\small $\ccc$};\node at (1.2,-7.4){\small $\ddd$};\node at (1.8,-7.9){\small $\ddd$};\node at (2.3,-7.4){\small $\bbb$};\node at (3.7,-7.4){\small $\aaa$};\node at (4.3,-7.6){\small $\ddd$};
			\node at (0,-9.3){\small $\aaa$};\node at (0.3,-10.2){\small $\aaa$};\node at (1.7,-9.8){\small $\ccc$};\node at (2.3,-9.6){\small $\ccc$};\node at (3.7,-9.6){\small $\ddd$};\node at (4.3,-9.6){\small $\aaa$};\node at (6.7,-9.6){\small $\bbb$};\node at (6.6,-8.3){\small $\ccc$};
			\node at (0.3,-11.2){\small $\bbb$};\node at (0.8,-11.7){\small $\ccc$};\node at (2.1,-10.4){\small $\bbb$};\node at (3.7,-10.4){\small $\aaa$};\node at (4.3,-10.4){\small $\aaa$};\node at (6.7,-10.4){\small $\bbb$};\node at (7.3,-9.2){\small $\ccc$};\node at (3.7,-11.6){\small $\ddd$};\node at (4.3,-11.6){\small $\ddd$};\node at (6.7,-11.6){\small $\ccc$};

			\node[draw,shape=circle, inner sep=0.2] at (1,-6) {\small $1$};
			\node[draw,shape=circle, inner sep=0.2] at (1,-4) {\small $2$};
			\node[draw,shape=circle, inner sep=0.2] at (3,-4) {\small $3$};
			\node[draw,shape=circle, inner sep=0.2] at (3,-6) {\small $4$};
			\node[draw,shape=circle, inner sep=0.2] at (5.2,-5) {\small $5$};
			\node[draw,shape=circle, inner sep=0.2] at (-1,-4) {\small $6$};
			\node[draw,shape=circle, inner sep=0.2] at (-1,-6) {\small $7$};
			\node[draw,shape=circle, inner sep=0.2] at (0,-8.5) {\small $8$};
			\node[draw,shape=circle, inner sep=0.2] at (1,-9.5) {\small $9$};
			\node[draw,shape=circle, inner sep=0.2] at (3,-8.5) {\footnotesize $10$};
			\node[draw,shape=circle, inner sep=0.2] at (2.5,-11.2) {\footnotesize $11$};
			\node[draw,shape=circle, inner sep=0.2] at (6.5,-6.5) {\footnotesize $12$};
			\node[draw,shape=circle, inner sep=0.2] at (5.5,-8.5) {\footnotesize $13$};
			\node[draw,shape=circle, inner sep=0.2] at (5.5,-11.2) {\footnotesize $14$};
			\node[draw,shape=circle, inner sep=0.2] at (0,-1.6) {\footnotesize $15$};
			\node[draw,shape=circle, inner sep=0.2] at (2,-1.6) {\footnotesize $16$};
			\node[draw,shape=circle, inner sep=0.2] at (4,-1.6) {\footnotesize $17$};
			\node[draw,shape=circle, inner sep=0.2] at (6.5,-2.5) {\footnotesize $18$};
			\node[draw,shape=circle, inner sep=0.2] at (8.5,-2.5) {\footnotesize $19$};
			\node[draw,shape=circle, inner sep=0.2] at (9,-6.5) {\footnotesize $20$};
			
		\end{tikzpicture}
		\caption{ $f=36$, $\text{AVC}=\{\bbb\ccc^2,\aaa^3\ddd,\aaa^2\bbb^2,\aaa\bbb\ddd^3\}$ admit no tiling.} \label{f=36-2}
	\end{figure}
		
	Therefore, $\ddd^6$ is a vertex. We have the unique AAD for $\ddd^6=\thick\ddd_1\thin\ddd_2\thick\cdots$  which determines $T_1,T_2,T_3,T_4,T_5,T_6$. Then $\ccc_1\ccc_2\cdots=\bbb_7\ccc_1\ccc_2$ determines $T_7$. So $\aaa_2\aaa_3\cdots=\aaa^2\bbb^2$ or $\aaa^3\ddd$,  shown in Fig.\,\ref{f=36 1} and \ref{f=36} respectively.
		
	In Fig.\,\ref{f=36 1}, $\aaa_2\aaa_3\cdots=\aaa_2\aaa_3\bbb_8\bbb_9$. Then $\bbb_2\ccc_7\cdots=\bbb_2\ccc_7\ccc_8$ determines $T_8$. By $\bbb_9$, $\aaa_8\cdots=\aaa_8\aaa_9\cdots$ determines $T_9$. Then $\ccc_3\ccc_4\cdots=\bbb_{10}\ccc_3\ccc_4,\bbb_3\ccc_9\cdots=\bbb_3\ccc_9\ccc_{10}$ determine $T_{10}$; $\thick\aaa_8\thin\aaa_9\thick\cdots=\aaa_8\aaa_9\aaa_{12}\ddd_{11}$ determines $T_{11},T_{12}$; $\bbb_{12}\ccc_{11}\cdots$ $=\bbb_{12}\ccc_{11}\ccc_{13}$. By $\ccc_{13}$, $\ccc_{12}\cdots=\bbb_{13}\ccc_{12}\cdots$ determines $T_{13}$. We have $\aaa_{10}\bbb_4\cdots=\aaa^2\bbb^2$ or $\aaa\bbb\ddd^3$. If $\aaa_{10}\bbb_4\cdots=\aaa\bbb\ddd^3$, then we get $\aaa_4\aaa_5\ccc\cdots$, contradicting the AVC. Therefore,  $\aaa_{10}\bbb_4\cdots=\aaa_{10}\aaa_{14}\bbb_4\bbb_{15}$. This determines $T_{14}$. By $\bbb_{15}$, $\aaa_4\aaa_5\cdots=\aaa_4\aaa_5\aaa_{15}\ddd_{16}$ determines $T_{15},T_{16}$; $\ddd_9\ddd_{10}\ddd_{12}\ddd_{14}\cdots=\ddd_9\ddd_{10}\ddd_{12}\ddd_{14}\ddd_{17}\ddd_{18}$ determines $T_{17},T_{18}$; $\bbb_{14}\ccc_{15}\cdots=\bbb_{14}\ccc_{15}\ccc_{19}$, $\ccc_{14}\ccc_{18}\cdots=\bbb_{19}\ccc_{14}\ccc_{18}$ determine $T_{19}$; $\aaa_{16}\ddd_{15}\ddd_{19}\cdots$ $=$ $\aaa_{16}\ddd_{15}\ddd_{19}\ddd_{20}\cdots$ determines $T_{20}$; $\aaa_{19}\aaa_{20}\bbb_{18}\cdots$ $=$ $\aaa_{19}\aaa_{20}\bbb_{18}\\\bbb_{21}$. By $\bbb_{21}$, $\aaa_{17}\aaa_{18}\cdots=\aaa_{17}\aaa_{18}\aaa_{21}\ddd_{22}$ determines $T_{21},T_{22}$. Then we get $\aaa_{13}\bbb_{17}\ccc_{22}\cdots$, contradicting the AVC.
		
	\begin{figure}[htp]
		\centering
		\begin{tikzpicture}[>=latex,scale=0.47] 
					\draw (0,-2)--(26,-2)
					(4,0)--(4,-2)
					(12,0)--(12,-2)
					(22,0)--(22,-2)
					(6,-2)--(4,-5)--(6,-8)--(10,-8)--(10,-13)--(19,-13)
					(8,-8)--(8,-2)
					(10,-11)--(12,-11)
					(10,-8)--(12,-5)--(10,-2)
					(12,-5)--(14,-7)--(16,-10)--(12,-11)
					(14,-2)--(16,-5)--(14,-7)--(16,-10)--(19,-10)--(19,-13)
					(16,-5)--(20,-5)--(20,-2)
					(18,-5)--(19,-7)--(19,-10);
					\draw[line width=1.5] (0,0)--(0,-2)
					(8,0)--(8,-2)
					(18,0)--(18,-5)--(16,-10)
					(26,0)--(26,-2)
					(2,-2)--(4,-5)--(12,-5)--(14,-2)
					(12,-5)--(12,-11)--(19,-13)
					(6,-8)--(10,-11);
					\node at (2,0){\small $\ddd$}; \node at (6,0){\small $\ddd$}; \node at (10,0){\small $\ddd$}; \node at (15,0){\small $\ddd$}; \node at (20,0){\small $\ddd$}; \node at (24,0){\small $\ddd$}; 
					\node at (0.4,-1.6){\small $\aaa$}; \node at (7.6,-1.6){\small $\aaa$}; \node at (8.4,-1.6){\small $\aaa$}; \node at (17.6,-1.6){\small $\aaa$}; \node at (18.4,-1.6){\small $\aaa$}; \node at (25.6,-1.6){\small $\aaa$}; 
					\node at (2,-1.6){\small $\bbb$}; \node at (6,-1.6){\small $\bbb$}; \node at (10,-1.6){\small $\bbb$}; \node at (14,-1.6){\small $\bbb$}; \node at (20,-1.6){\small $\bbb$}; \node at (24,-1.6){\small $\bbb$}; 
					\node at (3.6,-1.6){\small $\ccc$}; \node at (4.4,-1.6){\small $\ccc$}; \node at (11.6,-1.6){\small $\ccc$}; \node at (12.4,-1.6){\small $\ccc$}; \node at (21.6,-1.6){\small $\ccc$}; \node at (22.4,-1.6){\small $\ccc$}; 
					\node at (6.2,-2.4){\small $\ccc$}; \node at (7.6,-2.4){\small $\bbb$}; \node at (8.4,-2.4){\small $\bbb$}; \node at (9.8,-2.4){\small $\ccc$}; \node at (10.6,-2.4){\small $\ccc$}; \node at (12,-2.4){\small $\bbb$}; \node at (13.3,-2.4){\small $\aaa$}; \node at (14,-2.6){\small $\aaa$}; \node at (14.8,-2.4){\small $\bbb$}; \node at (17.6,-2.4){\small $\aaa$}; \node at (18.4,-2.4){\small $\ddd$}; \node at (19.6,-2.4){\small $\ccc$}; 
					\node at (4.8,-4.6){\small $\ddd$}; \node at (7.6,-4.6){\small $\aaa$}; \node at (8.4,-4.6){\small $\aaa$}; \node at (11.3,-4.6){\small $\ddd$}; \node at (12,-4.3){\small $\ddd$}; \node at (12.6,-5){\small $\ddd$}; \node at (15.4,-5){\small $\bbb$}; \node at (16.2,-4.6){\small $\ccc$}; \node at (17.6,-4.6){\small $\ddd$}; \node at (18.4,-4.6){\small $\aaa$}; \node at (19.6,-4.6){\small $\bbb$}; 
			    	\node at (4.6,-5.4){\small $\aaa$}; \node at (7.6,-5.4){\small $\ddd$}; \node at (8.4,-5.4){\small $\aaa$}; \node at (11.2,-5.4){\small $\ddd$}; \node at (11.6,-6.2){\small $\ddd$}; 
			    	\node at (12.4,-6){\small $\ddd$}; \node at (16.1,-5.4){\small $\ccc$}; \node at (17.5,-5.4){\small $\ddd$}; \node at (18,-6){\small $\ddd$}; 
			    	\node at (6.3,-7.6){\small $\bbb$}; \node at (7.6,-7.6){\small $\ccc$}; \node at (8.4,-7.6){\small $\bbb$}; \node at (9.8,-7.6){\small $\ccc$}; \node at (10.3,-8.2){\small $\ccc$}; \node at (6.3,-7.6){\small $\bbb$}; \node at (14,-6.6){\small $\ccc$}; \node at (13.7,-7.4){\small $\ccc$}; 
			    	\node at (14.6,-7){\small $\bbb$}; \node at (18.6,-7.2){\small $\ccc$}; 
					\node at (7,-8.35){\small $\ddd$}; \node at (8,-8.4){\small $\ccc$}; \node at (9.6,-8.5){\small $\bbb$}; \node at (9.6,-10.2){\small $\aaa$}; \node at (10.4,-10.6){\small $\bbb$}; \node at (11.6,-10.6){\small $\aaa$}; \node at (12.4,-10.4){\small $\aaa$}; \node at (13,-11){\small $\aaa$};
					\node at (11.8,-11.4){\small $\ddd$};  \node at (10.4,-11.4){\small $\ccc$}; \node at (10.4,-12.6){\small $\bbb$}; \node at (16,-12.6){\small $\aaa$}; \node at (15.2,-9.7){\small $\bbb$}; \node at (15.9,-9.1){\small $\aaa$}; \node at (16.6,-9.6){\small $\aaa$}; \node at (16,-10.4){\small $\bbb$}; \node at (18.6,-9.6){\small $\bbb$}; \node at (18.6,-10.4){\small $\ccc$}; \node at (18.6,-12.3){\small $\ddd$}; 
					
					\node at (4,-2.4){\small $\bbb$}; \node at (4,-4.4){\small $\ddd$}; \node at (2.7,-2.4){\small $\aaa$}; \node at (5.2,-2.4){\small $\ccc$};

					\node[draw,shape=circle, inner sep=0.5] at (2,-0.8) {\small $1$};
					\node[draw,shape=circle, inner sep=0.5] at (6,-0.8) {\small $2$};
					\node[draw,shape=circle, inner sep=0.5] at (10,-0.8) {\small $3$};
					\node[draw,shape=circle, inner sep=0.5] at (15,-0.8) {\small $4$};
					\node[draw,shape=circle, inner sep=0.5] at (20,-0.8) {\small $5$};
					\node[draw,shape=circle, inner sep=0.5] at (24,-0.8) {\small $6$};
                    \node[draw,shape=circle, inner sep=0.5] at (4,-3.5) {\small $7$};
                    \node[draw,shape=circle, inner sep=0.5] at (6.5,-3.5) {\small $8$};
                    \node[draw,shape=circle, inner sep=0.5] at (9.5,-3.5) {\small $9$};
                    \node[draw,shape=circle, inner sep=0.5] at (12,-3.5) {\footnotesize $10$};
                    \node[draw,shape=circle, inner sep=0.5] at (6.5,-6.5) {\footnotesize $11$};
                    \node[draw,shape=circle, inner sep=0.5] at (9.5,-6.5) {\footnotesize $12$};
                    \node[draw,shape=circle, inner sep=0.5] at (8.8,-9.2) {\footnotesize $13$};
                    \node[draw,shape=circle, inner sep=0.5] at (16.5,-3.5) {\footnotesize $15$};
                    \node[draw,shape=circle, inner sep=0.5] at (19,-3.5) {\footnotesize $16$};
                    \node[draw,shape=circle, inner sep=0.5] at (14,-4.5) {\footnotesize $14$};
                    \node[draw,shape=circle, inner sep=0.5] at (11,-9.2) {\footnotesize $17$};
                    \node[draw,shape=circle, inner sep=0.5] at (13.5,-9.2) {\footnotesize $18$};
                    \node[draw,shape=circle, inner sep=0.5] at (16,-7) {\footnotesize $19$};
                    \node[draw,shape=circle, inner sep=0.5] at (18,-8) {\footnotesize $20$};
                    \node[draw,shape=circle, inner sep=0.5] at (17,-11.5) {\footnotesize $21$};
                    \node[draw,shape=circle, inner sep=0.5] at (12,-12.3) {\footnotesize $22$};
			
		\end{tikzpicture}
		\caption{Case $\aaa_2\aaa_3\cdots=\aaa^2\bbb^2$ admits no tiling.} \label{f=36 1}
	\end{figure}
	
	In Fig.\,\ref{f=36}, $\aaa_2\aaa_3\cdots=\aaa^3\ddd$.  We have $\bbb_2\ccc_7\cdots=\bbb_2\ccc_7\ccc_8$. By $\ccc_8$, $\aaa_2\aaa_3\cdots=\aaa_2\aaa_3\aaa_9\ddd_8$ determines $T_8,T_9$. Then $\bbb_3\bbb_9\cdots=\aaa_{10}\aaa_{11}\bbb_3\bbb_9$ determines $T_{10},T_{11}$; $\bbb_{11}\ccc_9\cdots=\bbb_{11}\ccc_9\ccc_{12}$. By $\ccc_{12}$, $\ccc_{11}\cdots=\bbb_{12}\ccc_{11}\cdots$ determines $T_{12}$. Then $\aaa_8\ddd_9\ddd_{12}\cdots=\aaa_8\bbb_{14}\ddd_9\ddd_{12}\ddd_{13}$ determines $T_{13}$. By $\bbb_{14}$, $\ccc_{13}\cdots=\bbb_{15}\ccc_{13}\ccc_{14}$ determines $T_{14}$. By $\bbb_{15}$, $\ddd_{14}\cdots=\aaa_{15}\ddd_{14}\cdots$ determines $T_{15}$. We have $\aaa_{12}\aaa_{13}\cdots=\aaa^3\ddd$ or $\aaa^2\bbb^2$. If $\aaa_{12}\aaa_{13}\cdots=\aaa^3\ddd$, then we get $\bbb_{13}\ccc_{15}\bbb\cdots$ or $\bbb_{12}\ccc_{11}\bbb\cdots$, contradicting the AVC. Therefore, $\aaa_{12}\aaa_{13}\cdots=\aaa^2\bbb^2$. Then $\bbb_{13}\ccc_{15}\cdots=\bbb_{13}\ccc_{15}\ccc_{16},\bbb_{12}\ccc_{11} \cdots=\bbb_{12}\ccc_{11}\ccc_{17}$ determine $T_{16},T_{17}$. The argument started at $T_7$ can be repeated at $T_{10}$. Two repetitions give a unique tiling  $T(18\bbb\ccc^2,6\aaa^3\ddd,6\aaa^2\bbb^2,6\aaa\bbb\ddd^3,2\ddd^6)$. This is Case $(5, 4, 7, 3)/9$ in Table \ref{Tab-1.1}. \label{discrete-16}

	\begin{figure}[htp]
		\centering
		\begin{tikzpicture}[>=latex,scale=0.5] 
			\foreach \a in {0,1,2}
			{
				\begin{scope}[xshift=8*\a cm] 
					\draw (0,0)--(8,0)
					(2,0)--(2,2)
					(0,2)--(8,2)
					(4,2)--(4,4)
					(4,0)--(4,-2)--(2,-2)--(2,-4)
					(6,0)--(6,2)
					(4,-2)--(10,-2)
					(8,0)--(8,-2)
					(10,-2)--(10,-4)
					(0,0)--(0,-2)--(2,-2);
					\draw[line width=1.5] (0,4)--(0,0)--(4,-2)
					(2,2)--(6,0)--(6,-4)
					(8,4)--(8,0);
					\node at (2,4){\small $\ddd$}; \node at (6,4){\small $\ddd$};
					\node at (0.4,2.4){\small $\aaa$}; \node at (2,2.4){\small $\bbb$};\node at (3.5,2.4){\small $\ccc$};\node at (4.4,2.4){\small $\ccc$};\node at (6,2.4){\small $\bbb$};
					\node at (7.6,2.4){\small $\aaa$};\node at (0.4,1.6){\small $\aaa$};\node at (1.6,1.6){\small $\bbb$};\node at (2.4,1.4){\small $\aaa$}; \node at (3.3,1.73){\small $\aaa$};
					\node at (4.1,1.6){\small $\bbb$};\node at (5.6,1.6){\small $\ccc$};\node at (6.4,1.6){\small $\ccc$};\node at (7.6,1.6){\small $\ddd$};\node at (0.4,0.4){\small $\ddd$};
					\node at (1.6,0.4){\small $\ccc$};\node at (2.4,0.4){\small $\bbb$};\node at (3.8,0.4){\small $\ccc$};\node at (4.5,0.35){\small $\ddd$};
					\node at (5.6,0.6){\small $\ddd$}; \node at (6.4,0.4){\small $\bbb$};\node at (7.6,0.4){\small $\aaa$};
					
					\node at (1.3,-0.32){\small $\ddd$};\node at (2.2,-0.35){\small $\ccc$};\node at (3.7,-0.4){\small $\bbb$}; \node at (3.7,-1.5){\small $\aaa$};\node at (4.4,-0.4){\small $\ccc$};\node at (5.6,-0.4){\small $\ddd$};\node at (6.4,-0.4){\small $\aaa$};\node at (7.6,-0.4){\small $\bbb$};\node at (0.3,-0.6){\small $\ddd$};
					\node at (4.4,-1.6){\small $\bbb$};\node at (5.6,-1.6){\small $\aaa$};\node at (6.4,-1.6){\small $\ddd$};\node at (7.6,-1.6){\small $\ccc$};\node at (0.3,-1.6){\small $\ccc$};\node at (1.8,-1.6){\small $\bbb$};\node at (2.8,-1.8){\small $\aaa$};
					\node at (2.4,-2.4){\small $\ccc$};\node at (4,-2.4){\small $\bbb$};\node at (5.6,-2.4){\small $\aaa$};\node at (6.4,-2.4){\small $\aaa$};\node at (8,-2.4){\small $\bbb$};\node at (9.6,-2.4){\small $\ccc$};\node at (4,-4){\small $\ddd$};\node at (8,-4){\small $\ddd$};
					
				\end{scope}
			}
		
		    \node[draw,shape=circle, inner sep=0.2] at (2,3.25) {\small $1$};
		    \node[draw,shape=circle, inner sep=0.2] at (6,3.25) {\small $2$};
		    \node[draw,shape=circle, inner sep=0.2] at (10,3.25) {\small $3$};
		    \node[draw,shape=circle, inner sep=0.2] at (14,3.25) {\small $4$};
		    \node[draw,shape=circle, inner sep=0.2] at (18,3.25) {\small $5$};
		    \node[draw,shape=circle, inner sep=0.2] at (22,3.25) {\small $6$};
		    \node[draw,shape=circle, inner sep=0.2] at (4.8,1.3) {\small $7$};
		    \node[draw,shape=circle, inner sep=0.2] at (7,1) {\small $8$};
		    \node[draw,shape=circle, inner sep=0.2] at (9,1) {\small $9$};
		    \node[draw,shape=circle, inner sep=0.1] at (13,1.3) {\footnotesize $10$};
		    \node[draw,shape=circle, inner sep=0.1] at (11.1,0.8) {\footnotesize $11$};
		    \node[draw,shape=circle, inner sep=0.1] at (11,-0.9) {\footnotesize $12$};
		    \node[draw,shape=circle, inner sep=0.1] at (9.1,-1.3) {\footnotesize $13$};
		    \node[draw,shape=circle, inner sep=0.1] at (7,-1) {\footnotesize $14$};
		    \node[draw,shape=circle, inner sep=0.1] at (8,-3.2)  {\footnotesize $15$};
		    \node[draw,shape=circle, inner sep=0.1] at (12,-3.2)   {\footnotesize $16$};
		    \node[draw,shape=circle, inner sep=0.1] at (13,-1)  {\footnotesize $17$};
		    \node[draw,shape=circle, inner sep=0.1] at (3.1,0.8) {\footnotesize $18$};
		    \node[draw,shape=circle, inner sep=0.1] at (5,-1) {\footnotesize $19$};
		    \node[draw,shape=circle, inner sep=0.1] at (1,1) {\footnotesize $20$};
		    \node[draw,shape=circle, inner sep=0.1] at (3,-0.9) {\footnotesize $21$};
		    \node[draw,shape=circle, inner sep=0.1] at (23,1) {\footnotesize $22$};
		    \node[draw,shape=circle, inner sep=0.1] at (21,1.3) {\footnotesize $23$};
		    \node[draw,shape=circle, inner sep=0.1] at (19.1,0.8) {\footnotesize $24$};
		    \node[draw,shape=circle, inner sep=0.1] at (17,1) {\footnotesize $25$};
		    \node[draw,shape=circle, inner sep=0.1] at (15,1) {\footnotesize $26$};
		    \node[draw,shape=circle, inner sep=0.1] at (15,-1) {\footnotesize $27$};
            \node[draw,shape=circle, inner sep=0.1] at (17.1,-1.3) {\footnotesize $28$};
            \node[draw,shape=circle, inner sep=0.1] at (19,-0.9) {\footnotesize $29$};
		    \node[draw,shape=circle, inner sep=0.1] at (1.1,-1.3){\footnotesize $30$};
		    \node[draw,shape=circle, inner sep=0.1] at (23,-1) {\footnotesize $31$};
		    \node[draw,shape=circle, inner sep=0.1] at (21,-1) {\footnotesize $32$};
		    \node[draw,shape=circle, inner sep=0.1] at (16,-3.2)  {\footnotesize $33$};
		    \node[draw,shape=circle, inner sep=0.1] at (20,-3.2)  {\footnotesize $34$};
		    \node[draw,shape=circle, inner sep=0.1] at (24,-3.2) {\footnotesize $35$};
		    \node[draw,shape=circle, inner sep=0.1] at (4,-3.2) {\footnotesize $36$};

		\end{tikzpicture}
		\caption{Case $\aaa_2\aaa_3\cdots=\aaa^3\ddd$ admits $T(18\bbb\ccc^2,6\aaa^3\ddd,6\aaa^2\bbb^2,6\aaa\bbb\ddd^3,2\ddd^6)$.} \label{f=36}
	\end{figure}

	
	
	\section{Concave case $\aaa>1$}
	\label{sec-concave-a}
	
	An $a^3b$-quadrilateral with $\aaa>1,\bbb,\ccc,\ddd<1$ is shown in Fig.\,\ref{quadrilateral-2}, where $\phi=\angle ACB=\angle BAC$ and $\psi=\angle BDC=\angle CBD$.  We first prove some basic facts. Recall that Lemma \ref{geometry4} implies $\bbb+2\ccc>1$.
	\begin{figure}[htp]
		\centering
		\begin{tikzpicture}[>=latex,scale=0.45] 
			\draw (5,-4)--(0,0)--(-5,-4)--(1,-3);

			\draw[line width=1.5] (5,-4)--(1,-3);
			
			\draw[dotted] (0,0)--(1,-3)
			(5,-4)--(-5,-4);
			
			\node at (0,0.5){\small $C$};  \node at (0.8,-3.5){\small $A$};  \node at (-5.05,-4.45){\small $B$};
			\node at (4.9,-4.45){\small $D$};   \node at (0,-1.5){\small $x$}; \node at (0,-4.5){\small $y$};

			\node at (-3.3,-3.2){\small $\bbb$};  \node at (3.3,-3.2){\small $\ddd$};
			
			\node at (-0.2,-0.5){\small $\phi$};  \node at (0.5,-2.7){\small $\phi$};
			
			\node at (-4.1,-3.7){\small $\psi$}; \node at (4.1,-3.7){\small $\psi$};

			\node at (-2.5,-1.5){\small $a$}; \node at (2.3,-2.9){\small $b$}; \node at (2.5,-1.5){\small $a$};\node at (-1.8,-3.1){\small $a$};

		\end{tikzpicture}
		\caption{$a^3b$-quadrilateral with $\aaa>1,\bbb,\ccc,\ddd<1$.} \label{quadrilateral-2}
\end{figure}

\begin{lemma}\label{proposition-4}
	In an $a^3b$-tiling with $\aaa>1$, we have $a>b$, $\aaa>1>\ccc>\bbb>\ddd$, $\ccc>\frac{1}{3}$ and $\ddd<\frac12$.

\end{lemma}
\begin{proof}
$\aaa>1>\ccc$ implies $\angle CAD=\aaa-\phi>\ccc-\phi=\angle ACD$. So $a>b$. Then $\angle ABD<\angle ADB$. By $\angle CBD=\angle BDC$, we get $\bbb>\ddd$. By Lemma \ref{geometry1}, we have $\bbb<\ccc$. By Lemma \ref{geometry4}, we have $\bbb+2\ccc>1$, so $\ccc>\frac13$.

If $\ddd\ge\frac12$, by $\aaa>1>\ccc>\bbb>\ddd$, the sum of $\aaa$ with any two angles is $>2$ and there is no $\aaa\cdots$ vertex, a contradiction. 
\end{proof}

	\begin{lemma}\label{proposition-5}
		In an $a^3b$-tiling with $\aaa>1$, if $\aaa\bbb\ddd$ appears, then $\ccc=\frac{2}{3},\frac12$ or $\frac25$.
	\end{lemma}

	\begin{proof}
		By $\aaa+\bbb+\ddd=2$, we get $\ccc=\frac{4}{f}$ with $f\ge6$ being even. By Lemma \ref{proposition-4}, $\ccc>\frac{1}{3}$. So $f<12$. Then $f=6,8,10$ and $\ccc=\frac{2}{3},\frac12$ or $\frac25$.
	\end{proof}
    
    To find rational $a^3b$-quadrilaterals by solving \eqref{4-7} or \eqref{4-8} via Proposition \ref{proposition-6}, we have to transform $\aaa-\frac{\ccc}{2},\frac\bbb2,\frac\ccc2,\ddd-\frac\bbb2,\aaa+\frac\ccc2,-\ddd-\frac\bbb2$ to the range $[0,\frac12]$.
	For \eqref{4-7}, by Lemma \ref{proposition-4},  we have $\frac12<\aaa-\frac{\ccc}{2}<2$, $0<\frac{\bbb}{2},\frac{\ccc}{2}<\frac12$, $-\frac12<\ddd-\frac{\bbb}{2}<\frac12$.
	For \eqref{4-8}, by Lemma \ref{proposition-4}, we have  $0<\frac{\bbb}{2},\frac{\ccc}{2}<\frac12,0<\ddd+\frac{\bbb}{2}<1$, which implies $\sin (\aaa+\frac\ccc2)<0$. By $1<\aaa+\frac{\ccc}{2}<\frac52$, we get $1<\aaa+\frac{\ccc}{2}<2$.  		
	Thus we have to consider the following seven choices:
	\begin{small}
		\begin{equation}
			\{ x_1, x_2, x_3, x_4\}=\begin{cases}  
				\{1-\aaa+\frac{\ccc}{2},\frac{\bbb}{2},\frac{\ccc}{2},\ddd-\frac{\bbb}{2}\}, \\ 
				\{-1+\aaa-\frac{\ccc}{2},\frac{\bbb}{2},\frac{\ccc}{2},-\ddd+\frac{\bbb}{2}\}, \\
				\{2-\aaa+\frac{\ccc}{2}, \frac{\bbb}{2}, \frac{\ccc}{2}, -\ddd+\frac{\bbb}{2}\},\\
				\{-1+\aaa+\frac{\ccc}{2},\frac{\bbb}{2},\frac{\ccc}{2},\ddd+\frac{\bbb}{2}\}, \\
				\{-1+\aaa+\frac{\ccc}{2},\frac{\bbb}{2},\frac{\ccc}{2},1-\ddd-\frac{\bbb}{2}\}, \\
				\{2-\aaa-\frac{\ccc}{2},\frac{\bbb}{2},\frac{\ccc}{2},\ddd+\frac{\bbb}{2}\}, \\
				\{2-\aaa-\frac{\ccc}{2},\frac{\bbb}{2},\frac{\ccc}{2},1-\ddd-\frac{\bbb}{2}\}. \\
			\end{cases} \label{4-9}
		\end{equation}
	\end{small}
	
	We will match these choices with four cases of solutions in Proposition \ref{proposition-6} as follows.

	\subsection*{Case $1$: $x_1  x_2=x_3 x_4=0$}
	By $-\frac12<\ddd-\frac\bbb2<\frac12$, $\frac12<\aaa-\frac\ccc2<2$ and $0<\ddd+\frac\bbb2<1$, the only solution of $x_1  x_2=x_3 x_4=0$ for \eqref{4-9} comes from $\aaa-\frac\ccc2=1$ and $\ddd-\frac\bbb2=0$.
	By Lemma \ref{proposition-5'}, we know that $\aaa\bbb\ddd$ or $\aaa\ccc\ddd$ is a vertex.

	If $\aaa\bbb\ddd$ is a vertex, we get three subcases by Lemma \ref{proposition-5}:
	
	$1$. $\aaa=\frac43,\bbb=\frac49,\ccc=\frac23,\ddd=\frac29$; (Case $(12,4,6,2)/9$ in Table \ref{Tab-1.1})
	
	$2$. $\aaa=\frac54,\bbb=\ccc=\frac12,\ddd=\frac14$;
	
	$3$. $\aaa=\frac65,\bbb=\frac{8}{15},\ccc=\frac25,\ddd=\frac{4}{15}$.
	
	For the second and third subcases, we have $\bbb\ge\ccc$, contradicting $\bbb<\ccc$ in Lemma \ref{proposition-4}. In the first subcase, we have $\aaa\cdots=\aaa\bbb\ddd$ or $\aaa\ddd^3$. By $\#\aaa=\#\ddd$, we get $\aaa\cdots=\aaa\bbb\ddd$. There is only one solution satisfying Balance Lemma \ref{balance}: $\{6\aaa\bbb\ddd,2\ccc^3\}$, and it gives a $2$-layer earth map tiling by Lemma \ref{proposition-7}. \label{discrete-3}
	
	If $\aaa\ccc\ddd$ is a vertex, then we get $\aaa=\frac43-\frac{2}{3f},\bbb=\frac{4}{f},\ccc=\frac23-\frac{4}{3f},\ddd=\frac{2}{f}$. 
	By Lemma \ref{proposition-4}, $\frac4f=\bbb<\ccc=\frac23-\frac{4}{3f}$. This implies $f>8$.	
	By the angle values, Parity Lemma and Lemma \ref{proposition-5'}, we get \[\text{AVC}\sub\{\aaa\ccc\ddd,\bbb\ccc^3,\aaa\bbb^{\frac{f-2}{6}}\ddd,\bbb^{\frac{f+4}{6}}\ccc^2,\bbb^{\frac{f+1}{3}}\ccc,\bbb^{\frac{f}{2}}\}.\]

	When $f=6k(k\ge2)$ or $6k+4(k\ge1)$, we have $\text{AVC}\sub\{\aaa\ccc\ddd,\bbb\ccc^3,\bbb^{\frac{f}{2}}\}$, and the only solution satisfying Balance Lemma \ref{balance} is $\{ f\aaa\ccc\ddd,2\bbb^{\frac{f}{2}}\}$  which gives a $2$-layer earth map tiling by Lemma $\ref{proposition-7}'$.  \label{a>1 b=2d 1}
	
	
	When $f=6k+2$ $(k\ge2)$, we have $\ccc=k\bbb$, $\aaa+\ddd=(2k+1)\bbb$ and  \[\text{AVC}\sub\{\aaa\ccc\ddd,\bbb\ccc^3,\aaa\bbb^k\ddd,\bbb^{k+1}\ccc^2,\bbb^{2k+1}\ccc,\bbb^{\frac{f}{2}}\}.\]  
		
		By AVC, we know $\aaa^2\cdots,\ddd^2\cdots$ and $\aaa\thin\ddd\cdots$ are not vertices. So we have the unique AAD for any  $\bbb^x\ccc^y=\thin^{\bbb}\ccc^{\ddd}\thin^{\ccc}\bbb^{\aaa}\thin\cdots\thin^{\bbb}\ccc^{\ddd}\thin^{\ccc}\bbb^{\aaa}\thin$. We will discuss all possible $\bbb$-vertices in any tiling as follows. 
	
	If $\bbb^{\frac{f}{2}}$ appears, the tiling is a $2$-layer earth map tiling by Lemma $\ref{proposition-7}'$. 

	If $\bbb^{2k+1}\ccc$ appears ($\bbb^{\frac{f}{2}}$ never appears), then $R(\ccc_2\cdots)=\bbb^{2k+1}$ in the first picture of Fig.\,\ref{fig 4-2} and this $\bbb^{2k+1}$  determines $2k+1$ time zones ($4k+2$ or $\frac{2f+2}{3}$ tiles).  Then $R(\aaa_1\ddd_4\cdots)=\bbb^{k}$ and this $\bbb^{k}$  determines $k$ time zones ($2k$ or $\frac{f-2}{3}$ tiles). We obtain a unique tiling $T(6k\aaa\ccc\ddd,2\aaa\bbb^{k}\ddd,2\bbb^{2k+1}\ccc)$ which can be viewed as the first flip modification of the $2$-layer earth map tilings.
	
	If $\bbb^{k+1}\ccc^2$ appears ($\bbb^{2k+1}\ccc,\bbb^{\frac{f}{2}}$ never appear), the tilings are shown in the second picture of Fig.\,\ref{fig 4-2}. Depending on the space between two flips, there are $\lfloor \frac{k+3}{2} \rfloor$ or $\lfloor \frac{f+16}{12} \rfloor$ different tilings with the same set of vertices. 
	
	If $\bbb\ccc^3$ appears ($\bbb^{k+1}\ccc^2,\bbb^{2k+1}\ccc,\bbb^{\frac{f}{2}}$ never appear), the tiling is shown in Fig.\,\ref{fig 4-4}. We obtain a unique tiling $T((6k-4)\aaa\ccc\ddd,2\bbb\ccc^3,6\aaa\bbb^{k}\ddd)$ which can be obtained by applying the first flip modification three times.   
	
	All of the above tilings belong to the third infinite sequence in Table \ref{Tab-1.2} after
	interchanging $\aaa\leftrightarrow\ddd$ and $\bbb\leftrightarrow\ccc$ to keep consistent the AVC for $2$-layer earth map tilings.  \label{a>1 b=2d 2}\label{a>1 b=2d 3} \label{a>1 b=2d 4}	
	
	\begin{figure}[htp]
		\centering
		\begin{tikzpicture}[>=latex,scale=0.48] 
			\foreach \a in {0,-6}
			{
				\begin{scope}[xshift=\a cm] 
					\draw (4,0)--(4,-9)
					(6,0)--(6,-9);
					\draw[line width=1.5] (6,-2)--(4,-7);
					\node at (5,0){\small $\bbb$};
					\node at (4.3,-2){\small $\ccc$};
					\node at (5.7,-1.9){\small $\aaa$};
					\node at (5.7,-3.4){\small $\ddd$};
					\node at (4.3,-5.6){\small $\ddd$};
					\node at (5.7,-7){\small $\ccc$};
					\node at (4.4,-7.1){\small $\aaa$};
					\node at (5,-9){\small $\bbb$};
				\end{scope}
			}
			
			\draw (0,0)--(0,-9)
			(4,0)--(4,-9)
			(0,-2)--(3,-2)--(3,-4)--(4,-7)--(2,-5)
			(0,-2)--(2,-4)
			(0,-2)--(1,-5)--(1,-7)--(4,-7);
			\draw[line width=1.5] (3,-2)--(4,-2)
			(2,-4)--(3,-4)
			(1,-5)--(2,-5)
			(0,-7)--(1,-7);
			
			\node at (2,0){\small $\ccc$}; 
			\node at (0.3,-1.6){\small $\bbb$};
			\node at (3.7,-1.6){\small $\ddd$};
			\node at (3,-1.6){\small $\aaa$};
			\node at (1,-2.4){\small $\bbb$};
			\node at (2.7,-2.4){\small $\ccc$};
			\node at (3.2,-2.4){\small $\ddd$};
			\node at (3.75,-2.4){\small $\aaa$};
			\node at (2.1,-3.6){\small $\aaa$};
			\node at (2.75,-3.6){\small $\ddd$};
			\node at (3.25,-4){\small $\ccc$};
			\node at (0.25,-3.8){\small $\bbb$};
			\node at (3.75,-5.2){\small $\bbb$};
			\node at (0.7,-5){\small $\ccc$};
			\node at (0.25,-6.6){\small $\aaa$};
			\node at (0.75,-6.6){\small $\ddd$};
			\node at (1.3,-5.4){\small $\ddd$};
			\node at (2,-5.4){\small $\aaa$};
			\node at (1.3,-6.6){\small $\ccc$};
			\node at (3,-6.6){\small $\bbb$};
			\node at (0.25,-7.4){\small $\ddd$};
			\node at (1,-7.4){\small $\aaa$};
			\node at (3.7,-7.4){\small $\bbb$};
			\node at (2,-9){\small $\ccc$};

			\fill (0.56,-3.21) circle (0.05);
			\fill (0.68,-3.13) circle (0.05);
			\fill (0.80,-3.05) circle (0.05);
			
			\fill (0.56+2.55,-3.21-2.6) circle (0.05);
			\fill (0.68+2.55,-3.13-2.6) circle (0.05);
			\fill (0.80+2.55,-3.05-2.6) circle (0.05);
			
			\fill (7-0.25,-5) circle (0.05);
			\fill (7.3-0.25,-5) circle (0.05);
			\fill (7.6-0.25,-5) circle (0.05);
			
			\node[draw,shape=circle, inner sep=0.5] at (-1,-1) {\small $1$};
			\node[draw,shape=circle, inner sep=0.5] at (2,-1) {\small $2$};
			\node[draw,shape=circle, inner sep=0.5] at (5,-1) {\small $3$};
			\node[draw,shape=circle, inner sep=0.5] at (-1,-8) {\small $4$};
			\node[draw,shape=circle, inner sep=0.5] at (2,-8) {\small $5$};
			\node[draw,shape=circle, inner sep=0.5] at (5,-8) {\small $6$};
			
			\fill (0,-2) circle (0.15); \fill (4,-7) circle (0.15);
			\fill (6,-7) circle (0.04);
			\node at (2,-10.5) {\small $T(6k\aaa\ccc\ddd,2\aaa\bbb^{k}\ddd,2\bbb^{2k+1}\ccc)$};
			
		\begin{scope}[xshift=10cm]
				\foreach \a in {0,1}
			{
				\begin{scope}[xshift=8*\a cm] 
					\draw (4,0)--(4,-9)
					(6,0)--(6,-9);
					\draw[line width=1.5] (6,-2)--(4,-7);
					\node at (5,0){\small $\bbb$};
					\node at (4.3,-2){\small $\ccc$};
					\node at (5.7,-1.9){\small $\aaa$};
					\node at (5.7,-3.4){\small $\ddd$};
					\node at (4.3,-5.6){\small $\ddd$};
					\node at (5.7,-7){\small $\ccc$};
					\node at (4.4,-7.1){\small $\aaa$};
					\node at (5,-9){\small $\bbb$};
				\end{scope}
			}
			\foreach \b in {0,1}
			{
				\begin{scope}[xshift=8*\b cm] 
					\draw (0,0)--(0,-9)
					(4,0)--(4,-9)
					(0,-2)--(3,-2)--(3,-4)--(4,-7)--(2,-5)
					(0,-2)--(2,-4)
					(0,-2)--(1,-5)--(1,-7)--(4,-7);
					\draw[line width=1.5] (3,-2)--(4,-2)
					(2,-4)--(3,-4)
					(1,-5)--(2,-5)
					(0,-7)--(1,-7);
					
					\node at (2,0){\small $\ccc$}; 
					\node at (0.3,-1.6){\small $\bbb$};
					\node at (3.7,-1.6){\small $\ddd$};
					\node at (3,-1.6){\small $\aaa$};
					\node at (1,-2.4){\small $\bbb$};
					\node at (2.7,-2.4){\small $\ccc$};
					\node at (3.2,-2.4){\small $\ddd$};
					\node at (3.75,-2.4){\small $\aaa$};
					\node at (2.1,-3.6){\small $\aaa$};
					\node at (2.75,-3.6){\small $\ddd$};
					\node at (3.25,-4){\small $\ccc$};
					\node at (0.25,-3.8){\small $\bbb$};
					\node at (3.75,-5.2){\small $\bbb$};
					\node at (0.7,-5){\small $\ccc$};
					\node at (0.25,-6.6){\small $\aaa$};
					\node at (0.75,-6.6){\small $\ddd$};
					\node at (1.3,-5.4){\small $\ddd$};
					\node at (2,-5.4){\small $\aaa$};
					\node at (1.3,-6.6){\small $\ccc$};
					\node at (3,-6.6){\small $\bbb$};
					\node at (0.25,-7.4){\small $\ddd$};
					\node at (1,-7.4){\small $\aaa$};
					\node at (3.7,-7.4){\small $\bbb$};
					\node at (2,-9){\small $\ccc$};

					\fill (0.56,-3.21) circle (0.05);
					\fill (0.68,-3.13) circle (0.05);
					\fill (0.80,-3.05) circle (0.05);
					
					\fill (0.56+2.55,-3.21-2.6) circle (0.05);
					\fill (0.68+2.55,-3.13-2.6) circle (0.05);
					\fill (0.80+2.55,-3.05-2.6) circle (0.05);
					
					\fill (7-0.25,-5) circle (0.05);
					\fill (7.3-0.25,-5) circle (0.05);
					\fill (7.6-0.25,-5) circle (0.05);
					\fill (0,-2) circle (0.15); \fill (4,-7) circle (0.15);
					
					\fill (6,-7) circle (0.04);
				\end{scope}
			}
			\node at (7,-10.5){\small $\{(6k-2)\aaa\ccc\ddd,4\aaa\bbb^{k}\ddd,2\bbb^{k+1}\ccc^2\}$};
		\end{scope}	
		\end{tikzpicture}
		\caption{Tilings by flipping once, or twice with different spacing.} \label{fig 4-2}
	\end{figure}
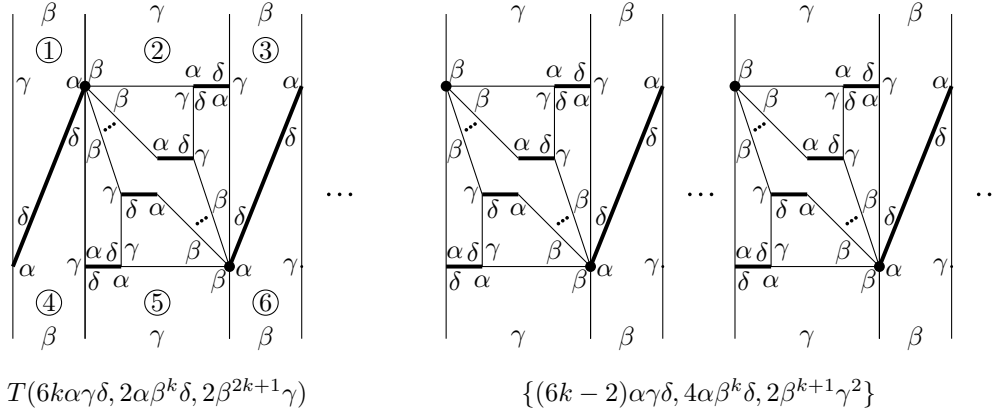

	\begin{figure}[htp]
		\centering
		\begin{tikzpicture}[>=latex,scale=0.52] 
			
			\begin{scope}[xshift=8 cm] 
				\draw (4,0)--(4,-9)
				(7,0)--(7,-9);
				\draw[line width=1.5] (7,-2)--(4,-7);
				\node at (5.5,0){\small $\bbb$};
				\node at (4.3,-2){\small $\ccc$};
				\node at (6.7,-1.9){\small $\aaa$};
				\node at (6.7,-3){\small $\ddd$};
				\node at (4.3,-6){\small $\ddd$};
				\node at (6.7,-7){\small $\ccc$};
				\node at (4.4,-7.1){\small $\aaa$};
				\node at (5.5,-9){\small $\bbb$};
			\end{scope}
			\foreach \b in {0,1,2}
			{
				\begin{scope}[xshift=4*\b cm] 
					\draw (0,0)--(0,-9)
					(4,0)--(4,-9)
					(0,-2)--(3,-2)--(3,-4)--(4,-7)--(2,-5)
					(0,-2)--(2,-4)
					(0,-2)--(1,-5)--(1,-7)--(4,-7);
					\draw[line width=1.5] (3,-2)--(4,-2)
					(2,-4)--(3,-4)
					(1,-5)--(2,-5)
					(0,-7)--(1,-7);
					
					\node at (2,0){\small $\ccc$}; 
					\node at (0.3,-1.6){\small $\bbb$};
					\node at (3.7,-1.6){\small $\ddd$};
					\node at (3,-1.6){\small $\aaa$};
					\node at (1,-2.4){\small $\bbb$};
					\node at (2.7,-2.4){\small $\ccc$};
					\node at (3.2,-2.4){\small $\ddd$};
					\node at (3.75,-2.4){\small $\aaa$};
					\node at (2.1,-3.6){\small $\aaa$};
					\node at (2.75,-3.6){\small $\ddd$};
					\node at (3.25,-4){\small $\ccc$};
					\node at (0.25,-3.8){\small $\bbb$};
					\node at (3.75,-5.2){\small $\bbb$};
					\node at (0.7,-5){\small $\ccc$};
					\node at (0.25,-6.6){\small $\aaa$};
					\node at (0.75,-6.6){\small $\ddd$};
					\node at (1.3,-5.4){\small $\ddd$};
					\node at (2,-5.4){\small $\aaa$};
					\node at (1.3,-6.6){\small $\ccc$};
					\node at (3,-6.6){\small $\bbb$};
					\node at (0.25,-7.4){\small $\ddd$};
					\node at (1,-7.4){\small $\aaa$};
					\node at (3.7,-7.4){\small $\bbb$};
					\node at (2,-9){\small $\ccc$};

					\fill (0.56,-3.21) circle (0.05);
					\fill (0.68,-3.13) circle (0.05);
					\fill (0.80,-3.05) circle (0.05);
					
					\fill (0.56+2.55,-3.21-2.6) circle (0.05);
					\fill (0.68+2.55,-3.13-2.6) circle (0.05);
					\fill (0.80+2.55,-3.05-2.6) circle (0.05);
					
					\fill (0,-2) circle (0.15); \fill (4,-7) circle (0.15);
				\end{scope}
			}
			
			
			\fill (15,-7) circle (0.05);
			
		\end{tikzpicture}
		\caption{$T((6k-4)\aaa\ccc\ddd,2\bbb\ccc^3,6\aaa\bbb^{k}\ddd)$ obtained by flipping $3$ times.} \label{fig 4-4}
	\end{figure}
	
	If $\text{AVC}\sub\{\aaa\ccc\ddd,\aaa\bbb^{k}\ddd\}$, there is no solution satisfying Balance Lemma \ref{balance}.

	\subsection*{Case $2$: $\{x_1, x_2\}=\{ x_3, x_4\}$}
	By Proposition \ref{symmetric}, $\bbb\neq\ccc$. So the only possibility is that  $x_1=x_3$ and $x_2=x_4$ in \eqref{4-9}. After an easy check of the seven choices in \eqref{4-9}, only the last one might hold: $2-\aaa-\frac\ccc2=\frac\ccc2$ and $\frac\ccc2=1-\ddd-\frac\bbb2$. Then we get $\aaa+\bbb+\ccc+\ddd=3>\frac83$, a contradiction.

	\begin{table*}[htp]                        
		\centering     		 
		\resizebox{\textwidth}{85mm}{\begin{tabular}{cccc|cccc|c}	 
				
				$\theta$&$\frac16$&$\frac{\theta}{2}$&$\frac{1}{2}-\frac{\theta}{2}$&$\aaa$&$\bbb$&$\ccc$&$\ddd$&$\aaa>1>\ccc>\bbb>\ddd$ (Lemma \ref{proposition-4}) \\
				\hline 
				$1-\aaa+\frac{\ccc}{2}$&$\frac{\bbb}{2}$&$\frac{\ccc}{2}$&$\ddd-\frac{\bbb}{2}$&$1-\frac{\theta}{2}$&$\frac13$&$\theta$&$\frac23-\frac{\theta}{2}$&$\aaa<1$ \\
				$-1+\aaa-\frac{\ccc}{2}$&$\frac{\bbb}{2}$&$\frac{\ccc}{2}$&$-\ddd+\frac{\bbb}{2}$&$1+\frac{3\theta}{2}$&$\frac13$&$\theta$&$-\frac13+\frac{\theta}{2}$&$\ddd<0$ \\ 
				$2-\aaa+\frac{\ccc}{2}$&$\frac{\bbb}{2}$&$\frac{\ccc}{2}$&$-\ddd+\frac{\bbb}{2}$&$2-\frac{\theta}{2}$&$\frac13$&$\theta$&$-\frac13+\frac{\theta}{2}$&$\ddd<0 $ \\
				$-1+\aaa+\frac{\ccc}{2}$&$\frac{\bbb}{2}$&$\frac{\ccc}{2}$&$\ddd+\frac{\bbb}{2}$&$1+\frac{\theta}{2}$&$\frac13$&$\theta$&$\frac13-\frac{\theta}{2}$&\textbf{$\surd$ \,\,\, Subcase} $1$\\
				$-1+\aaa+\frac{\ccc}{2}$&$\frac{\bbb}{2}$&$\frac{\ccc}{2}$&$1-\ddd-\frac{\bbb}{2}$&$1+\frac{\theta}{2}$&$\frac13$&$\theta$&$\frac13+\frac{\theta}{2}$& $\bbb<\ddd$   \\
				$2-\aaa-\frac{\ccc}{2}$&$\frac{\bbb}{2}$&$\frac{\ccc}{2}$&$\ddd+\frac{\bbb}{2}$&$2-\frac{3\theta}{2}$&$\frac13$&$\theta$&$\frac13-\frac{\theta}{2}$&\textbf{$\surd$ \,\,\, Subcase} $2$ \\
				$2-\aaa-\frac{\ccc}{2}$&$\frac{\bbb}{2}$&$\frac{\ccc}{2}$&$1-\ddd-\frac{\bbb}{2}$&$2-\frac{3\theta}{2}$&$\frac13$&$\theta$&$\frac13+\frac{\theta}{2}$& $\bbb<\ddd$   \\
				\multicolumn{9}{c}{ }\\
				$\theta$&$\frac16$&$\frac{1}{2}-\frac{\theta}{2}$&$\frac{\theta}{2}$&$\aaa$&$\bbb$&$\ccc$&$\ddd$&$\aaa>1>\ccc>\bbb>\ddd$ (Lemma \ref{proposition-4})  \\
				\hline 
				$1-\aaa+\frac{\ccc}{2}$&$\frac{\bbb}{2}$&$\frac{\ccc}{2}$&$\ddd-\frac{\bbb}{2}$&$\frac{3}{2}-\frac{3\theta}{2}$&$\frac13$&$1-\theta$&$\frac16+\frac{\theta}{2}$&\textbf{$\surd$ \,\,\, Subcase} $3$\\
				$-1+\aaa-\frac{\ccc}{2}$&$\frac{\bbb}{2}$&$\frac{\ccc}{2}$&$-\ddd+\frac{\bbb}{2}$&$\frac{3}{2}+\frac{\theta}{2}$&$\frac13$&$1-\theta$&$\frac16-\frac{\theta}{2}$&$\ddd>0\Rightarrow\theta<\frac13$ but\\
				&&&&&&&&$\aaa+\bbb+\ccc+\ddd\le\frac83\Rightarrow\theta\ge\frac13$\\
				$2-\aaa+\frac{\ccc}{2}$&$\frac{\bbb}{2}$&$\frac{\ccc}{2}$&$-\ddd+\frac{\bbb}{2}$&$\frac{5}{2}-\frac{3\theta}{2}$&$\frac13$&$1-\theta$&$\frac16-\frac{\theta}{2}$&$\aaa<2\Rightarrow\theta>\frac13$ but\\
				&&&&&&&&$\ddd>0\Rightarrow\theta<\frac13$\\
				$-1+\aaa+\frac{\ccc}{2}$&$\frac{\bbb}{2}$&$\frac{\ccc}{2}$&$\ddd+\frac{\bbb}{2}$&$\frac12+\frac{3\theta}{2}$&$\frac13$&$1-\theta$&$-\frac16+\frac{\theta}{2}$&\textbf{$\surd$ \,\,\, Subcase} $2$\\
				$-1+\aaa+\frac{\ccc}{2}$&$\frac{\bbb}{2}$&$\frac{\ccc}{2}$&$1-\ddd-\frac{\bbb}{2}$&$\frac12+\frac{3\theta}{2}$&$\frac13$&$1-\theta$&$\frac{5}{6}-\frac{\theta}{2}$&$\bbb<\ddd$ \\
				$2-\aaa-\frac{\ccc}{2}$&$\frac{\bbb}{2}$&$\frac{\ccc}{2}$&$\ddd+\frac{\bbb}{2}$&$\frac{3}{2}-\frac{\theta}{2}$&$\frac13$&$1-\theta$&$-\frac16+\frac{\theta}{2}$&\textbf{$\surd$ \,\,\, Subcase} $1$ \\
				$2-\aaa-\frac{\ccc}{2}$&$\frac{\bbb}{2}$&$\frac{\ccc}{2}$&$1-\ddd-\frac{\bbb}{2}$&$\frac{3}{2}-\frac{\theta}{2}$&$\frac13$&$1-\theta$&$\frac{5}{6}-\frac{\theta}{2}$&$\bbb<\ddd$  \\			
				\multicolumn{9}{c}{ }\\
				$\frac16$&$\theta$&$\frac{1}{2}-\frac{\theta}{2}$&$\frac{\theta}{2}$&$\aaa$&$\bbb$&$\ccc$&$\ddd$&$\aaa>1>\ccc>\bbb>\ddd$ (Lemma \ref{proposition-4}) \\
				\hline 
				$1-\aaa+\frac{\ccc}{2}$&$\frac{\bbb}{2}$&$\frac{\ccc}{2}$&$\ddd-\frac{\bbb}{2}$&$\frac{4}{3}-\frac{\theta}{2}$&$2\theta$&$1-\theta$&$\frac{3\theta}{2}$&\textbf{$\surd$ \,\,\, Subcase} $4$\\ 
				$-1+\aaa-\frac{\ccc}{2}$&$\frac{\bbb}{2}$&$\frac{\ccc}{2}$&$-\ddd+\frac{\bbb}{2}$&$\frac{5}{3}-\frac{\theta}{2}$&$2\theta$&$1-\theta$&$\frac{\theta}{2}$&$\aaa+\bbb+\ccc+\ddd>\frac83$  \\
				$2-\aaa+\frac{\ccc}{2}$&$\frac{\bbb}{2}$&$\frac{\ccc}{2}$&$-\ddd+\frac{\bbb}{2}$&$\frac73-\frac{\theta}{2}$&$2\theta$&$1-\theta$&$\frac{\theta}{2}$&$\aaa+\bbb+\ccc+\ddd>\frac83$  \\
				$-1+\aaa+\frac{\ccc}{2}$&$\frac{\bbb}{2}$&$\frac{\ccc}{2}$&$\ddd+\frac{\bbb}{2}$&$\frac23+\frac{\theta}{2}$&$2\theta$&$1-\theta$&$-\frac{\theta}{2}$&$\ddd<0$ \\
				$-1+\aaa+\frac{\ccc}{2}$&$\frac{\bbb}{2}$&$\frac{\ccc}{2}$&$1-\ddd-\frac{\bbb}{2}$&$\frac23+\frac{\theta}{2}$&$2\theta$&$1-\theta$&$1-\frac{3\theta}{2}$&$\aaa<1$ \\
				$2-\aaa-\frac{\ccc}{2}$&$\frac{\bbb}{2}$&$\frac{\ccc}{2}$&$\ddd+\frac{\bbb}{2}$&$\frac{4}{3}+\frac{\theta}{2}$&$2\theta$&$1-\theta$&$-\frac{\theta}{2}$&$\ddd<0$ \\
				$2-\aaa-\frac{\ccc}{2}$&$\frac{\bbb}{2}$&$\frac{\ccc}{2}$&$1-\ddd-\frac{\bbb}{2}$&$\frac{4}{3}+\frac{\theta}{2}$&$2\theta$&$1-\theta$&$1-\frac{3\theta}{2}$&$\aaa+\bbb+\ccc+\ddd=\frac{10}{3}>\frac{8}{3}$\\
				\multicolumn{9}{c}{ }\\
				$\frac{\theta}{2}$&$\frac{1}{2}-\frac{\theta}{2}$&$\theta$&$\frac{1}{6}$&$\aaa$&$\bbb$&$\ccc$&$\ddd$&$\aaa>1>\ccc>\bbb>\ddd$ (Lemma \ref{proposition-4}) \\
				\hline 
				$1-\aaa+\frac{\ccc}{2}$&$\frac{\bbb}{2}$&$\frac{\ccc}{2}$&$\ddd-\frac{\bbb}{2}$&$1+\frac{\theta}{2}$&$1-\theta$&$2\theta$&$\frac23-\frac{\theta}{2}$&$\aaa+\bbb+\ccc+\ddd>\frac83$ \\
				$-1+\aaa-\frac{\ccc}{2}$&$\frac{\bbb}{2}$&$\frac{\ccc}{2}$&$-\ddd+\frac{\bbb}{2}$&$1+\frac{3\theta}{2}$&$1-\theta$&$2\theta$&$\frac13-\frac{\theta}{2}$&$\bbb<\ccc\Rightarrow\theta>\frac13$ but\\
				&&&&&&&&$ \aaa+\bbb+\ccc+\ddd\le\frac83\Rightarrow\theta\le\frac16$\\
				$2-\aaa+\frac{\ccc}{2}$&$\frac{\bbb}{2}$&$\frac{\ccc}{2}$&$-\ddd+\frac{\bbb}{2}$&$2+\frac{\theta}{2}$&$1-\theta$&$2\theta$&$\frac13-\frac{\theta}{2}$&$\aaa>2$ \\
				$-1+\aaa+\frac{\ccc}{2}$&$\frac{\bbb}{2}$&$\frac{\ccc}{2}$&$\ddd+\frac{\bbb}{2}$&$1-\frac{\theta}{2}$&$1-\theta$&$2\theta$&$-\frac13+\frac{\theta}{2}$&$\aaa<1$ \\
				$-1+\aaa+\frac{\ccc}{2}$&$\frac{\bbb}{2}$&$\frac{\ccc}{2}$&$1-\ddd-\frac{\bbb}{2}$&$1-\frac{\theta}{2}$&$1-\theta$&$2\theta$&$\frac13+\frac{\theta}{2}$&$\aaa<1$ \\
				$2-\aaa-\frac{\ccc}{2}$&$\frac{\bbb}{2}$&$\frac{\ccc}{2}$&$\ddd+\frac{\bbb}{2}$&$2-\frac{3\theta}{2}$&$1-\theta$&$2\theta$&$-\frac13+\frac{\theta}{2}$&$\ddd<0$  \\
				$2-\aaa-\frac{\ccc}{2}$&$\frac{\bbb}{2}$&$\frac{\ccc}{2}$&$1-\ddd-\frac{\bbb}{2}$&$2-\frac{3\theta}{2}$&$1-\theta$&$2\theta$&$\frac13+\frac{\theta}{2}$&$\aaa+\bbb+\ccc+\ddd=\frac{10}{3}>\frac83$  \\			
				\multicolumn{9}{c}{ }\\
				$\frac{1}{2}-\frac{\theta}{2}$&$\frac{\theta}{2}$&$\theta$&$\frac{1}{6}$&$\aaa$&$\bbb$&$\ccc$&$\ddd$& $\aaa>1>\ccc>\bbb>\ddd$ (Lemma \ref{proposition-4})\\
				\hline 
				$1-\aaa+\frac{\ccc}{2}$&$\frac{\bbb}{2}$&$\frac{\ccc}{2}$&$\ddd-\frac{\bbb}{2}$&$\frac12+\frac{3\theta}{2}$&$\theta$&$2\theta$&$\frac16+\frac{\theta}{2}$&\textbf{$\surd$ \,\,\, Subcase} $5$\\ 
				$-1+\aaa-\frac{\ccc}{2}$&$\frac{\bbb}{2}$&$\frac{\ccc}{2}$&$-\ddd+\frac{\bbb}{2}$&$\frac{3}{2}+\frac{\theta}{2}$&$\theta$&$2\theta$&$-\frac16+\frac{\theta}{2}$&$\ddd>0\Rightarrow\theta>\frac13$ but \\
				&&&&&&&&$ \aaa+\bbb+\ccc+\ddd\le\frac83\Rightarrow\theta\le\frac13$\\
				$2-\aaa+\frac{\ccc}{2}$&$\frac{\bbb}{2}$&$\frac{\ccc}{2}$&$-\ddd+\frac{\bbb}{2}$&$\frac{3}{2}+\frac{3\theta}{2}$&$\theta$&$2\theta$&$-\frac16+\frac{\theta}{2}$&$\ddd>0\Rightarrow\theta>\frac13$ but\\
				&&&&&&&&$ \aaa+\bbb+\ccc+\ddd\le\frac83\Rightarrow\theta\le\frac{4}{15}$\\
				$-1+\aaa+\frac{\ccc}{2}$&$\frac{\bbb}{2}$&$\frac{\ccc}{2}$&$\ddd+\frac{\bbb}{2}$&$\frac{3}{2}-\frac{3\theta}{2}$&$\theta$&$2\theta$&$\frac16-\frac{\theta}{2}$&$\aaa>1\Rightarrow\theta<\frac13$ but\\
				&&&&&&&&$\aaa+\bbb+\ccc+\ddd>2\Rightarrow\theta>\frac13$\\
				$-1+\aaa+\frac{\ccc}{2}$&$\frac{\bbb}{2}$&$\frac{\ccc}{2}$&$1-\ddd-\frac{\bbb}{2}$&$\frac{3}{2}-\frac{3\theta}{2}$&$\theta$&$2\theta$&$\frac{5}{6}-\frac{\theta}{2}$&$\aaa>1\Rightarrow\theta<\frac13$ but\\
				&&&&&&&&$\ccc>\ddd\Rightarrow\theta>\frac13$\\
				$2-\aaa-\frac{\ccc}{2}$&$\frac{\bbb}{2}$&$\frac{\ccc}{2}$&$\ddd+\frac{\bbb}{2}$&$\frac{3}{2}-\frac{\theta}{2}$&$\theta$&$2\theta$&$\frac16-\frac{\theta}{2}$&\textbf{$\surd$ \,\,\, Subcase} $6$\\
				$2-\aaa-\frac{\ccc}{2}$&$\frac{\bbb}{2}$&$\frac{\ccc}{2}$&$1-\ddd-\frac{\bbb}{2}$&$\frac{3}{2}-\frac{\theta}{2}$&$\theta$&$2\theta$&$\frac{5}{6}-\frac{\theta}{2}$&$\bbb<\ddd$ \\

		\end{tabular}}
	\caption{Case $\{x_1, x_2\}=\{ \frac{1}{6},\theta\}$ and $\{x_3,x_4\}=\{\frac{\theta}{2},\frac{1}{2}-\frac{\theta}{2}\}$ or $\{x_3,x_4\}=\{\frac{1}{6},\theta\}$ and $\{x_1,x_2\}=\{\frac{\theta}{2},\frac{1}{2}-\frac{\theta}{2}\}$, for some $0<\theta\le\frac{1}{2}$}\label{Tab-4.1}        
	\end{table*}

	\subsection*{Case $3$: $\{x_1, x_2\}=\{ \frac{1}{6},\theta\}$ and $\{x_3,x_4\}=\{\frac{\theta}{2},\frac{1}{2}-\frac{\theta}{2}\}$, or $\{x_3,x_4\}$ $=$ $\{\frac{1}{6},\theta\}$ and $\{x_1,x_2\}=\{\frac{\theta}{2},\frac{1}{2}-\frac{\theta}{2}\}$, for some $0<\theta\le\frac{1}{2}$}
	
	In the seven choices of \eqref{4-9}, if $\frac{\ccc}{2}=\frac16$, then $\ccc=\frac13$, contradicting Lemma \ref{proposition-4}; if $\frac\ccc2=\frac \theta2$ and $\frac{\bbb}{2}=\theta$, then $\ccc<\bbb$, contradicting Lemma \ref{proposition-4}. Therefore, we only have $5\times 7=35$ options to consider.	
	It turns out $27$ of these options are dismissed by Lemma \ref{anglesum} and Lemma \ref{proposition-4}. We list the corresponding details in the right hand column of Table \ref{Tab-4.1}. 
	The remaining $8$ options are summarized as the following six subcases:
	
	$1$. $\aaa=1+\frac{\ccc}{2},\,\,\quad\bbb=\frac{1}{3},\quad\quad\,\,\quad \ddd=\frac{1}{3}-\frac{\ccc}{2},\quad\,\,\frac13<\ccc<\frac23$;
	
	$2$. $\aaa=2-\frac{3\ccc}{2}, \quad\bbb=\frac{1}{3},\quad\quad\quad\,\, \ddd=\frac{1}{3}-\frac{\ccc}{2},\quad\,\,\frac13<\ccc<\frac23$;
	
	$3$. $\aaa=\frac{3\ccc}{2}, \quad\quad\,\,\,\,\bbb=\frac{1}{3},\quad\quad\quad\,\, \ddd=\frac{2}{3}-\frac{\ccc}{2},\quad\,\,\frac23<\ccc<1$;

	$4$. $\aaa=\frac56+\frac{\ccc}{2},\quad\,\, \bbb=2-2\ccc,\quad \ddd=\frac{3}{2}-\frac{3\ccc}{2},\quad\frac56\le\ccc<1$;
	
	$5$. $\aaa=\frac12+\frac{3\ccc}{4},\quad \bbb=\frac{\ccc}{2},\quad\quad\quad\, \ddd=\frac{1}{6}+\frac{\ccc}{4},\quad\,\,\frac23<\ccc\le\frac45$;
	
	$6$. $\aaa=\frac32-\frac{\ccc}{4},\quad\,\, \bbb=\frac{\ccc}{2},\quad\quad\quad\, \ddd=\frac{1}{6}-\frac{\ccc}{4},\quad\,\,\frac13<\ccc<\frac23$.
	
	For the first, second and sixth subcases, we have $\bbb,\ccc<R(\aaa\ddd\cdots)$; for the third, fourth and fifth subcases, we have $\bbb<R(\aaa\ddd\cdots)<\ccc$. So neither $\aaa\bbb\ddd$ nor $\aaa\ccc\ddd$ is a vertex, contradicting  Lemma \ref{proposition-5'}.

	\subsection*{Case $4$: $\{ x_1,x_2,x_3,x_4\}$ are in Table \ref{Tab-3.1}.}
	
	\begin{table}[htp]                        
		\centering     	      
		\begin{minipage}{0.3\textwidth}
			\scriptsize{\begin{tabular}{cc} 
					
					$(\aaa,\bbb,\ccc,\ddd)$&$f$ \\
					\hline 
					(35,16,18,11)/30&6 \\     
					(35,16,18,3)/30&10 \\   
					(33,16,22,1)/30&10 \\   
					(19,7,9,1)/15&10 \\
					\hline 
					\multicolumn{1}{|c}{(41,10,16,3)/30}&\multicolumn{1}{c|}{12}  \\     
					\multicolumn{1}{|c}{(17,5,9,4)/15}&\multicolumn{1}{c|}{12} \\
					\hline
					(19,3,11,2)/15&12  \\   
					(67,12,50,11)/60&12  \\     
					(71,8,54,7)/60&12  \\   
					(41,8,18,3)/30&12  \\   
					
			\end{tabular}}
		\end{minipage}
		\raisebox{0em}{\begin{minipage}{0.3\textwidth}
				\scriptsize{\begin{tabular}{cc}
						$(\aaa,\bbb,\ccc,\ddd)$&$f$ \\
						\hline
						(55,16,18,7)/42&14  \\    
						(49,16,30,1)/42&14  \\
						\hline    
						\multicolumn{1}{|c}{(43,6,16,1)/30}&\multicolumn{1}{c|}{20} \\ 
						\hline  
						(43,4,18,1)/30&20  \\   
						(83,16,18,13)/60&24 \\
						(71,16,42,1)/60&24   \\   
						(23,3,5,1)/15&30     \\ 
						(41,8,10,5)/30&30   \\  
						(37,8,18,1)/30&30    \\  
						(67,16,42,1)/60&40   \\   
				\end{tabular}}
		\end{minipage}}
	\raisebox{0.4em}{\begin{minipage}{0.3\textwidth}
			\scriptsize{\begin{tabular}{cc}
				$(\aaa,\bbb,\ccc,\ddd)$&$f$ \\
				\hline				 
				(79,16,18,13)/60&40  \\    
				(43,6,8,5)/30&60     \\ 
				(39,8,10,5)/30&60    \\  
				(35,8,18,1)/30&60    \\  
				(49,4,6,3)/30&60    \\ 
				(39,6,16,1)/30&60    \\  
				(47,4,10,1)/30&60    \\  
				(77,10,36,1)/60&60  \\   
				(59,6,20,1)/42&84  \\
			\end{tabular}}
	\end{minipage}}
    \caption{ $29$ solutions induced from Table \ref{Tab-3.1}.}\label{Tab-4.3}  
	\end{table}

	\label{discrete-6} \label{discrete-8} \label{discrete-10}	
	There are $8\times7\times15=840$ subcases to consider, but most are ruled out by violating $2>\aaa>1>\ccc>\bbb>\ddd>0$, $\ccc>\frac13$, $\ddd<\frac12$ or $f$ being even integer. Such computations can be carried out efficiently by any spreadsheet program. 
	Only $29$ subcases are left in Table \ref{Tab-4.3}. But $26$ of them are ruled out by Lemma \ref{proposition-5'}: there is neither $\aaa\bbb\ddd$ nor $\aaa\ccc\ddd$. There are only three subcases left: $(17,5,9,4)/15$, $(41,10,16,3)/30$, $(43,6,16,1)/30$. They all imply
	$\aaa\cdots=\aaa\ccc\ddd$ by the angle values and Parity Lemma. There is only one solution satisfying Balance Lemma \ref{balance}: $\{f\aaa\ccc\ddd,2\bbb^{\frac{f}{2}}\}$, and it gives three $2$-layer earth map tilings in Table \ref{Tab-1.1} after interchanging $\aaa\leftrightarrow\ddd$ and $\bbb\leftrightarrow\ccc$ by Lemma \ref{proposition-7}.

	

	\section{Degenerate case $\aaa=1$}
	\label{sec-degenerate-a}	
	
If $\aaa=1$, the quadrilateral degenerates to an isosceles triangle in Fig.\,\ref{fig a-pi}. 
	
	\begin{figure}[htp]
		\centering
		\begin{tikzpicture}[>=latex,scale=0.75]

			\draw (0,0)--(2,2)--(4,0)
			(0,0)--(2.5,0);
			\draw[line width=1.5] (2.5,0)--(4,0);
			\node at (0.6,0.25){\small $\bbb$};
			\node at (2,1.6){\small $\ccc$};
			\node at (3.4,0.25){\small $\ddd$};
			\fill (2.5,0) circle (0.05);
			\node at (1,1.3){\small $a$};  	\node at (3,1.3){\small $a$};		\node at (2,-0.3){\small $a+b$};  		
		\begin{scope}[xshift=3 cm]	
			\draw (5.26,1.8)--(6.8,0.2);
			\draw[line width=1.5] (6.73,1.8) arc (18:-13:3); 
			\draw[line width=1.5] (5.26,1.8) arc (108+55:139+55:3); 
			\node at (6,2.6){\small $\bbb$}; \node at (6,-0.8){\small $\bbb$};\node at (6.6,0){\small $\ccc$};
			\node at (6.7,0.8){\small $\ddd$}; \node at (5.6,1.9){\small $\ccc$}; \node at (5.4,1.2){\small $\ddd$};		
			\draw (6,3) arc (45:-45:3); \draw (6,3) arc (135:225:3);
			
			\fill (6.73,1.8) circle (0.05); \fill (5.21,0.2) circle (0.05);
		\end{scope}
		\end{tikzpicture}
		\caption{Degenerate case $\aaa=1$ and the subcase $\ccc+\ddd=1$.} \label{fig a-pi}
	\end{figure}
Then $\bbb=\ddd$, and Lemma \ref{geometry1} implies $\bbb<\ccc$. By Lemma \ref{proposition-5'}, exactly one of $\aaa\bbb\ddd$ or $\aaa\ccc\ddd$ must be a vertex in any spherical tiling by congruent such quadrilaterals. 
	
    
    \subsection*{Subcase $\aaa\bbb\ddd$ is a vertex} \label{discrete-1}
     By Lemma \ref{anglesum}, we get $\bbb=\ddd=\frac12,\ccc=\frac{4}{f}$. Then  $\bbb<\ccc$ implies $f=6$ and $\ccc=\frac23$. So the $\text{AVC}=\{\aaa\bbb\ddd,\ccc^3\}$, and it gives a $2$-layer earth map tiling by Lemma \ref{proposition-7}. This is Case $(6,3,4,3)/6$ in Table \ref{Tab-1.1}. 
	
	\subsection*{Subcase $\aaa\ccc\ddd$ is a vertex}
	 By Lemma \ref{anglesum}, we get $\bbb=\ddd=\frac{4}{f},\ccc=1-\frac{4}{f}$. Then  $\bbb<\ccc$ implies $f>8$ and $\ccc>\frac12$. By the angle values and Parity Lemma, we get $\text{AVC}\sub\{\aaa\ccc\ddd,\ccc^3,\bbb^2\ccc^2,\aaa\bbb^{\frac{f-4}{4}}\ddd,\bbb^{\frac{f+4}{4}}\ccc,\bbb^{\frac{f}{2}}\}$. 
	  
	 When $f=4k+2(k\ge2)$, we have $\text{AVC}\sub\{\aaa\ccc\ddd,\ccc^3,\bbb^2\ccc^2,\bbb^{\frac{f}{2}}\}$, and the only solution satisfying Balance Lemma \ref{balance} is $\{ f\aaa\ccc\ddd,2\bbb^{\frac{f}{2}}\}$  which gives a $2$-layer earth map tiling by Lemma $\ref{proposition-7}'$.   
	 
	 
	 When $f=4k$ $(k\ge3)$, we have $\ccc=(k-1)\bbb$, $\aaa+\ddd=(k+1)\bbb$ and  \[\text{AVC}\sub\{\aaa\ccc\ddd,\ccc^3,\bbb^2\ccc^2,\aaa\bbb^{k-1}\ddd, \bbb^{k+1}\ccc,\bbb^{\frac{f}{2}}\}.\]
	 Trying out all possible $\bbb$-vertices as the previous section, there are always $4$ tilings as shown in Fig.\,$\ref{1-1}'$, the second picture of Fig.\,\ref{fig a=1 2} (flip once), Fig.\,\ref{fig a=1 3} (flip twice with different spacing). \label{a=1} 
	 Only when $f=12$, we can apply the first flip modification in Fig.\,\ref{flip1} (after interchanging $\aaa\leftrightarrow\ddd$ and $\bbb\leftrightarrow\ccc$) three times, as shown in the first picture of Fig.\,\ref{fig a=1 2}. This is because $3(k-1)>2k$ when $k\ge4$.
	 
	 All above tilings belong to the first infinite sequence in Table \ref{Tab-1.2}.

	\begin{figure}[htp]
		\centering
		\begin{tikzpicture}[>=latex,scale=0.48] 
			\foreach \a in {0,1}
			{
				\begin{scope}[xshift=3*\a cm] 
					\draw (4,0)--(4,-9)
					(7,0)--(7,-9);
					\draw[line width=1.5] (7,-2)--(4,-7);
					\node at (5.5,0){\small $\bbb$};
					\node at (4.3,-2){\small $\ccc$};
					\node at (6.7,-1.9){\small $\aaa$};
					\node at (6.7,-3){\small $\ddd$};
					\node at (4.3,-6){\small $\ddd$};
					\node at (6.7,-7){\small $\ccc$};
					\node at (4.4,-7.1){\small $\aaa$};
					\node at (5.5,-9){\small $\bbb$};
				\end{scope}
			}
			
			\draw (0,0)--(0,-9)
			(4,0)--(4,-9)
			(0,-2)--(3,-2)--(3,-4)--(4,-7)--(2,-5)
			(0,-2)--(2,-4)
			(0,-2)--(1,-5)--(1,-7)--(4,-7);
			\draw[line width=1.5] (3,-2)--(4,-2)
			(2,-4)--(3,-4)
			(1,-5)--(2,-5)
			(0,-7)--(1,-7);
			
			\node at (2,0){\small $\ccc$}; 
			\node at (0.3,-1.6){\small $\bbb$};
			\node at (3.7,-1.6){\small $\ddd$};
			\node at (3,-1.6){\small $\aaa$};
			\node at (1,-2.4){\small $\bbb$};
			\node at (2.7,-2.4){\small $\ccc$};
			\node at (3.2,-2.4){\small $\ddd$};
			\node at (3.75,-2.4){\small $\aaa$};
			\node at (2.1,-3.6){\small $\aaa$};
			\node at (2.75,-3.6){\small $\ddd$};
			\node at (3.25,-4){\small $\ccc$};
			\node at (0.25,-3.8){\small $\bbb$};
			\node at (3.75,-5.2){\small $\bbb$};
			\node at (0.7,-5){\small $\ccc$};
			\node at (0.25,-6.6){\small $\aaa$};
			\node at (0.75,-6.6){\small $\ddd$};
			\node at (1.3,-5.4){\small $\ddd$};
			\node at (2,-5.4){\small $\aaa$};
			\node at (1.3,-6.6){\small $\ccc$};
			\node at (3,-6.6){\small $\bbb$};
			\node at (0.25,-7.4){\small $\ddd$};
			\node at (1,-7.4){\small $\aaa$};
			\node at (3.7,-7.4){\small $\bbb$};
			\node at (2,-9){\small $\ccc$};

			\fill (0.56,-3.21) circle (0.05);
			\fill (0.68,-3.13) circle (0.05);
			\fill (0.80,-3.05) circle (0.05);
			
			\fill (0.56+2.55,-3.21-2.6) circle (0.05);
			\fill (0.68+2.55,-3.13-2.6) circle (0.05);
			\fill (0.80+2.55,-3.05-2.6) circle (0.05);
			
			\fill (11,-5) circle (0.05);
			\fill (11.3,-5) circle (0.05);
			\fill (11.6,-5) circle (0.05);

			\fill (0,-2) circle (0.15); \fill (4,-7) circle (0.15);
			
			\fill (10,-7) circle (0.04);
			
		\node at (5.1,-10.6){\small $T((4k-2)\aaa\ccc\ddd,2\aaa\bbb^{k-1}\ddd,2\bbb^{k+1}\ccc)$};
	    \begin{scope}[xshift=-15.5cm] 
	    	
	    	\foreach \b in {0,1,2}
	    	{
	    		\begin{scope}[xshift=4*\b cm] 
	    			\draw (0,0)--(0,-7-2)
	    			(4,0)--(4,-7-2)
	    			(0,-2-1)--(2,-2-1)--(2,-5-1)--(4,-5-1);
	    			\draw[line width=1.5] (2,-2-1)--(4,-2-1)
	    			(0,-5-1)--(2,-5-1);
	    			
	    			\node at (2,1-1){\small $\ccc$}; 
	    			\node at (0.4,-1.4-1){\small $\bbb$};
	    			\node at (3.6,-1.4-1){\small $\ddd$};
	    			\node at (2,-1.4-1){\small $\aaa$};
	    			\node at (0.4,-2.6-1){\small $\bbb$};
	    			\node at (1.6,-2.6-1){\small $\ccc$};
	    			\node at (2.4,-2.6-1){\small $\ddd$};
	    			\node at (3.6,-2.6-1){\small $\aaa$};
	    			\node at (0.4,-4.4-1){\small $\aaa$};
	    			\node at (1.6,-4.4-1){\small $\ddd$};
	    			\node at (2.4,-4.4-1){\small $\ccc$};
	    			\node at (3.6,-4.4-1){\small $\bbb$};
	    			\node at (0.4,-5.6-1){\small $\ddd$};
	    			\node at (2,-5.6-1){\small $\aaa$};
	    			\node at (3.6,-5.6-1){\small $\bbb$};
	    			\node at (2,-8-1){\small $\ccc$};

	    		\end{scope}
	    	}
	    	\node at (6,-9.6-1){\small $T(6\aaa\ccc\ddd,2\ccc^3,6\aaa\bbb^2\ddd)$};
	    \end{scope}
       \end{tikzpicture}
		\caption{Two degenerate $a^3b$-tilings.} \label{fig a=1 2}
	\end{figure}
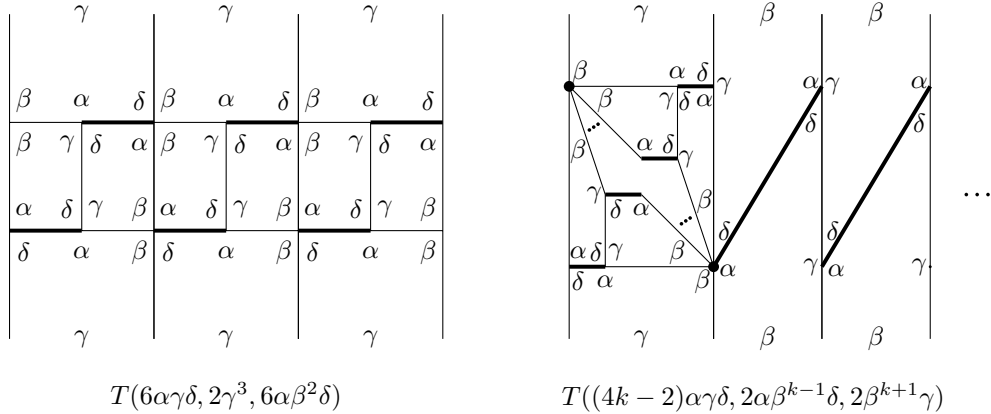	
	\begin{figure}[htp]
		\centering
		\begin{tikzpicture}[>=latex,scale=0.48] 
			\foreach \a in {0,1}
			{
				\begin{scope}[xshift=2*\a cm] 
					\draw (4,0)--(4,-9)
					(6,0)--(6,-9);
					\draw[line width=1.5] (6,-2)--(4,-7);
					\node at (5,0){\small $\bbb$};
					\node at (4.3,-2){\small $\ccc$};
					\node at (5.7,-1.9){\small $\aaa$};
					\node at (5.7,-3.4){\small $\ddd$};
					\node at (4.3,-5.6){\small $\ddd$};
					\node at (5.7,-7){\small $\ccc$};
					\node at (4.4,-7.1){\small $\aaa$};
					\node at (5,-9){\small $\bbb$};
				\end{scope}
			}
			\foreach \b in {0,1}
			{
				\begin{scope}[xshift=8*\b cm] 
					\draw (0,0)--(0,-9)
					(4,0)--(4,-9)
					(0,-2)--(3,-2)--(3,-4)--(4,-7)--(2,-5)
					(0,-2)--(2,-4)
					(0,-2)--(1,-5)--(1,-7)--(4,-7);
					\draw[line width=1.5] (3,-2)--(4,-2)
					(2,-4)--(3,-4)
					(1,-5)--(2,-5)
					(0,-7)--(1,-7);
					
					\node at (2,0){\small $\ccc$}; 
					\node at (0.3,-1.6){\small $\bbb$};
					\node at (3.7,-1.6){\small $\ddd$};
					\node at (3,-1.6){\small $\aaa$};
					\node at (1,-2.4){\small $\bbb$};
					\node at (2.7,-2.4){\small $\ccc$};
					\node at (3.2,-2.4){\small $\ddd$};
					\node at (3.75,-2.4){\small $\aaa$};
					\node at (2.1,-3.6){\small $\aaa$};
					\node at (2.75,-3.6){\small $\ddd$};
					\node at (3.25,-4){\small $\ccc$};
					\node at (0.25,-3.8){\small $\bbb$};
					\node at (3.75,-5.2){\small $\bbb$};
					\node at (0.7,-5){\small $\ccc$};
					\node at (0.25,-6.6){\small $\aaa$};
					\node at (0.75,-6.6){\small $\ddd$};
					\node at (1.3,-5.4){\small $\ddd$};
					\node at (2,-5.4){\small $\aaa$};
					\node at (1.3,-6.6){\small $\ccc$};
					\node at (3,-6.6){\small $\bbb$};
					\node at (0.25,-7.4){\small $\ddd$};
					\node at (1,-7.4){\small $\aaa$};
					\node at (3.7,-7.4){\small $\bbb$};
					\node at (2,-9){\small $\ccc$};

					\fill (0.56,-3.21) circle (0.05);
					\fill (0.68,-3.13) circle (0.05);
					\fill (0.80,-3.05) circle (0.05);
					
					\fill (0.56+2.55,-3.21-2.6) circle (0.05);
					\fill (0.68+2.55,-3.13-2.6) circle (0.05);
					\fill (0.80+2.55,-3.05-2.6) circle (0.05);

					\fill (0,-2) circle (0.15); \fill (4,-7) circle (0.15);

				\end{scope}
			}
			

	    \begin{scope}[xshift=15cm]
	    	\foreach \a in {0,1}
	    	{
	    		\begin{scope}[xshift=6*\a cm] 
	    			\draw (4,0)--(4,-9)
	    			(6,0)--(6,-9);
	    			\draw[line width=1.5] (6,-2)--(4,-7);
	    			\node at (5,0){\small $\bbb$};
	    			\node at (4.3,-2){\small $\ccc$};
	    			\node at (5.7,-1.9){\small $\aaa$};
	    			\node at (5.7,-3.4){\small $\ddd$};
	    			\node at (4.3,-5.6){\small $\ddd$};
	    			\node at (5.7,-7){\small $\ccc$};
	    			\node at (4.4,-7.1){\small $\aaa$};
	    			\node at (5,-9){\small $\bbb$};
	    		\end{scope}
	    	}
	    	\foreach \b in {0,1}
	    	{
	    		\begin{scope}[xshift=6*\b cm] 
	    			\draw (0,0)--(0,-9)
	    			(4,0)--(4,-9)
	    			(0,-2)--(3,-2)--(3,-4)--(4,-7)--(2,-5)
	    			(0,-2)--(2,-4)
	    			(0,-2)--(1,-5)--(1,-7)--(4,-7);
	    			\draw[line width=1.5] (3,-2)--(4,-2)
	    			(2,-4)--(3,-4)
	    			(1,-5)--(2,-5)
	    			(0,-7)--(1,-7);
	    			
	    			\node at (2,0){\small $\ccc$}; 
	    			\node at (0.3,-1.6){\small $\bbb$};
	    			\node at (3.7,-1.6){\small $\ddd$};
	    			\node at (3,-1.6){\small $\aaa$};
	    			\node at (1,-2.4){\small $\bbb$};
	    			\node at (2.7,-2.4){\small $\ccc$};
	    			\node at (3.2,-2.4){\small $\ddd$};
	    			\node at (3.75,-2.4){\small $\aaa$};
	    			\node at (2.1,-3.6){\small $\aaa$};
	    			\node at (2.75,-3.6){\small $\ddd$};
	    			\node at (3.25,-4){\small $\ccc$};
	    			\node at (0.25,-3.8){\small $\bbb$};
	    			\node at (3.75,-5.2){\small $\bbb$};
	    			\node at (0.7,-5){\small $\ccc$};
	    			\node at (0.25,-6.6){\small $\aaa$};
	    			\node at (0.75,-6.6){\small $\ddd$};
	    			\node at (1.3,-5.4){\small $\ddd$};
	    			\node at (2,-5.4){\small $\aaa$};
	    			\node at (1.3,-6.6){\small $\ccc$};
	    			\node at (3,-6.6){\small $\bbb$};
	    			\node at (0.25,-7.4){\small $\ddd$};
	    			\node at (1,-7.4){\small $\aaa$};
	    			\node at (3.7,-7.4){\small $\bbb$};
	    			\node at (2,-9){\small $\ccc$};

	    			\fill (0.56,-3.21) circle (0.05);
	    			\fill (0.68,-3.13) circle (0.05);
	    			\fill (0.80,-3.05) circle (0.05);
	    			
	    			\fill (0.56+2.55,-3.21-2.6) circle (0.05);
	    			\fill (0.68+2.55,-3.13-2.6) circle (0.05);
	    			\fill (0.80+2.55,-3.05-2.6) circle (0.05);

	    			\fill (0,-2) circle (0.15); \fill (4,-7) circle (0.15);

	    		\end{scope}
	    	}
	    	
	    	\fill (7,-7) circle (0.04);
	    	
	    \end{scope}
	    \end{tikzpicture}
		\caption{Two tilings for $\{(4k-4)\aaa\ccc\ddd,2\bbb^2\ccc^2,4\aaa\bbb^{k-1}\ddd\}$.} \label{fig a=1 3}
	\end{figure}

	

	\section{Concave case $\bbb>1$}
	\label{sec-concave-b}
		
    The quadrilateral with $\bbb>1$ is shown in Fig.\,\ref{quadrilateral-1}. We first prove some basic facts. Recall that Lemma $\ref{geometry4}'$ implies $\ccc+2\ddd>1$.
	
	\begin{figure}[htp]
		\centering
		\begin{tikzpicture}[>=latex,scale=0.4] 
			\draw (-4,5)--(-6,-1)--(0,0)--(6,-1);

			\draw[line width=1.5] (-4,5)--(6,-1);
			
			\draw[dotted] (-6,-1)--(6,-1)
			(0,0)--(-4,5);
			
			\node at (-6,-1.6){\small $C$};  \node at (5.8,-1.6){\small $A$};  \node at (0,-0.55){\small $B$};
			\node at (-4.4,5.5){\small $D$};   \node at (0,-1.7){\small $x$};    \node at (-2.5,2.2){\small $y$};

			\node at (-0.9,0.35){\small $\psi$}; \node at (4.2,-0.35){\small $\aaa$};
			\node at (-5.3,-0.4){\small $\ccc$}; \node at (-3.8,3.9){\small $\psi$};
			
			\node at (3,-0.8){\small $\phi$};  \node at (-3,-0.8){\small $\phi$};

			\node at (-5.5,2){\small $a$}; \node at (2.1,2){\small $b$}; \node at (-2,0.1){\small $a$};\node at (2,0.1){\small $a$};

		\end{tikzpicture}
		\caption{$a^3b$-quadrilateral with $\bbb>1,\aaa,\ccc,\ddd<1$.} \label{quadrilateral-1}
	\end{figure}
	
	\begin{lemma}\label{proposition-2}
	In an $a^3b$-tiling with $\bbb>1$, we have $a<b$, $\aaa<\ccc,\ddd$ and $\aaa<\frac12$. 
	\end{lemma}
	
	\begin{proof}
		 
		 In Fig.\,\ref{quadrilateral-1}, by $\bbb>\ddd$, $\angle ABD=\bbb-\psi>\ddd-\psi=\angle ADB$. This implies $a<b$. Then $\angle CAD<\angle ACD$, i.e. $\aaa+\phi<\ccc+\phi$. So $\aaa<\ccc$. By Lemma $\ref{geometry1}'$, $\aaa<\ddd$. 	
		 If $\aaa\ge\frac12$, then $\ccc,\ddd>\frac12$, and there is no $\bbb\cdots$ vertex. So $\aaa<\frac12$.		
	\end{proof}

	\begin{lemma}\label{proposition-3}
		In an $a^3b$-tiling with $\bbb>1$, $\bbb\ddd\cdots$ is a vertex and $\bbb+\ddd<2$.
	\end{lemma}
	
	\begin{proof}
		If $\bbb\ddd\cdots$ is not a vertex, by $\bbb>1$ and Parity Lemma, we get $\bbb\cdots=\aaa^x\bbb,\aaa^y\bbb\ccc^z,\bbb\ccc^w(x,y,w\ge2,z\ge1)$. Then $\#\aaa+\#\ccc\ge2\#\bbb=2f$, and there is only one solution satisfying Balance Lemma \ref{balance}: $\{\frac f2\,\aaa^2\bbb,\frac f2\, \bbb\ccc^2,\frac{f}{k}\,\ddd^k\}$. But this implies $\aaa=\ccc$, contradicting Lemma $\ref{geometry3}'$. Therefore, $\bbb\ddd\cdots$ is a vertex.						
	\end{proof}
	To find rational $a^3b$-quadrilaterals by solving \eqref{4-7} or \eqref{4-8} via Proposition \ref{proposition-6}, we have to transform $\aaa-\frac{\ccc}{2},\frac\bbb2,\frac\ccc2,\ddd-\frac\bbb2,\aaa+\frac\ccc2,-\ddd-\frac\bbb2$ to the range $[0,\frac12]$.
    For \eqref{4-7}, by Lemma \ref{proposition-2}, we have $-\frac12<\aaa-\frac{\ccc}{2}<\frac12,\frac12<\frac{\bbb}{2}<1,0<\frac{\ccc}{2}<\frac12$ and $-1<\ddd-\frac{\bbb}{2}<\frac12$.  
    For \eqref{4-8}, by Lemma \ref{proposition-2}, we have $0<\aaa+\frac{\ccc}{2}<1,\frac12<\frac\bbb2<1,0<\ccc<\frac12$, which implies $\sin(\ddd+\frac\bbb2)<0$. By $\frac12<\ddd+\frac{\bbb}{2}<2$, we get $1<\ddd+\frac{\bbb}{2}<2$. 	
	If $\frac32\le\ddd+\frac{\bbb}{2}<2$, we get $\bbb+\ddd>2$, contradicting Lemma \ref{proposition-3}.
	So for the equation \eqref{4-8}, we have $0<\aaa+\frac{\ccc}{2}<1,\frac12<\frac{\bbb}{2}<1,0<\frac{\ccc}{2}<\frac12$ and $1<\ddd+\frac{\bbb}{2}<\frac32$. 				
    Thus we have to consider the following five choices:
	\begin{small}
		\begin{equation}
			\{ x_1, x_2, x_3, x_4\}=\begin{cases}  
				\{\aaa-\frac{\ccc}{2},1-\frac{\bbb}{2},\frac{\ccc}{2},\ddd-\frac{\bbb}{2}\}, \\ 
				\{-\aaa+\frac{\ccc}{2},1-\frac{\bbb}{2},\frac{\ccc}{2},1+\ddd-\frac{\bbb}{2}\}, \\
				\{-\aaa+\frac{\ccc}{2},1-\frac{\bbb}{2},\frac{\ccc}{2},-\ddd+\frac{\bbb}{2}\}, \\
				\{\aaa+\frac{\ccc}{2},1-\frac{\bbb}{2},\frac{\ccc}{2},-1+\ddd+\frac{\bbb}{2}\}, \\
				\{1-\aaa-\frac{\ccc}{2},1-\frac{\bbb}{2},\frac{\ccc}{2},-1+\ddd+\frac{\bbb}{2}\}. \\
			\end{cases} \label{3-12}
		\end{equation}
	\end{small}

    We will match these choices with four cases of solutions in Proposition \ref{proposition-6} as follows. 

	\subsection*{Case $1$: $x_1  x_2=x_3 x_4=0$}			
	    By $-1<\ddd-\frac{\bbb}{2}<\frac12$, $-\frac12<\aaa-\frac{\ccc}{2}<\frac12$ and $1<\ddd+\frac\bbb2<\frac32$, the only solution of $x_1  x_2=x_3 x_4=0$ for \eqref{3-12} comes from $\aaa=\frac{\ccc}{2},\ddd=\frac{\bbb}{2}$. Then we get $3\aaa+\bbb+\ddd>2$. By $R(\bbb\ddd\cdots)<3\aaa,\bbb,\ddd$, Parity Lemma and Lemma \ref{proposition-3}, we deduce that $\aaa\bbb\ddd$ is a vertex. This implies $\aaa=\frac{2}{f},\bbb=\frac43-\frac{4}{3f},\ccc=\frac{4}{f},\ddd=\frac23-\frac{2}{3f}$. 
        By the angle values and Parity Lemma, we get $\text{AVC}\sub\{\aaa\bbb\ddd,\aaa\ddd^3,\aaa^x\bbb\ccc^{\frac{f-3x+2}{6}},\aaa^y\ccc^{\frac{f-3y+2}{6}}\ddd^2,\aaa^z\ccc^{\frac{2f-3z+1}{6}}\ddd,\aaa^w\ccc^{\frac{f-w}{2}}\}$.  Then we know there is no $\bbb^2\cdots$ vertex, which further implies that (by AAD)  \[\text{AVC}\sub\{\aaa\bbb\ddd,\aaa\ddd^3,\aaa^2\bbb\ccc^{\frac{f-4}{6}},\aaa^2\ccc^{\frac{f-4}{6}}\ddd^2,\bbb\ccc^{\frac{f+2}{6}} ,\ccc^{\frac{f+2}{6}}\ddd^2,\aaa\ccc^{\frac{f-1}{3}}\ddd,\ccc^{\frac{f}{2}}\}.\]
        
        When $f=6k$ or $6k+2$ $(k\ge1)$, we have $\text{AVC}\sub\{\aaa\bbb\ddd,\aaa\ddd^3,\ccc^{\frac{f}{2}}\}$, and the only solution satisfying Balance Lemma \ref{balance} is $\{ f\aaa\bbb\ddd,2\ccc^{\frac{f}{2}}\}$  which gives a $2$-layer earth map tiling by Lemma \ref{proposition-7}.   \label{b=2d-1}
        
        
        When $f=6k+4$ $(k\ge1)$, we have $\bbb=(2k+1)\ccc$, $\aaa+\ddd=(k+1)\ccc$ and  \[\text{AVC}\sub\{\aaa\bbb\ddd,\aaa\ddd^3,\aaa^2\bbb\ccc^{k} ,\aaa^2\ccc^{k}\ddd^2,\bbb\ccc^{k+1},\ccc^{k+1}\ddd^2,
        \aaa\ccc^{2k+1}\ddd,\ccc^{\frac{f}{2}}\}.\]

		We will discuss all possible $\ccc$-vertices in any tiling as follows.  Whenever $\ccc^{\frac{f}{2}}$ is a vertex, the tiling must be a $2$-layer earth map tiling by Lemma \ref{proposition-7}. If $\ccc^{\frac{f}{2}}$ never appears, we have the following subcases. 
		
		\subsubsection*{Subcase $\aaa^2\bbb\ccc^{k}$ appears ($\ccc^{\frac{f}{2}}$ never appears)} \label{b=2d-2}
			
			By AVC, $\bbb^2\cdots$ is never a vertex.
			Then $\aaa^2\bbb\ccc^{k}$ has only two possible AAD.   
			In Fig.\,\ref{b=2d 1}, $\aaa^2\bbb\ccc^{k}=\thin^{\bbb}\ccc^{\ddd}\thin^{\bbb}\aaa^{\ddd}\thick^{\ddd}\aaa^{\bbb}\thin^{\ddd}\ccc^{\bbb}\thin\cdots$ determines $T_1,T_2,T_3,T_4$. Then $\bbb_2\ddd_1\cdots=\aaa_5\bbb_2\ddd_1$ determines $T_5$. So $\thin\bbb_5\thin\ccc_2\thin\cdots$ $=$ $\thin\bbb_5\thin\ccc_2\thin\aaa\thick\cdots$, $\thin\bbb_5\thin\ccc_2\thin^{\bbb}\ccc^{\ddd}\thin\cdots$ or $\thin\bbb_5\thin\ccc_2\thin^{\ddd}\ccc^{\bbb}\thin\cdots$. If $\thin\bbb_5\thin\ccc_2\thin\cdots$ $=$ $\thin\bbb_5\thin\ccc_2\thin\aaa\thick\cdots$ or $\thin\bbb_5\thin\ccc_2\thin^{\bbb}\ccc^{\ddd}\thin\cdots$, we get $\ddd_2\ddd_3\cdots$ $=$ $\bbb\ddd_2\ddd_3\cdots$, contradicting the AVC. So we have $\thin\bbb_5\thin\ccc_2\thin\cdots$ $=$ $\thin\bbb_5\thin\ccc_2\thin^{\ddd}\ccc^{\bbb}\thin\cdots$ which determines $T_6$. Similarly, we can determine $T_7,T_8$. Then we get $\ddd_2\ddd_3\ddd_6\ddd_7\cdots$, contradicting the AVC. Therefore,  $\aaa^2\bbb\ccc^{k}$ $=$ $\thin^{\bbb}\ccc^{\ddd}\thin^{\bbb}\aaa^{\ddd}\thick^{\ddd}\aaa^{\bbb}\thin\bbb\thin\cdots$.		
			
			\begin{figure}[htp]
				\centering
				\begin{tikzpicture}[>=latex,scale=0.45] 
					
					\draw (-8,0)--(-8,-2)--(-6,-2)--(-6,-5)--(-2,-2)--(-4,-2)--(-4,0)
					(-2,-2)--(-3,-5)--(-1,-5)
					(8,0)--(8,-2)--(6,-2)--(6,-5)--(2,-2)--(4,-2)--(4,0)
					(2,-2)--(3,-5)--(1,-5)
					(-2,-2)--(2,-2);
					\draw[line width=1.5] (-4,-2)--(-6,-2)
					(0,-2)--(-1,-5)
					(4,-2)--(6,-2)
					(0,-2)--(1,-5)
					(0,0)--(0,-2);
					\node at (-6,0){\small $\ccc$};\node at (-2,0){\small $\aaa$};\node at (-7.6,-1.6){\small $\bbb$};\node at (-6,-1.6){\small $\aaa$};\node at (-4.4,-1.6){\small $\ddd$};\node at (-3.6,-1.6){\small $\bbb$};\node at (-2,-1.6){\small $\ccc$};\node at (-0.4,-1.6){\small $\ddd$};\node at (-5.6,-2.5){\small $\ddd$};\node at (-4,-2.5){\small $\aaa$};\node at (-3,-2.45){\small $\bbb$};\node at (-5.7,-4.4){\small $\ccc$};\node at (-1.8,-2.4){\small $\ccc$};\node at (-0.5,-2.4){\small $\ddd$};\node at (-2.5,-4.6){\small $\bbb$};\node at (-1.2,-4.6){\small $\aaa$};
					\node at (6,0){\small $\ccc$};\node at (2,0){\small $\aaa$};\node at (7.6,-1.6){\small $\bbb$};\node at (6,-1.6){\small $\aaa$};\node at (4.4,-1.6){\small $\ddd$};\node at (3.6,-1.6){\small $\bbb$};\node at (2,-1.6){\small $\ccc$};\node at (0.4,-1.6){\small $\ddd$};\node at (5.6,-2.5){\small $\ddd$};\node at (4,-2.4){\small $\aaa$};\node at (3.2,-2.45){\small $\bbb$};\node at (5.7,-4.4){\small $\ccc$};\node at (1.7,-2.4){\small $\ccc$};\node at (0.5,-2.4){\small $\ddd$};\node at (2.5,-4.6){\small $\bbb$};\node at (1.3,-4.6){\small $\aaa$};

					\fill (9,-1) circle (0.05); \fill (9.6,-1) circle (0.05);
					\fill (9.3,-1) circle (0.05);
					
					\node[draw,shape=circle, inner sep=0.5] at (-6,-0.8) {\small $1$};
					\node[draw,shape=circle, inner sep=0.5] at (-2,-0.8) {\small $2$};
					\node[draw,shape=circle, inner sep=0.5] at (2,-0.8) {\small $3$};
					\node[draw,shape=circle, inner sep=0.5] at (6,-0.8) {\small $4$};
					\node[draw,shape=circle, inner sep=0.5] at (-5,-3.5) {\small $5$};
					\node[draw,shape=circle, inner sep=0.5] at (-1.5,-3.5) {\small $6$};
					\node[draw,shape=circle, inner sep=0.5] at (1.5,-3.5) {\small $7$};
					\node[draw,shape=circle, inner sep=0.5] at (5,-3.5) {\small $8$};
					
				\end{tikzpicture}
				\caption{One possible AAD of  $\aaa^2\bbb\ccc^{k}=\protect\thin^{\bbb}\ccc^{\ddd}\protect\thin^{\bbb}\aaa^{\ddd}\protect\thick^{\ddd}\aaa^{\bbb}\protect\thin^{\ddd}\ccc^{\bbb}\protect\thin\cdots$.} \label{b=2d 1}
			\end{figure}
			
			The AAD $\aaa^2\bbb\ccc^{\frac{f-4}{6}}=\thin^{\bbb}\ccc_1^{\ddd}\thin^{\bbb}\aaa_2^{\ddd}\thick^{\ddd}\aaa_3^{\bbb}\thin\bbb\thin\cdots$ determines $T_1,T_2,T_3$ in Fig.\,\ref{b=2d 1''}. Then $R(\thin\aaa\thick\aaa\thin\bbb\thin\cdots)=\ccc^{k}$ and this $\ccc^{k}$ determines $k$ time zones ($2k$ or $\frac{f-4}{3}$ tiles). We have $\bbb_6\ccc_2\cdots=\thin\bbb_6\thin\ccc_2\thin^{\bbb}\aaa^{\ddd}\thick\cdots$, $\thin\bbb_6\thin\ccc_2\thin^{\bbb}\ccc^{\ddd}\thin\cdots$ or $\thin\bbb_6\thin\ccc_2\thin^{\ddd}\ccc^{\bbb}\thin\cdots$. If $\bbb_6\ccc_2\cdots=\thin\bbb_6\thin\ccc_2\thin^{\bbb}\aaa^{\ddd}\thick\cdots$ or $\thin\bbb_6\thin\ccc_2\thin^{\bbb}\ccc^{\ddd}\thin\cdots$, then we get $\bbb\ddd_2\ddd_3\cdots$, contradicting the AVC. Therefore,  $\bbb_6\ccc_2\cdots=\thin\bbb_6\thin\ccc_2\thin^{\ddd}\ccc^{\bbb}\thin\cdots$. This determines $T_7$. Then $\ddd_2\ddd_3\ddd_7\cdots=\aaa_8\ddd_2\ddd_3\ddd_7$ determines $T_8$. We have $\thin\bbb_6\thin\ccc_2\thin\ccc_7\thin\cdots=\aaa^2\bbb\ccc^{k}$ or $\bbb\ccc^{k+1}$. 			
			If $\thin\bbb_6\thin\ccc_2\thin\ccc_7\thin\cdots=\aaa^2\bbb\ccc^{k}=\thin^{\bbb}\aaa\thick\aaa\thin\bbb_6\thin\ccc_2\thin^{\ddd}\ccc_7^{\bbb}\thin\cdots$, this gives a vertex $\bbb^2\cdots$, contradicting the AVC. Therefore,  $\thin\bbb_6\thin\ccc_2\thin\ccc_7\thin\cdots=\bbb\ccc^{k+1}=\thin\bbb_6\thin\ccc_2\thin^{\ddd}\ccc_7^{\bbb}\thin^{\ddd}\ccc^{\bbb}\thin\cdots$ determines $T_9$. Then $R(\bbb_6\ccc_2\cdots)=\ccc^{k}$ and this $\ccc^{k}$ determines $k$ time zones ($2k$ or $\frac{f-4}{3}$ tiles). Similarly, we get $\bbb_8\ccc_3\cdots=\thin\bbb_8\thin^{\ddd}\ccc_3^{\bbb}\thin^{\ddd}\ccc_{10}^{\bbb}\thin\cdots=\bbb\ccc^{k+1}$ which determines $T_{10}$. Then $R(\bbb_8\ccc_3\ccc_{10}\cdots)=\ccc^{k-1}$ and this $\ccc^{k-1}$ determines $k-1$ time zones ($2k-2$ or $\frac{f-10}{3}$ tiles). Then $\bbb_3\ddd_{10}\cdots=\aaa_{12}\bbb_3\ddd_{10}$ determines $T_{12}$. So only two tiles are undetermined. By checking all possibilities, it turns out there are $3$ different ways to arrange these last two tiles, and Fig.\,\ref{b=2d 1''} shows one way with $\bbb_{11}\ccc^{\frac{f-4}{6}}\cdots=\bbb\ccc^{\frac{f+2}{6}}$. Then $\aaa_9\ddd\cdots=\aaa_9\bbb_{14}\ddd$ determines $T_{14}$; $\aaa_{14}\bbb_{9}\ccc^{\frac{f-4}{6}}\cdots=\aaa_{14}\aaa_{15}\bbb_{9}\ccc^{\frac{f-4}{6}},\aaa_{11}\ddd_{14}\ddd\cdots=\aaa_{11}\ddd_{14}\ddd_{15}\ddd,\bbb_4\ccc^{\frac{f-4}{6}}\cdots=\bbb_4\ccc^{\frac{f-4}{6}}\ccc_{15},\aaa_4\ddd_{5}\cdots=\aaa_4\bbb_{15}\ddd_{5}$ determine $T_{15}$. Centering $T_{14}, T_{15}$ in Fig.\,\ref{b=2d 4}, it becomes clear that they form a hexagon with $3$-fold symmetry, and the other two ways are obtained by rotating the $b$-edge $120^{\circ}$ and $240^{\circ}$ respectively. 
			
			\begin{figure}[htp]
				\centering
				\begin{tikzpicture}[>=latex,scale=0.6] 
					
					\draw (0,0)--(0,-10)
					(2,0)--(2,-10)
					(13,0)--(13,-10)
					(15,0)--(15,-10)
					(17,0)--(17,-10)
					(2,-6)--(5,-2)--(7,-4)--(9,-2)--(9,0)
					(2,-6)--(4,-6)--(5,-5)--(8,-8)--(6,-9)--(6,-10)
					(2,-6)--(4,-8)
					(7,-4)--(9,-4)--(11,-2)--(13,-2)--(13,0)
					(13,-2)--(11,-4)
					(7,-4)--(9,-6)
					(7,-4)--(8,-8)--(11,-7)--(13,-2)
					(4,-8)--(5.5,-7.6)--(8,-8);
					\draw[line width=1.5] (0,-8)--(2,-2)
					(5,0)--(5,-5)
					(4,-8)--(6,-9)
					(9,-2)--(11,-2)
					(9,-4)--(11,-4)
					(9,-6)--(11,-7)
					(13,-8)--(15,-2)
					(15,-8)--(17,-2)
					(11,-7)--(11,-10);

					\node at (1,0){\small $\ccc$};\node at (3.5,0){\small $\aaa$};
					\node at (7,0){\small $\aaa$};\node at (11,0){\small $\bbb$};\node at (14,0){\small $\ccc$};\node at (16,0){\small $\ccc$};
					
					\node at (0.3,-2){\small $\bbb$};\node at (1.7,-2){\small $\ddd$}; \node at (2.3,-2){\small $\bbb$};\node at (4.6,-2){\small $\ddd$};\node at (5.4,-2){\small $\ddd$};\node at (8.6,-2){\small $\bbb$};\node at (9.4,-1.7){\small $\aaa$};
					\node at (11,-1.7){\small $\ddd$};\node at (12.6,-1.7){\small $\ccc$};\node at (13.3,-2){\small $\bbb$};\node at (14.7,-2){\small $\ddd$};\node at (15.3,-2){\small $\bbb$};\node at (16.7,-2){\small $\ddd$};
					\node at (0.3,-6){\small $\aaa$};\node at (1.7,-4){\small $\aaa$};\node at (1.7,-6){\small $\bbb$};\node at (0.3,-8){\small $\ddd$};\node at (1,-10){\small $\ccc$};
					\node at (13.3,-6){\small $\aaa$};\node at (14.7,-4){\small $\aaa$};\node at (14.7,-8){\small $\bbb$};\node at (13.3,-8){\small $\ddd$};\node at (14,-10){\small $\ccc$};
					\node at (15.3,-6){\small $\aaa$};\node at (16.7,-4){\small $\aaa$};\node at (16.7,-8){\small $\bbb$};\node at (15.3,-8){\small $\ddd$};\node at (16,-10){\small $\ccc$};
					\node at (2.3,-5.2){\small $\ccc$};\node at (2.6,-5.7){\small $\ccc$};\node at (2.3,-6.8){\small $\ccc$};\node at (4.7,-3){\small $\ddd$};\node at (5.3,-3){\small $\aaa$};\node at (4,-5.6){\small $\bbb$};\node at (4.7,-4.8){\small $\aaa$};
					\node at (5.3,-4.8){\small $\ddd$};\node at (7,-3.6){\small $\ccc$};\node at (6.8,-4.3){\small $\bbb$};\node at (7.7,-3.7){\small $\ccc$};\node at (7.5,-5.1){\small $\ccc$};\node at (9,-2.3){\small $\ddd$};\node at (9,-3.7){\small $\bbb$};\node at (10.4,-2.3){\small $\aaa$};\node at (11,-2.4){\small $\bbb$};\node at (12.4,-2.3){\small $\ccc$};\node at (11,-3.7){\small $\ddd$};\node at (9.7,-3.7){\small $\aaa$};
					
					\node at (3.7,-8.2){\small $\ddd$};\node at (5.7,-9.3){\small $\aaa$};\node at (4,-10){\small $\bbb$};\node at (7.4,-7){\small $\ccc$};\node at (8.2,-7.6){\small $\bbb$};\node at (8.6,-6){\small $\ddd$};\node at (10.2,-6.95){\small $\aaa$};\node at (6.3,-9.3){\small $\bbb$};\node at (12.7,-3.5){\small $\ccc$};
					
					\node at (8,-8.4){\small $\ccc$};\node at (11.3,-7.3){\small $\ddd$};\node at (12.7,-8){\small $\bbb$};\node at (8.5,-10){\small $\aaa$};\node at (12,-10){\small $\aaa$};\node at (10.7,-7.5){\small $\ddd$};
					
					\node at (6,-8.6){\small $\ddd$};\node at (5.5,-8){\small $\bbb$};\node at (4.6,-8.05){\small $\aaa$};\node at (7.2,-8.15){\small $\ccc$};

					\fill (18,-5) circle (0.05); \fill (18.6,-5) circle (0.05);
					\fill (18.3,-5) circle (0.05);
					
					\fill (3,-6.65) circle (0.05); \fill (3.15,-6.45) circle (0.05); \fill (3.3,-6.25) circle (0.05);
					
					\fill (3+3.6,-6.65-1) circle (0.05); \fill (3.15+3.6,-6.45-1) circle (0.05); \fill (3.3+3.6,-6.25-1) circle (0.05);
					
					\fill (8.1,-4.75) circle (0.05); \fill (8.25,-4.5) circle (0.05); \fill (8.4,-4.25) circle (0.05);
					
					\fill (11.5,-3.8) circle (0.05); \fill (11.7,-3.9) circle (0.05); \fill (11.9,-4) circle (0.05);

					\node[draw,shape=circle, inner sep=0.5] at (1,-1) {\small $1$};
					\node[draw,shape=circle, inner sep=0.5] at (3.5,-1) {\small $2$};
					\node[draw,shape=circle, inner sep=0.5] at (7,-1) {\small $3$};
					\node[draw,shape=circle, inner sep=0.5] at (14,-1) {\small $4$};
					
					\node[draw,shape=circle, inner sep=0.5] at (14,-9) {\small $5$};
					\node[draw,shape=circle, inner sep=0.5] at (1,-9) {\small $6$};
					\node[draw,shape=circle, inner sep=0.5] at (3.8,-4.4) {\small $7$};
					\node[draw,shape=circle, inner sep=0.5] at (6.3,-5.1) {\small $8$};
					\node[draw,shape=circle, inner sep=0.5] at (4,-9) {\small $9$};
					\node[draw,shape=circle, inner sep=0.5] at (9,-3) {\footnotesize $10$};
					\node[draw,shape=circle, inner sep=0.5] at (9,-7) {\footnotesize $11$};
					\node[draw,shape=circle, inner sep=0.5] at (11,-1) {\footnotesize $12$};
					\node[draw,shape=circle, inner sep=0.5] at (11,-3) {\footnotesize $13$};
					\node[draw,shape=circle, inner sep=0.5] at (8.5,-9) {\footnotesize $14$};
					\node[draw,shape=circle, inner sep=0.5] at (12,-9) {\footnotesize $15$};

				\end{tikzpicture}
				\caption{ One special tiling for $\{(6k-2)\aaa\bbb\ddd,2\aaa\ddd^3,2\aaa^2\bbb\ccc^{k},4\bbb\ccc^{k+1}\}$.} \label{b=2d 1''}
			\end{figure}

			\begin{figure}[htp]
				\centering								
				\begin{tikzpicture}[>=latex,scale=0.28]

					\draw (0,3.2*10.2/4)--(-1*10.2/4,-0.8*10.2/4)--(-4*10.2/4,{(-sqrt(8)-0.8)*10.2/4})--(0,{(-sqrt(3)-0.8)*10.2/4})--(4*10.2/4,{(-sqrt(8)-0.8)*10.2/4})--(1*10.2/4,-0.8*10.2/4)--(0,3.2*10.2/4);
					
					\draw[line width=1.5] (0,3.15*10.2/4)--(0,{(-sqrt(3)-0.8)*10.2/4});
				
					\draw[dotted] (-1*10.2/4,-0.8*10.2/4)--(4*10.2/4,{(-sqrt(8)-0.8)*10.2/4})
					(1*10.2/4,-0.8*10.2/4)--(-4*10.2/4,{(-sqrt(8)-0.8)*10.2/4});
					
					\node at (0.2*10.2/4,-2.3*10.2/4){\small $\ddd$};
					\node at (-0.2*10.2/4,-2.3*10.2/4){\small $\ddd$};					
					
					\node at (-0.23*10.2/4,1.5*10.2/4){\small $\aaa$};
					\node at (0.23*10.2/4,1.5*10.2/4){\small $\aaa$};
					
					\node at (-0.8*10.2/4,-0.9*10.2/4){\small $\bbb$};					
					\node at (0.8*10.2/4,-0.9*10.2/4){\small $\bbb$};
					
					\node at (-3.1*10.2/4,-3.15*10.2/4){\small $\ccc$};
					\node at (3.1*10.2/4,-3.15*10.2/4){\small $\ccc$};
					
					\node at (-1.4*10.2/4,-0.6*10.2/4){\small $\aaa\ddd$};
					\node at (1.4*10.2/4,-0.6*10.2/4){\small $\aaa\ddd$};
					\node at (0*10.2/4,-2.9*10.2/4){\small $\aaa\ddd$};
					
					\node at (0,3.5*10.2/4){\small $\bbb\ccc^{k}$};
					\node at (-4.4*10.2/4,-3.5*10.2/4){\small $\bbb\ccc^{k}$};				
					\node at (4.5*10.2/4,-3.5*10.2/4){\small $\bbb\ccc^{k}$};

				\begin{scope}[xshift=22cm]
					\coordinate (O) at (0,0);
					\def\m{6}
					\pgfmathsetmacro\i{\m-1}
					\foreach \x in {0,...,\i}
					{
						\def\pointname{\x}
						\coordinate (\pointname) at ($(0,0) +(\x*360/\m:4cm)$)  ;
						
					}
					
					\draw (0)
					\foreach \x in {0,...,\i}
					{ -- (\x) } -- cycle;

					\draw (-2,3.45)--(-2,7)
					(-6,8)--(-6,4)--(-8,1)--(-6,-3)--(-2,-3.45)--(-4,0)
					(-6,4)--(-2,3.45);
					
					\draw[line width=1.5] (-2,3.45)--(2,-3.45)
					(-2,7)--(-6,8)
					(-4,0)--(-8,1);
					
					\draw[rotate =120] (0,0) [line width=1.5] (-2,7)--(-6,8)
					(-4,0)--(-8,1);
					
					\draw[rotate =120] (0,0) (-2,3.45)--(-2,7)
					(-6,8)--(-6,4)--(-8,1)--(-6,-3)--(-2,-3.45)--(-4,0)
					(-6,4)--(-2,3.45);

					\draw[rotate =240] (0,0) [line width=1.5] (-2,7)--(-6,8)
					(-4,0)--(-8,1);
					
					\draw[rotate =240] (0,0) (-2,3.45)--(-2,7)
					(-6,8)--(-6,4)--(-8,1)--(-6,-3)--(-2,-3.45)--(-4,0)
					(-6,4)--(-2,3.45);

					\node at (-2,2.5){\small $\aaa$}; \node at (-1.1,3){\small $\aaa$};
					\node at (1.8,2.7){\small $\bbb$}; \node at (-3.4,0){\small $\bbb$};
					\node at (3.4,0){\small $\ccc$}; \node at (-1.8,-2.8){\small $\ccc$}; 
					\node at (2,-2.5){\small $\ddd$}; \node at (1.1,-3){\small $\ddd$};
					\node at (-3,2.8){\small $\bbb$};\node at (-2.5,4){\small $\ccc$};
					\node at (-2.5,6.4){\small $\ddd$};\node at (-5.5,7.2){\small $\aaa$};
					\node at (-5.5,4.7){\small $\bbb$};\node at (-5.8,3.3){\small $\ccc$};
					\node at (-7.2,1.4){\small $\ddd$};\node at (-4.3,0.6){\small $\aaa$};
					\node at (-4.3,-0.6){\small $\ddd$};\node at (-7,0.2){\small $\aaa$};
					\node at (-5.6,-2.4){\small $\bbb$};\node at (-3.3,-2.4){\small $\ccc$};
					\node at (-0.7,3.9){\small $\ccc$}; \node at (0.5,5.9){\small $\bbb$};
					\node at (3.6,6){\small $\aaa$};\node at (1.8,4){\small $\ddd$};
					\node at (2.7,3.3){\small $\aaa$}; \node at (4.7,5.6){\small $\ddd$};
					\node at (4.1,1){\small $\bbb$};\node at (5.8,3.1){\small $\ccc$};
					\node at (6.5,2.3){\small $\bbb$};\node at (9,1){\small $\aaa$};
					\node at (6.9,-1){\small $\ddd$};\node at (4.8,0.1){\small $\ccc$};
					\node at (3.8,-1.3){\small $\ccc$}; \node at (2.5,-3.6){\small $\ddd$};
					\node at (5,-3.9){\small $\bbb$};\node at (3.3,-6.1){\small $\aaa$};
					\node at (1.5,-4){\small $\aaa$}; \node at (2.3,-6.8){\small $\ddd$};
					\node at (-1.1,-4.1){\small $\bbb$};\node at (-0.2,-6.6){\small $\ccc$};
					\node at (-2.3,-4.3){\small $\ccc$};\node at (-1.2,-7){\small $\bbb$};
					\node at (-3.6,-8.2){\small $\aaa$};\node at (-4.4,-5.5){\small $\ddd$};
					
					\fill (-1.7,6) circle (0.1); \fill (-1.3,5.9) circle (0.1);
					\fill (-0.9,5.8) circle (0.1);
					
					\fill (5.4,-2.4) circle (0.1); \fill (5.7,-2.1) circle (0.1);
					\fill (6,-1.8) circle (0.1);
					
					\fill (-5,-3.5) circle (0.1); \fill (-4.8,-3.9) circle (0.1); 
					\fill (-4.6,-4.3) circle (0.1); 
					
				\end{scope}	
				\end{tikzpicture}
				\caption{Three special tilings for $\{(6k-2)\aaa\bbb\ddd,2\aaa\ddd^3,2\aaa^2\bbb\ccc^{k},4\bbb\ccc^{k+1}\}$.} \label{b=2d 4}
			\end{figure}

			\subsubsection*{Subcase $\aaa\ccc^{2k+1}\ddd$ appears ($\aaa^2\bbb\ccc^k,\ccc^{\frac{f}{2}}$ never appear)} \label{b=2d-3}
			
			If $\aaa^2\bbb\ccc^k,\ccc^{\frac{f}{2}}$ never appear, Balance Lemma \ref{balance} implies \[\text{AVC}\sub\{\aaa\bbb\ddd,\aaa^2\ccc^{k}\ddd^2,\bbb\ccc^{k+1},\aaa\ccc^{2k+1}\ddd\}.\]				
			In Fig.\,\ref{b=2d 2}, we have the unique AAD $\aaa\ccc^{2k+1}\ddd=\thin^{\bbb}\ccc_1^{\ddd}\thin^{\bbb}\aaa_2^{\ddd}\thick^{\aaa}\ddd_3^{\ccc}\thin^{\bbb}\ccc_4^{\ddd}\thin\cdots$ which determines $T_2,T_3$. Then $R(\aaa_2\ddd_3\cdots)=\ccc^{2k+1}$ and this $\ccc^{2k+1}$ determines $2k+1$ time zones ($4k+2$ or $\frac{2f-2}{3}$ tiles).  Then $\thin\bbb_4\thin^{\ddd}\ccc_3^{\bbb}\thin\cdots=\thin\bbb_4\thin^{\ddd}\ccc_3^{\bbb}\thin^{\ddd}\ccc^{\bbb}\thin\cdots=\bbb\ccc^{k+1}$. Then $R(\bbb_4\ccc_3\cdots)=\ccc^{k}$ and this $\ccc^{k}$ determines $k$ time zones ($2k$ or $\frac{f-4}{3}$ tiles).  This tiling is exactly the second flip modification in Fig.\,\ref{flip1}.

				\begin{figure}[htp]
					\centering
					\begin{tikzpicture}[>=latex,scale=0.48] 
						
						\draw (0,0)--(0,-8)
						(6,0)--(6,-8)
						(0,-6)--(2,-2)--(6,-2)--(4,-6)--(0,-6);
						\draw[line width=1.5] (2,0)--(2,-2)
						(4,-6)--(4,-8);
						\node at (1,0){\small $\aaa$};
						\node at (4,0){\small $\ddd$};
						\node at (0.3,-2){\small $\bbb$};
						\node at (1.7,-2){\small $\ddd$};
						\node at (2.3,-1.7){\small $\aaa$};
						\node at (4,-1.6){\small $\bbb$};
						\node at (5.7,-1.7){\small $\ccc$};
						\node at (0.25,-5.1){\small $\ccc$};
						\node at (5.75,-3.1){\small $\ccc$};
						\node at (0.3,-6.3){\small $\ccc$};
						\node at (2,-6.4){\small $\bbb$};
						\node at (3.7,-6.3){\small $\aaa$};
						\node at (4.3,-6){\small $\ddd$};
						\node at (5.7,-6){\small $\bbb$};
						\node at (2,-8){\small $\ddd$};
						\node at (5,-8){\small $\aaa$};
						
						\fill (4,-2) circle (0.04);
						\fill (2,-6) circle (0.04);

						\fill (1,-5.8) circle (0.05);
						\fill (0.8,-5.6) circle (0.05);
						\fill (0.6,-5.4) circle (0.05);
						
						\fill (1+4.4,-5.8+3.1) circle (0.05);
						\fill (0.8+4.4,-5.6+3.1) circle (0.05);
						\fill (0.6+4.4,-5.4+3.1) circle (0.05);
						
						\fill (0,-2) circle (0.04);

						\foreach \a in {-4,0,1}
						{
							\begin{scope}[xshift=2*\a cm] 
								\draw (8,0)--(8,-8)
								(6,0)--(6,-8);
								\draw[line width=1.5] (6,-6)--(8,-2);
								\node at (7,0){\small $\ccc$};
								\node at (6.3,-2){\small $\bbb$};
								\node at (7.7,-2){\small $\ddd$};
								\node at (6.3,-4.8){\small $\aaa$};
								\node at (7.7,-3.2){\small $\aaa$};
								\node at (6.3,-6){\small $\ddd$};
								\node at (7.7,-6){\small $\bbb$};
								\node at (7,-8){\small $\ccc$};
								
							\end{scope}
						}

						\fill (11.5,-4) circle (0.05);
						\fill (11.8,-4) circle (0.05);
						\fill (12.1,-4) circle (0.05);
						
						\fill (10,-6) circle (0.04);

						\node[draw,shape=circle, inner sep=0.5] at (-1,-0.8) {\small $1$};
						\node[draw,shape=circle, inner sep=0.5] at (1,-0.8) {\small $2$};
						\node[draw,shape=circle, inner sep=0.5] at (4,-0.8) {\small $3$};
						\node[draw,shape=circle, inner sep=0.5] at (7,-0.8) {\small $4$};
						\node[draw,shape=circle, inner sep=0.5] at (7,-7.2) {\small $5$};
						\node[draw,shape=circle, inner sep=0.5] at (-1,-7.2) {\small $6$};
						\node[draw,shape=circle, inner sep=0.5] at (2,-7.2) {\small $8$};
						\node[draw,shape=circle, inner sep=0.5] at (5,-7.2) {\small $7$};
						
					\end{tikzpicture}
					\caption{ $T((6k+2)\aaa\bbb\ddd,2\bbb\ccc^{k+1},2\aaa\ccc^{2k+1}\ddd)$.} \label{b=2d 2}
				\end{figure}
				
				\subsubsection*{Subcase $\aaa^2\ccc^{k}\ddd^2$ appears ($\aaa^2\bbb\ccc^k,\aaa\ccc^{2k+1}\ddd,\ccc^{\frac{f}{2}}$ never appear)} \label{b=2d-4}

					By AVC, $\bbb^2\cdots$ is never a vertex. Then $\aaa^2\ccc^{k}\ddd^2$ has only two possible AAD.
					In Fig.\,\ref{b=2d 3'}, $\aaa^2\ccc^{k}\ddd^2=\thin^{\bbb}\aaa_1^{\ddd}\thick^{\ddd}\aaa_2^{\bbb}\thin^{\ddd}\ccc_3^{\bbb}\thin\cdots$ determines $T_1,T_2,T_3$. Then $\bbb_2\ddd_3\cdots=\aaa_4\bbb_2\ddd_3$ determines $T_4$; $\bbb_4\ccc_2\cdots=\bbb_4\ccc_2\ccc_5\cdots=\bbb\ccc^{k+1}$. By $\ccc_5$, we get $\bbb\ddd_1\ddd_2\cdots$ or $\ddd_1\ddd_2\ddd\cdots$, contradicting the AVC. Therefore, $\aaa^2\ccc^{k}\ddd^2=\thin\aaa\thick\ddd\thin\cdots\thin\aaa\thick\ddd\thin\cdots$.

					\begin{figure}[htp]
						\centering	    	
						\begin{tikzpicture}[>=latex,scale=0.5] 
							
							\draw (0,0)--(0,-2)--(7,-2)--(7,0)
							(4,-4)--(5,-2)--(7,-4)--(8.5,-2)--(10,-2)--(10,0);
							\draw[line width=1.5] (3,0)--(3,-2)
							(7,-2)--(8.5,-2);
							\draw[dotted] (3,-2)--(2,-4)--(4,-4);		    		
							
							\node at (1.5,0){\small $\aaa$};\node at (5,0){\small $\aaa$};\node at (8.5,0){\small $\ccc$};\node at (0.3,-1.6){\small $\bbb$};\node at (1.5,-1.6){\small $\ccc$};\node at (2.7,-1.6){\small $\ddd$};\node at (3.3,-1.6){\small $\ddd$};\node at (5,-1.6){\small $\ccc$};\node at (6.7,-1.6){\small $\bbb$};\node at (7.3,-1.6){\small $\ddd$};\node at (8.5,-1.6){\small $\aaa$};\node at (9.7,-1.6){\small $\bbb$};\node at (4.5,-2.4){\small $\ccc$};\node at (6,-2.4){\small $\bbb$};\node at (7,-2.3){\small $\aaa$};\node at (7.9,-2.4){\small $\ddd$};\node at (7,-3.6){\small $\ccc$};

							\fill (1.5,-2) circle (0.04);
							
							\fill (11,-1) circle (0.05);
							\fill (11.3,-1) circle (0.05);
							\fill (11.6,-1) circle (0.05); 
							
							\fill (4.9,-2.7) circle (0.05);
							\fill (5.1,-2.7) circle (0.05);
							\fill (5.3,-2.7) circle (0.05);

							\node[draw,shape=circle, inner sep=0.5] at (1.5,-0.8) {\small $1$};
							\node[draw,shape=circle, inner sep=0.5] at (5,-0.8) {\small $2$};
							\node[draw,shape=circle, inner sep=0.5] at (8.5,-0.8) {\small $3$};
							\node[draw,shape=circle, inner sep=0.5] at (7,-2.9) {\small $4$};
							\node[draw,shape=circle, inner sep=0.5] at (3.5,-3) {\small $5$};

						\end{tikzpicture}
						\caption{One possible AAD of $\aaa^2\ccc^{k}\ddd^2=\protect\thin\aaa\protect\thick\aaa\protect\thin\cdots$.} \label{b=2d 3'}
					\end{figure}
					
					In Fig.\,\ref{b=2d 3}, $\aaa^2\ccc^{k}\ddd^2=\thin^{\bbb}\aaa_1^{\ddd}\thick^{\aaa}\ddd_2^{\ccc}\thin\cdots\thin^{\bbb}\ccc^{\ddd}\thin\cdots\thin^{\bbb}\aaa_3^{\ddd}\thick^{\aaa}\ddd_4^{\ccc}\thin\cdots\thin^{\bbb}\ccc^{\ddd}\thin\cdots$ determines $T_1,T_2,T_3,T_4$. Then $R(\aaa^2\ddd^2\cdots)=\ccc^{k}$ determines $k$ time zones ($2k$ or $\frac{f-4}{3}$ tiles). Then $R(\bbb\ccc_2\cdots)=R(\bbb\ccc_4\cdots)=\ccc^{k}$ and each of these two $\ccc^{k}$ determines $k$ time zones ($2k$ or $\frac{f-4}{3}$ tiles). This tiling can also be obtained by applying the second flip modification in Fig.\,\ref{flip1} two times.

					\begin{figure}[htp]
						\centering	    	
						\begin{tikzpicture}[>=latex,scale=0.5] 
							\foreach \a in {0,1}
							{
								\begin{scope}[xshift=12*\a cm] 
									\draw (0,0)--(0,-8)
									(6,0)--(6,-8)
									(0,-6)--(2,-2)--(6,-2)--(4,-6)--(0,-6);
									\draw[line width=1.5] (2,0)--(2,-2)
									(4,-6)--(4,-8);
									\node at (1,0){\small $\aaa$};
									\node at (4,0){\small $\ddd$};
									\node at (0.3,-2){\small $\bbb$};
									\node at (1.7,-2){\small $\ddd$};
									\node at (2.3,-1.7){\small $\aaa$};
									\node at (4,-1.6){\small $\bbb$};
									\node at (5.7,-1.7){\small $\ccc$};
									\node at (0.25,-5.1){\small $\ccc$};
									\node at (5.75,-3.1){\small $\ccc$};
									\node at (0.3,-6.3){\small $\ccc$};
									\node at (2,-6.4){\small $\bbb$};
									\node at (3.7,-6.3){\small $\aaa$};
									\node at (4.3,-6){\small $\ddd$};
									\node at (5.7,-6){\small $\bbb$};
									\node at (2,-8){\small $\ddd$};
									\node at (5,-8){\small $\aaa$};
									
									\fill (4,-2) circle (0.04);
									\fill (2,-6) circle (0.04);

									\fill (1,-5.8) circle (0.05);
									\fill (0.8,-5.6) circle (0.05);
									\fill (0.6,-5.4) circle (0.05);
									
									\fill (1+4.4,-5.8+3.1) circle (0.05);
									\fill (0.8+4.4,-5.6+3.1) circle (0.05);
									\fill (0.6+4.4,-5.4+3.1) circle (0.05);
									
									\fill (6,-6) circle (0.04);
									\fill (0,-2) circle (0.04);
									\fill (8,-2) circle (0.04);
									\fill (10,-6) circle (0.04);
									
								\end{scope}
							}
							
							\foreach \b in {1,7}
							{
								\begin{scope}[xshift=2*\b cm] 
									\draw (6,0)--(6,-8)
									(8,0)--(8,-8);
									\draw[line width=1.5] (6,-6)--(8,-2);
									\node at (7,0){\small $\ccc$};
									\node at (6.3,-2){\small $\bbb$};
									\node at (7.7,-2){\small $\ddd$};
									\node at (6.3,-4.8){\small $\aaa$};
									\node at (7.7,-3.2){\small $\aaa$};
									\node at (6.3,-6){\small $\ddd$};
									\node at (7.7,-6){\small $\bbb$};
									\node at (7,-8){\small $\ccc$};
									
									\fill (8.6,-4) circle (0.05);
									\fill (8.9,-4) circle (0.05);
									\fill (9.2,-4) circle (0.05);
									
									\fill (4.6,-4) circle (0.05);
									\fill (4.9,-4) circle (0.05);
									\fill (5.2,-4) circle (0.05);
								\end{scope}
							}
							
							\fill (0,-2) circle (0.04);
							\fill (12,-6) circle (0.04);

							\node[draw,shape=circle, inner sep=0.5] at (1,-0.8) {\small $1$};
							\node[draw,shape=circle, inner sep=0.5] at (4,-0.8) {\small $2$};
							\node[draw,shape=circle, inner sep=0.5] at (13,-0.8) {\small $3$};
							\node[draw,shape=circle, inner sep=0.5] at (16,-0.8) {\small $4$};
							
						\end{tikzpicture}
						\caption{Many different tilings for $\{6k\aaa\bbb\ddd,2\aaa^2\ccc^{k}\ddd^2,4\bbb\ccc^{k+1}\}$.} \label{b=2d 3}
					\end{figure}							
					
					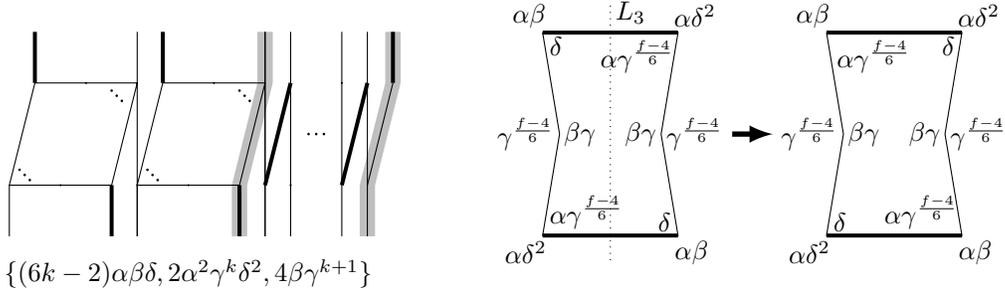
\begin{figure}[htp]
						\centering	    	
						\begin{tikzpicture}[>=latex,scale=0.34] 
							\begin{scope}[xshift=-12cm]						
							\fill[gray!50]  (1.7+15,0+4)--(1.7+15,-2+4)--(-0.3+16,-6+4)--(-0.3+16,-8+4)--(0.3+16,-8+4)--(0.3+16,-6+4)--(2.3+15,-2+4)--(2.3+15,0+4)--(1.7+15,0+4);
							\fill[gray!50]  (1.7+10,0+4)--(1.7+10,-2+4)--(-0.3+11,-6+4)--(-0.3+11,-8+4)--(0.3+11,-8+4)--(0.3+11,-6+4)--(2.3+10,-2+4)--(2.3+10,0+4)--(1.7+10,0+4);
							\foreach \a in {1,6}
							{
								\begin{scope}[xshift=\a cm] 
									\draw (1,-6+4)--(1,-8+4)
									(6,0+4)--(6,-8+4)
									(1,-6+4)--(2,-2+4)--(6,-2+4)--(5,-6+4)--(1,-6+4);
									\draw[line width=1.5] (2,0+4)--(2,-2+4)
									(5,-6+4)--(5,-8+4);

									\fill (4,-2+4) circle (0.04);
									\fill (3,-6+4) circle (0.04);
			
									\fill (1+0.8,-5.8+4) circle (0.05);
									\fill (0.8+0.8,-5.6+4) circle (0.05);
									\fill (0.6+0.8,-5.4+4) circle (0.05);
									
									\fill (1+4.4,-5.8+3.1+4) circle (0.05);
									\fill (0.8+4.4,-5.6+3.1+4) circle (0.05);
									\fill (0.6+4.4,-5.4+3.1+4) circle (0.05);
									
									\fill (6,-6+4) circle (0.04);

								\end{scope}
							}
							\fill (15,-2+4) circle (0.04);
							\fill (13,-6+4) circle (0.04);
							\foreach \b in {6,9}
							{
								\begin{scope}[xshift=\b cm] 
									\draw (6,0+4)--(6,-8+4)
									(7,0+4)--(7,-8+4);
									\draw[line width=1.5] (6,-6+4)--(7,-2+4);

								\end{scope}
							}
							\fill (4.7+9,-4+4) circle (0.05);
							\fill (5+9,-4+4) circle (0.05);
							\fill (5.3+9,-4+4) circle (0.05);
							\draw [line width=1.5] (17,0+4)--(17,-2+4);
							\draw (17,-2+4)--(16,-6+4);
							\node at (9,-9.45+4){\small $\{(6k-2)\aaa\bbb\ddd,2\aaa^2\ccc^{k}\ddd^2,4\bbb\ccc^{k+1}\}$};
						\end{scope}				    
			
						\begin{scope}[xshift=13.5 cm]
							\draw(-2*45/34,3*45/34)--(2*45/34,3*45/34)--(1.5*45/34,0)--(2*45/34,-3*45/34)--(-2*45/34,-3*45/34)--(-1.5*45/34,0)--(-2*45/34,3*45/34);
							\draw[line width=1.5] (-2*45/34,3*45/34)--(2*45/34,3*45/34)
								(2*45/34,-3*45/34)--(-2*45/34,-3*45/34);
								
							\draw[line width=2pt, ->](3.6*45/34,0)--(4.8*45/34,0);
				
							\draw[dotted] (0,3.8*45/34)--(0,-3.8*45/34);						  																													
							\node at (-1.6*45/34,2.6*45/34){\small $\ddd$};						
							\node at (1.6*45/34,-2.6*45/34){\small $\ddd$};
							
							\node at (0.8*45/34,2.3*45/34){\small $\aaa\ccc^{\frac{f-4}{6}}$};
							\node at (-0.7*45/34,-2.3*45/34){\small  $\aaa\ccc^{\frac{f-4}{6}}$};
														
							\node at (2.5*45/34,3.5*45/34){\small $\aaa\ddd^2$};												
							\node at (-2.5*45/34,-3.5*45/34){\small $\aaa\ddd^2$};
							
							\node at (0.9*45/34,0){\small $\bbb\ccc$};
							\node at (2.5*45/34,0){\small $\ccc^{\frac{f-4}{6}}$};
							\node at (-0.9*45/34,0){\small $\bbb\ccc$};
							\node at (-2.5*45/34,0){\small $\ccc^{\frac{f-4}{6}}$};

							\node at (-2.4*45/34,3.5*45/34){\small $\aaa\bbb$};						
							\node at (2.4*45/34,-3.5*45/34){\small $\aaa\bbb$};	
							
							\node at (0.6*45/34,3.6*45/34){\small $L_3$};
						\end{scope}
						\begin{scope}[xshift=24.6 cm]
							\draw(-2*45/34,3*45/34)--(2*45/34,3*45/34)--(1.5*45/34,0)--(2*45/34,-3*45/34)--(-2*45/34,-3*45/34)--(-1.5*45/34,0)--(-2*45/34,3*45/34);
							\draw[line width=1.5] (-2*45/34,3*45/34)--(2*45/34,3*45/34)
							(2*45/34,-3*45/34)--(-2*45/34,-3*45/34);				
							
							\node at (1.6*45/34,2.6*45/34){\small $\ddd$};						
							\node at (-1.6*45/34,-2.6*45/34){\small $\ddd$};
							
							\node at (-0.6*45/34,2.3*45/34){\small$\aaa\ccc^{\frac{f-4}{6}}$};
							\node at (0.8*45/34,-2.3*45/34){\small $\aaa\ccc^{\frac{f-4}{6}}$};
							
							\node at (2.5*45/34,3.5*45/34){\small $\aaa\ddd^2$};												
							\node at (-2.5*45/34,-3.5*45/34){\small $\aaa\ddd^2$};
							
							\node at (0.9*45/34,0){\small $\bbb\ccc$};
							\node at (2.5*45/34,0){\small $\ccc^{\frac{f-4}{6}}$};
							\node at (-0.9*45/34,0){\small $\bbb\ccc$};
							\node at (-2.5*45/34,0){\small $\ccc^{\frac{f-4}{6}}$};

							\node at (-2.4*45/34,3.5*45/34){\small $\aaa\bbb$};						
							\node at (2.4*45/34,-3.5*45/34){\small $\aaa\bbb$};	
							
						\end{scope}	
						\end{tikzpicture}
						\caption{A special flip modification ($\frac{f+2}{3}$ tiles flipped).} \label{flip2}
					\end{figure}
		If $\text{AVC}\sub\{\aaa\bbb\ddd,\bbb\ccc^{k+1}\}$, there is no solution satisfying Balance Lemma \ref{balance}. 
				
		In fact one special tiling in Fig.\,\ref{b=2d 3}, as shown in the first picture of Fig.\,\ref{flip2}, is related to Fig.\,\ref{b=2d 1''} by a special flip modification along $L_3$ in Fig.\,\ref{flip2}. 		
		
		All of the above tilings belong to the second infinite sequence in Table \ref{Tab-1.2}.
		
		\subsection*{Case $2$: $\{x_1, x_2\}=\{ x_3, x_4\}$}	
    We can fix $x_2=1- \frac\bbb2$ and $x_3=\frac\ccc2$ in \eqref{3-12}. Then either $x_1=x_3$, $x_4=x_2$ as listed in the left of Table \ref{Tab-6}, or $x_1=x_4, x_2=x_3$ as listed in the right. All solutions are ruled out by the fact listed in the other column of Table \ref{Tab-6} except one solution $-\aaa+\frac{\ccc}{2}=-\ddd+\frac\bbb2$, $1-\frac\bbb2=\frac\ccc2$. By $\bbb+\ccc=2$, we have $\bbb\cdots=\aaa^x\bbb\ddd^y$. By $\#\aaa+\#\ddd\ge2\#\bbb=2f$ and $\aaa\neq\ddd$, there is only one solution satisfying Balance Lemma \ref{balance}: $\{f\aaa\bbb\ddd,2\ccc^{\frac f2}\}$. Then $\aaa=-\frac12+\frac4f,\bbb=2-\frac4f,\ccc=\frac4f,\ddd=\frac12$. By $\aaa>0$, we get $f<8$ which forces $f=6$. This solution admits only a $2$-layer earth map tiling $T(6\aaa\bbb\ddd,2\ccc^3)$, and is Case $(1,8,4,3)/6$ in Table \ref{Tab-1.1}. 	\label{discrete-2'}
    
   \begin{table}[htp]  
	\centering        
	\begin{minipage}{0.55\textwidth}
		\scriptsize{\begin{tabular}{c|c|c}
				$x_1(=x_3=\frac\ccc2)$&$x_4(=x_2=1-\frac\bbb2)$& \\
				\hline 
				$\aaa-\frac\ccc2$&$\ddd-\frac\bbb2$&	$\ddd=1$	\\
				$-\aaa+\frac\ccc2$ & $1+\ddd-\frac\bbb2$&$\ddd=0$	\\
				$-\aaa+\frac\ccc2$ & $-\ddd+\frac\bbb2$&$\aaa=0$\\
				$\aaa+\frac\ccc2$ & $-1+\ddd+\frac\bbb2$&$\aaa=0$\\	
				$1-\aaa-\frac\ccc2$ & $1-\frac\bbb2$&$f=4$\\
		\end{tabular}}
	\end{minipage}
	\raisebox{0em}{\begin{minipage}{0.4\textwidth}
			\scriptsize{\begin{tabular}{c|c}
					$x_1=x_4$ $(\&\,1-\frac\bbb2=\frac\ccc2)$& \\
					\hline					
					$\aaa-\frac\ccc2=\ddd-\frac\bbb2$&$f<4$\\
					$-\aaa+\frac\ccc2=1+\ddd-\frac\bbb2$ &$\aaa+\ddd=0$\\
					$-\aaa+\frac\ccc2=-\ddd+\frac\bbb2$ &\\
					$\aaa+\frac\ccc2=-1+\ddd+\frac\bbb2$ &$\bbb+\ddd>2$\\
					$1-\aaa-\frac\ccc2=-1+\ddd+\frac\bbb2$ &$f=4$\\
			\end{tabular}}
	\end{minipage}}
    \caption{Ten subcases of $\{x_1, x_2\}=\{ x_3, x_4\}.$}\label{Tab-6}      
    \end{table}

		\subsection*{Case $3$: $\{x_1, x_2\}=\{ \frac{1}{6},\theta\}$ and $\{x_3,x_4\}=\{\frac{\theta}{2},\frac{1}{2}-\frac{\theta}{2}\}$, or $\{x_3,x_4\}$ $=$ $\{\frac{1}{6},\theta\}$ and $\{x_1,x_2\}=\{\frac{\theta}{2},\frac{1}{2}-\frac{\theta}{2}\}$, for some $0<\theta\le\frac{1}{2}$}

		\begin{table*}[htp]                        
			\centering             
			\resizebox{\textwidth}{90mm}{\begin{tabular}{cccc|cccc|c}	 
					
					$\theta$&$\frac16$&$\frac{\theta}{2}$&$\frac{1}{2}-\frac{\theta}{2}$&$\aaa$&$\bbb$&$\ccc$&$\ddd$&$\ccc+2\ddd>1,\bbb+\ddd<2$ (Lemma $\ref{geometry4}'$\&\ref{proposition-3}) \\
					\hline 
					$\aaa-\frac{\ccc}{2}$&$1-\frac{\bbb}{2}$&$\frac{\ccc}{2}$&$\ddd-\frac{\bbb}{2}$&$\frac{3\theta}{2}$&$\frac{5}{3}$&$\theta$&$\frac{4}{3}-\frac{\theta}{2}$&$\aaa+\bbb+\ccc+\ddd>\frac83$ \\
					$-\aaa+\frac{\ccc}{2}$&$1-\frac{\bbb}{2}$&$\frac{\ccc}{2}$&$1+\ddd-\frac{\bbb}{2}$&$-\frac{\theta}{2}$&$\frac{5}{3}$&$\theta$&$\frac{1}{3}-\frac{\theta}{2}$&$\aaa<0$  \\
					$-\aaa+\frac{\ccc}{2}$&$1-\frac{\bbb}{2}$&$\frac{\ccc}{2}$&$-\ddd+\frac{\bbb}{2}$&$-\frac{\theta}{2}$&$\frac{5}{3}$&$\theta$&$\frac{1}{3}+\frac{\theta}{2}$&$\aaa<0$  \\
					$\aaa+\frac{\ccc}{2}$&$1-\frac{\bbb}{2}$&$\frac{\ccc}{2}$&$-1+\ddd+\frac{\bbb}{2}$&$\frac{\theta}{2}$&$\frac{5}{3}$&$\theta$&$\frac{2}{3}-\frac{\theta}{2}$& $\bbb+\ddd>2$  \\
					$1-\aaa-\frac{\ccc}{2}$&$1-\frac{\bbb}{2}$&$\frac{\ccc}{2}$&$-1+\ddd+\frac{\bbb}{2}$&$1-\frac{3\theta}{2}$&$\frac{5}{3}$&$\theta$&$\frac{2}{3}-\frac{\theta}{2}$& $\aaa+\bbb+\ccc+\ddd>\frac83$   \\
					\multicolumn{9}{c}{ }\\
					$\theta$&$\frac16$&$\frac{1}{2}-\frac{\theta}{2}$&$\frac{\theta}{2}$&$\aaa$&$\bbb$&$\ccc$&$\ddd$&$\ccc+2\ddd>1,\bbb+\ddd<2$ (Lemma $\ref{geometry4}'$\&\ref{proposition-3})\\
					\hline 
					$\aaa-\frac{\ccc}{2}$&$1-\frac{\bbb}{2}$&$\frac{\ccc}{2}$&$\ddd-\frac{\bbb}{2}$&$\frac12+\frac{\theta}{2}$&$\frac{5}{3}$&$1-\theta$&$\frac{5}{6}+\frac{\theta}{2}$&$\aaa+\bbb+\ccc+\ddd=4>\frac83$ \\
					$-\aaa+\frac{\ccc}{2}$&$1-\frac{\bbb}{2}$&$\frac{\ccc}{2}$&$1+\ddd-\frac{\bbb}{2}$&$\frac12-\frac{3\theta}{2}$&$\frac{5}{3}$&$1-\theta$&$-\frac{1}{6}+\frac{\theta}{2}$&$\ccc+2\ddd=\frac23$  \\
					$-\aaa+\frac{\ccc}{2}$&$1-\frac{\bbb}{2}$&$\frac{\ccc}{2}$&$-\ddd+\frac{\bbb}{2}$&$\frac12-\frac{3\theta}{2}$&$\frac{5}{3}$&$1-\theta$&$\frac{5}{6}-\frac{\theta}{2}$&$\bbb+\ddd>2$  \\
					$\aaa+\frac{\ccc}{2}$&$1-\frac{\bbb}{2}$&$\frac{\ccc}{2}$&$-1+\ddd+\frac{\bbb}{2}$&$-\frac12+\frac{3\theta}{2}$&$\frac{5}{3}$&$1-\theta$&$\frac{1}{6}+\frac{\theta}{2}$&$\aaa>0\Rightarrow \theta>\frac13$ but \\
					&&&&&&&& $\bbb+\ddd<2 \Rightarrow \theta<\frac13$\\
					$1-\aaa-\frac{\ccc}{2}$&$1-\frac{\bbb}{2}$&$\frac{\ccc}{2}$&$-1+\ddd+\frac{\bbb}{2}$&$\frac12-\frac{\theta}{2}$&$\frac{5}{3}$&$1-\theta$&$\frac{1}{6}+\frac{\theta}{2}$&$\aaa<\ddd \Rightarrow \theta>\frac13$ but \\
					&&&&&&&& $\aaa+\bbb+\ccc+\ddd>2 \Rightarrow \theta<\frac13$\\
					\multicolumn{9}{c}{ }\\
					$\frac16$&$\theta$&$\frac{\theta}{2}$&$\frac{1}{2}-\frac{\theta}{2}$&$\aaa$&$\bbb$&$\ccc$&$\ddd$&$\ccc+2\ddd>1,\bbb+\ddd<2$ (Lemma $\ref{geometry4}'$\&\ref{proposition-3})\\
					\hline 
					$\aaa-\frac{\ccc}{2}$&$1-\frac{\bbb}{2}$&$\frac{\ccc}{2}$&$\ddd-\frac{\bbb}{2}$&$\frac16+\frac{\theta}{2}$&$2-2\theta$&$\theta$&$\frac{3}{2}-\frac{3\theta}{2}$&$\bbb>1 \Rightarrow \theta<\frac12$ but \\
					&&&&&&&& $\aaa+\bbb+\ccc+\ddd\le\frac83 \Rightarrow \theta\ge\frac12$\\
					$-\aaa+\frac{\ccc}{2}$&$1-\frac{\bbb}{2}$&$\frac{\ccc}{2}$&$1+\ddd-\frac{\bbb}{2}$&$-\frac16+\frac{\theta}{2}$&$2-2\theta$&$\theta$&$\frac{1}{2}-\frac{3\theta}{2}$&$\ccc+2\ddd<1$  \\
					$-\aaa+\frac{\ccc}{2}$&$1-\frac{\bbb}{2}$&$\frac{\ccc}{2}$&$-\ddd+\frac{\bbb}{2}$&$-\frac16+\frac{\theta}{2}$&$2-2\theta$&$\theta$&$\frac{1}{2}-\frac{\theta}{2}$&$\ccc+2\ddd=1$ \\
					
					$\aaa+\frac{\ccc}{2}$&$1-\frac{\bbb}{2}$&$\frac{\ccc}{2}$&$-1+\ddd+\frac{\bbb}{2}$&$\frac16-\frac{\theta}{2}$&$2-2\theta$&$\theta$&$\frac{1}{2}+\frac{\theta}{2}$&$\aaa>0 \Rightarrow \theta<\frac13$ but \\
					&&&&&&&& $\bbb+\ddd<2 \Rightarrow \theta>\frac13$\\
					$1-\aaa-\frac{\ccc}{2}$&$1-\frac{\bbb}{2}$&$\frac{\ccc}{2}$&$-1+\ddd+\frac{\bbb}{2}$&$\frac56-\frac{\theta}{2}$&$2-2\theta$&$\theta$&$\frac{1}{2}+\frac{\theta}{2}$&$\aaa+\bbb+\ccc+\ddd>\frac83$  \\
					\multicolumn{9}{c}{ }\\
					$\frac16$&$\theta$&$\frac{1}{2}-\frac{\theta}{2}$&$\frac{\theta}{2}$&$\aaa$&$\bbb$&$\ccc$&$\ddd$&$\ccc+2\ddd>1,\bbb+\ddd<2$ (Lemma $\ref{geometry4}'$\&\ref{proposition-3}) \\
					\hline 
					$\aaa-\frac{\ccc}{2}$&$1-\frac{\bbb}{2}$&$\frac{\ccc}{2}$&$\ddd-\frac{\bbb}{2}$&$\frac23-\frac{\theta}{2}$&$2-2\theta$&$1-\theta$&$1-\frac{\theta}{2}$&$\bbb>1 \Rightarrow \theta<\frac12$ but \\
					&&&&&&&& $\aaa+\bbb+\ccc+\ddd\le\frac83 \Rightarrow \theta\ge\frac12$\\
					$-\aaa+\frac{\ccc}{2}$&$1-\frac{\bbb}{2}$&$\frac{\ccc}{2}$&$1+\ddd-\frac{\bbb}{2}$&$\frac13-\frac{\theta}{2}$&$2-2\theta$&$1-\theta$&$-\frac{\theta}{2}$&$\ddd<0$ \\
					$-\aaa+\frac{\ccc}{2}$&$1-\frac{\bbb}{2}$&$\frac{\ccc}{2}$&$-\ddd+\frac{\bbb}{2}$&$\frac13-\frac{\theta}{2}$&$2-2\theta$&$1-\theta$&$1-\frac{3\theta}{2}$&\textbf{$\surd$ \,\,\, Subcase} $1$ \\
					$\aaa+\frac{\ccc}{2}$&$1-\frac{\bbb}{2}$&$\frac{\ccc}{2}$&$-1+\ddd+\frac{\bbb}{2}$&$-\frac13+\frac{\theta}{2}$&$2-2\theta$&$1-\theta$&$\frac{3\theta}{2}$&$\aaa<0$  \\
					$1-\aaa-\frac{\ccc}{2}$&$1-\frac{\bbb}{2}$&$\frac{\ccc}{2}$&$-1+\ddd+\frac{\bbb}{2}$&$\frac13+\frac{\theta}{2}$&$2-2\theta$&$1-\theta$&$\frac{3\theta}{2}$&$\aaa+\bbb+\ccc+\ddd>\frac83$  \\
					\multicolumn{9}{c}{ }\\
					$\frac{\theta}{2}$&$\frac{1}{2}-\frac{\theta}{2}$&$\frac{1}{6}$&$\theta$&$\aaa$&$\bbb$&$\ccc$&$\ddd$&$\ccc+2\ddd>1,\bbb+\ddd<2$ (Lemma $\ref{geometry4}'$\&\ref{proposition-3}) \\
					\hline 
					$\aaa-\frac{\ccc}{2}$&$1-\frac{\bbb}{2}$&$\frac{\ccc}{2}$&$\ddd-\frac{\bbb}{2}$&$\frac16+\frac{\theta}{2}$&$1+\theta$&$\frac13$&$\frac12+\frac{3\theta}{2}$&\textbf{$\surd$ \,\,\, Subcase} $2$ \\
					$-\aaa+\frac{\ccc}{2}$&$1-\frac{\bbb}{2}$&$\frac{\ccc}{2}$&$1+\ddd-\frac{\bbb}{2}$&$\frac16-\frac{\theta}{2}$&$1+\theta$&$\frac13$&$-\frac12+\frac{3\theta}{2}$&$\aaa>0\Rightarrow \theta<\frac13$ but \\
					&&&&&&&& $\ddd>0 \Rightarrow \theta>\frac13$\\
					$-\aaa+\frac{\ccc}{2}$&$1-\frac{\bbb}{2}$&$\frac{\ccc}{2}$&$-\ddd+\frac{\bbb}{2}$&$\frac16-\frac{\theta}{2}$&$1+\theta$&$\frac13$&$\frac12-\frac{\theta}{2}$& $\aaa+\bbb+\ccc+\ddd=2$ \\
					$\aaa+\frac{\ccc}{2}$&$1-\frac{\bbb}{2}$&$\frac{\ccc}{2}$&$-1+\ddd+\frac{\bbb}{2}$&$-\frac16+\frac{\theta}{2}$&$1+\theta$&$\frac13$&$\frac12+\frac{\theta}{2}$&$\aaa>0 \Rightarrow \theta>\frac13$ but \\
					&&&&&&&& $\bbb+\ddd<2 \Rightarrow \theta<\frac13$\\
					$1-\aaa-\frac{\ccc}{2}$&$1-\frac{\bbb}{2}$&$\frac{\ccc}{2}$&$-1+\ddd+\frac{\bbb}{2}$&$\frac56-\frac{\theta}{2}$&$1+\theta$&$\frac13$&$\frac12+\frac{\theta}{2}$&$\aaa+\bbb+\ccc+\ddd>\frac38$  \\
					\multicolumn{9}{c}{ }\\
					$\frac{\theta}{2}$&$\frac{1}{2}-\frac{\theta}{2}$&$\theta$&$\frac{1}{6}$&$\aaa$&$\bbb$&$\ccc$&$\ddd$&$\ccc+2\ddd>1,\bbb+\ddd<2$ (Lemma $\ref{geometry4}'$\&\ref{proposition-3})\\
					\hline 
					$\aaa-\frac{\ccc}{2}$&$1-\frac{\bbb}{2}$&$\frac{\ccc}{2}$&$\ddd-\frac{\bbb}{2}$&$\frac{3\theta}{2}$&$1+\theta$&$2\theta$&$\frac23+\frac{\theta}{2}$&\textbf{$\surd$ \,\,\, Subcase} $3$  \\
					$-\aaa+\frac{\ccc}{2}$&$1-\frac{\bbb}{2}$&$\frac{\ccc}{2}$&$1+\ddd-\frac{\bbb}{2}$&$\frac{\theta}{2}$&$1+\theta$&$2\theta$&$-\frac13+\frac{\theta}{2}$&$\ddd<0$ \\
					$-\aaa+\frac{\ccc}{2}$&$1-\frac{\bbb}{2}$&$\frac{\ccc}{2}$&$-\ddd+\frac{\bbb}{2}$&$\frac{\theta}{2}$&$1+\theta$&$2\theta$&$\frac13+\frac{\theta}{2}$&\textbf{$\surd$ \,\,\, Subcase} $4$  \\
					$\aaa+\frac{\ccc}{2}$&$1-\frac{\bbb}{2}$&$\frac{\ccc}{2}$&$-1+\ddd+\frac{\bbb}{2}$&$-\frac{\theta}{2}$&$1+\theta$&$2\theta$&$\frac23-\frac{\theta}{2}$&$\aaa<0$ \\
					$1-\aaa-\frac{\ccc}{2}$&$1-\frac{\bbb}{2}$&$\frac{\ccc}{2}$&$-1+\ddd+\frac{\bbb}{2}$&$1-\frac{3\theta}{2}$&$1+\theta$&$2\theta$&$\frac23-\frac{\theta}{2}$&$\aaa+\bbb+\ccc+\ddd>\frac83$  \\
					
			\end{tabular}}  
		\end{table*}

		\begin{table*}[htp]                        
			\centering            
			\resizebox{\textwidth}{30mm}{\begin{tabular}{cccc|cccc|c}	 
					
					$\frac{1}{2}-\frac{\theta}{2}$&$\frac{\theta}{2}$&$\frac{1}{6}$&$\theta$&$\aaa$&$\bbb$&$\ccc$&$\ddd$&$\ccc+2\ddd>1,\bbb+\ddd<2$ (Lemma $\ref{geometry4}'$\&\ref{proposition-3}) \\
					\hline 
					$\aaa-\frac{\ccc}{2}$&$1-\frac{\bbb}{2}$&$\frac{\ccc}{2}$&$\ddd-\frac{\bbb}{2}$&$\frac23-\frac{\theta}{2}$&$2-\theta$&$\frac13$&$1+\frac{\theta}{2}$&$\ddd>1$  \\
					$-\aaa+\frac{\ccc}{2}$&$1-\frac{\bbb}{2}$&$\frac{\ccc}{2}$&$1+\ddd-\frac{\bbb}{2}$&$-\frac13+\frac{\theta}{2}$&$2-\theta$&$\frac13$&$\frac{\theta}{2}$&$\aaa+\bbb+\ccc+\ddd=2$\\
					$-\aaa+\frac{\ccc}{2}$&$1-\frac{\bbb}{2}$&$\frac{\ccc}{2}$&$-\ddd+\frac{\bbb}{2}$&$-\frac13+\frac{\theta}{2}$&$2-\theta$&$\frac13$&$1-\frac{3\theta}{2}$&$\aaa<0$  \\
					$\aaa+\frac{\ccc}{2}$&$1-\frac{\bbb}{2}$&$\frac{\ccc}{2}$&$-1+\ddd+\frac{\bbb}{2}$&$\frac13-\frac{\theta}{2}$&$2-\theta$&$\frac13$&$\frac{3\theta}{2}$&$\bbb+\ddd>2$  \\
					$1-\aaa-\frac{\ccc}{2}$&$1-\frac{\bbb}{2}$&$\frac{\ccc}{2}$&$-1+\ddd+\frac{\bbb}{2}$&$\frac13+\frac{\theta}{2}$&$2-\theta$&$\frac13$&$\frac{3\theta}{2}$&$\aaa+\bbb+\ccc+\ddd>\frac83$ \\
					\multicolumn{9}{c}{ }\\
					$\frac{1}{2}-\frac{\theta}{2}$&$\frac{\theta}{2}$&$\theta$&$\frac{1}{6}$&$\aaa$&$\bbb$&$\ccc$&$\ddd$&$\ccc+2\ddd>1,\bbb+\ddd<2$ (Lemma $\ref{geometry4}'$\&\ref{proposition-3}) \\
					\hline 
					$\aaa-\frac{\ccc}{2}$&$1-\frac{\bbb}{2}$&$\frac{\ccc}{2}$&$\ddd-\frac{\bbb}{2}$&$\frac12+\frac{\theta}{2}$&$2-\theta$&$2\theta$&$\frac76-\frac{\theta}{2}$&$\aaa+\bbb+\ccc+\ddd>\frac83$ \\ 
					$-\aaa+\frac{\ccc}{2}$&$1-\frac{\bbb}{2}$&$\frac{\ccc}{2}$&$1+\ddd-\frac{\bbb}{2}$&$-\frac12+\frac{3\theta}{2}$&$2-\theta$&$2\theta$&$\frac16-\frac{\theta}{2}$&$\ccc+2\ddd<1$ \\
					$-\aaa+\frac{\ccc}{2}$&$1-\frac{\bbb}{2}$&$\frac{\ccc}{2}$&$-\ddd+\frac{\bbb}{2}$&$-\frac12+\frac{3\theta}{2}$&$2-\theta$&$2\theta$&$\frac56-\frac{\theta}{2}$&$\aaa>0 \Rightarrow \theta>\frac13$ but \\
					&&&&&&&& $\aaa+\cdots+\ddd\le\frac83 \Rightarrow \theta\le\frac16$\\
					$\aaa+\frac{\ccc}{2}$&$1-\frac{\bbb}{2}$&$\frac{\ccc}{2}$&$-1+\ddd+\frac{\bbb}{2}$&$\frac12-\frac{3\theta}{2}$&$2-\theta$&$2\theta$&$\frac16+\frac{\theta}{2}$&$\aaa>0 \Rightarrow \theta<\frac13$ but \\
					&&&&&&&& $\bbb+\ddd<2 \Rightarrow \theta>\frac13$\\
					$1-\aaa-\frac{\ccc}{2}$&$1-\frac{\bbb}{2}$&$\frac{\ccc}{2}$&$-1+\ddd+\frac{\bbb}{2}$&$\frac12-\frac{\theta}{2}$&$2-\theta$&$2\theta$&$\frac16+\frac{\theta}{2}$&$\aaa+\bbb+\ccc+\ddd>\frac83$  \\
					
			\end{tabular}} 
		\caption{Case $\{x_1, x_2\}=\{ \frac{1}{6},\theta\}$ and $\{x_3,x_4\}=\{\frac{\theta}{2},\frac{1}{2}-\frac{\theta}{2}\}$ or $\{x_3,x_4\}=\{\frac{1}{6},\theta\}$ and $\{x_1,x_2\}=\{\frac{\theta}{2},\frac{1}{2}-\frac{\theta}{2}\}$, for some $0<\theta\le\frac{1}{2}$. }\label{Tab-3.2}
		\end{table*}
								
		In Table \ref{Tab-3.2}  we list each of these $5 \times 8=40$ options.
		It turns out $36$ of these options are dismissed by Lemma \ref{anglesum}, Lemma $\ref{geometry4}'$ and Lemma \ref{proposition-3}. We list the corresponding details in the right hand column. 
		The remaining $4$ options are summarized as follows:
		
		$1$. $\aaa=-\frac16+\frac{\ccc}{2},\quad\bbb=2\ccc,\,\,\quad\quad\ddd=-\frac12+\frac{3\ccc}{2},\quad\frac{8}{15}<\ccc\le\frac23$;		
		
		$2$.  $\bbb=\frac23+2\aaa,\,\,\quad\ccc=\frac{1}{3},\,\,\,\,\,\quad\quad\ddd=3\aaa,\quad\quad\quad\,\,\,\,\,\frac16<\aaa\le\frac{5}{18}$;
		
		$3$.
		$\aaa=\frac{3\ccc}{4},\,\,\quad\quad\quad\bbb=1+\frac{\ccc}{2},\quad\ddd=\frac{2}{3}+\frac{\ccc}{4},\quad\quad\,\frac{2}{15}<\ccc\le\frac25$;
		
		$4$.
		$\aaa=\frac{\ccc}{4},\,\,\,\,\quad\quad\quad\bbb=1+\frac{\ccc}{2},\quad\ddd=\frac{1}{3}+\frac{\ccc}{4},\quad\quad\,\,\,\frac{1}{3}<\ccc\le\frac23$.
		
		\subsubsection*{Subcase $\aaa=-\frac16+\frac{\ccc}{2},\bbb=2\ccc,\ddd=-\frac12+\frac{3\ccc}{2},\frac{8}{15}<\ccc\le\frac23$} \label{discrete-2}
		
		By the angle values and Parity Lemma, only $\aaa\bbb\ddd$, $\bbb\ddd^2$ and $\ccc^3$ can be degree $3$ vertices. If $\bbb\ddd^2$ is a vertex, we have $\aaa=\frac{2}{15}$, $\bbb=\frac65$, $\ccc=\frac35$, $\ddd=\frac25$. Then $\bbb\cdots=\bbb\ddd^2$, $\aaa^3\bbb\ddd$ or $\aaa^6\bbb$. So $\#\aaa+\#\ddd>2\#\bbb=2f$, contradicting Balance Lemma \ref{balance}.  So $\aaa\bbb\ddd$ or $\ccc^3$ is a vertex. Both cases imply  $\aaa=\frac16$, $\bbb=\frac43$, $\ccc=\frac23$, $\ddd=\frac12$, and $f=6$. This implies all vertices have degree $3$.  There is only one solution satisfying Balance Lemma \ref{balance}: $\{6\aaa\bbb\ddd,2\ccc^3\}$, and it gives a $2$-layer earth map tiling by Lemma \ref{proposition-7}. This also gives Case $(1,8,4,3)/6$ in Table \ref{Tab-1.1} (see Remark \ref{remark23}).
		
		\subsubsection*{Subcase $\bbb=\frac23+2\aaa,\ccc=\frac{1}{3},\ddd=3\aaa,\frac16<\aaa\le\frac{5}{18}$} \label{discrete-4}
		By $R(\bbb\ddd\cdots)<3\aaa,\bbb,\ddd$, Parity Lemma and Lemma \ref{proposition-3}, we get $\aaa\bbb\ddd$ is a vertex. This implies $\aaa=\frac29,\bbb=\frac{10}{9},\ccc=\frac13,\ddd=\frac23$. Then $\bbb\cdots=\aaa\bbb\ddd,\aaa^4\bbb$. By $\#\bbb=\#\aaa$, we have  $\bbb\cdots=\aaa\bbb\ddd$. There is only one solution satisfying Balance Lemma \ref{balance}: $\{12\aaa\bbb\ddd,2\ccc^6\}$, and it gives a $2$-layer earth map tiling by Lemma \ref{proposition-7}. This is Case $(2,10,3,6)/9$ in Table \ref{Tab-1.1}. 
		
		\subsubsection*{Subcase $\aaa=\frac{3\ccc}{4},\bbb=1+\frac{\ccc}{2},\ddd=\frac{2}{3}+\frac{\ccc}{4},\frac{2}{15}<\ccc\le\frac25$} 
		
		By $R(\bbb\ddd\cdots)<3\aaa,\bbb,\ddd$, Parity Lemma and Lemma \ref{proposition-3}, we get $\aaa\bbb\ddd$ is a vertex.  This implies $\aaa=\frac16,\bbb=\frac{10}{9},\ccc=\frac29,\ddd=\frac{13}{18}$. Then  $\aaa+\ddd=4\ccc$ and $\bbb=5\ccc$. Then we get $f=18$. By the angle values and Parity Lemma, we get \[\text{AVC}\subset\{\aaa\bbb\ddd,\aaa^2\ccc\ddd^2,\bbb\ccc^4,\aaa^4\bbb\ccc,\aaa\ccc^5\ddd,\aaa^5\ccc^2\ddd,\ccc^9,\aaa^4\ccc^6,\aaa^8\ccc^3,\aaa^{12}\}.\] By $\#\ddd=\#\aaa$, we have $\aaa\cdots=\ddd\cdots=\aaa\bbb\ddd,\aaa^2\ccc\ddd^2$ or $\aaa\ccc^5\ddd$. Therefore,  \[\text{AVC}\subset\{\aaa\bbb\ddd,\aaa^2\ccc\ddd^2,\bbb\ccc^4,\aaa\ccc^5\ddd,\ccc^9\}.\] We will discuss all possible vertices containing $\ccc$ in any tiling as follows.
				
		If $\ccc^9$ appears, the tiling is a $2$-layer earth map tiling by Lemma \ref{proposition-7}. This is Case $(3,20,4,13)/18$ in Table \ref{Tab-1.1}. \label{discrete-12}	
		
		If $\aaa\ccc^5\ddd$ appears ($\ccc^9$ never appears), then
		$\aaa\ccc^5\ddd=\thin^{\bbb}\aaa_1^{\ddd}\thick^{\aaa}\ddd_2^{\ccc}\thin^{\bbb}\ccc_3^{\ddd}\thin\cdots$ determines $T_1,T_2,\dots,T_7$ in Fig.\,\ref{fig 3-2}. Then $\bbb_4\ddd_3\cdots=\aaa_{8}\bbb_4\ddd_3$ determines $T_{8}$. Similarly, we can determine $T_{9},T_{10},T_{11},T_{12}$. Then $\bbb_3\ccc_2\cdots=\bbb_3\ccc_2\ccc_{13}\ccc_{14}\ccc_{15}$. By $\ccc_{13}$, we get $\bbb_2\cdots=\bbb_2\ddd_{13}\cdots$ which determines $T_{13}$. Similarly, we can determine $T_{14},T_{15}$. Then $\bbb_2\ddd_{13}\cdots=\aaa_{16}\bbb_2\ddd_{13}$ and $\aaa_2\ddd_1\cdots=\aaa_2\bbb_{16}\ddd_1$ determine $T_{16}$. Similarly, we can determine $T_{17},T_{18}$. This tiling is exactly the second flip modification in Fig.\,\ref{flip1}.	\label{discrete-13}
		
		\begin{figure}[htp]
			\centering
			\begin{tikzpicture}[>=latex,scale=0.75] 
				
				\draw [line width=4pt, gray!50](0,0)--(0,-2)--(1,-3)--(1,-5)
				(6,0)--(6,-2)--(7,-3)--(7,-5);
				
				\foreach \a in {0,1,2,3,4}
				{
					\begin{scope}[xshift=2*\a cm] 
						\draw (6,0)--(6,-2)--(7,-3)--(7,-5)
						(8,0)--(8,-2)--(9,-3)--(9,-5);
						\draw[line width=1.5] (7,-3)--(8,-2);
						
						\node at (7,0){\small $\ccc$};
						\node at (6.3,-1.9){\small $\bbb$};
						\node at (7.7,-1.9){\small $\ddd$};
						\node at (7,-2.6){\small $\aaa$};
						\node at (8,-2.4){\small $\aaa$};
						\node at (8,-5){\small $\ccc$};
						\node at (7.3,-3.1){\small $\ddd$};
						\node at (8.7,-3.1){\small $\bbb$};
					\end{scope}
				}
				
				\draw (0,0)--(0,-2)--(1,-3)--(1,-5)
				(1,-3)--(2,-1.5)--(4,-1)--(6,-2)--(6,0)
				(1,-3)--(3,-2.5)--(4,-2.5)--(6,-2)
				(1,-3)--(3,-4)--(5,-3.5)--(6,-2)--(7,-3)--(7,-5);
				\draw[line width=1.5] (2,0)--(2,-1.5)
				(4,-1)--(3,-2.5)
				(3,-4)--(4,-2.5)
				(5,-3.5)--(5,-5);
				\node at (1,0){\small $\aaa$};
				\node at (1.7,-1.5){\small $\ddd$};
				\node at (0.3,-1.8){\small $\bbb$};
				\node at (1,-2.6){\small $\ccc$};
				\node at (4,0){\small $\ddd$};
				\node at (2.3,-1.2){\small $\aaa$};
				\node at (2.1,-1.8){\small $\bbb$};
				\node at (1.5,-2.65){\small $\ccc$};
				\node at (1.6,-3.1){\small $\ccc$};
				\node at (1.2,-3.4){\small $\ccc$};
				\node at (4,-0.65){\small $\bbb$};
				\node at (3.5,-1.35){\small $\aaa$};
				\node at (4,-1.35){\small $\ddd$};
				\node at (2.8,-2.3){\small $\ddd$};
				\node at (3.4,-2.3){\small $\aaa$};
				\node at (4,-2.25){\small $\bbb$};
				\node at (3,-2.8){\small $\bbb$};
				\node at (3.6,-2.7){\small $\aaa$};
				\node at (4.2,-2.7){\small $\ddd$};
				\node at (2.9,-3.7){\small $\ddd$};
				\node at (3.5,-3.7){\small $\aaa$};
				\node at (3,-4.3){\small $\bbb$};
				\node at (4.9,-3.3){\small $\bbb$};
				\node at (4.7,-3.8){\small $\aaa$};
				\node at (5.3,-3.6){\small $\ddd$};
				\node at (3,-5){\small $\ddd$};
				\node at (6.75,-3.3){\small $\bbb$};
				\node at (6,-5){\small $\aaa$};
				\node at (5.8,-1.7){\small $\ccc$};
				\node at (5.4,-1.95){\small $\ccc$};
				\node at (5.5,-2.4){\small $\ccc$};
				\node at (6,-2.4){\small $\ccc$};
				
				\node[draw,shape=circle, inner sep=0.5] at (1,-1) {\small $1$};
				\node[draw,shape=circle, inner sep=0.5] at (5,-1) {\small $2$};
				\node[draw,shape=circle, inner sep=0.5] at (7,-1) {\small $3$};
				\node[draw,shape=circle, inner sep=0.5] at (9,-1) {\small $4$};
				\node[draw,shape=circle, inner sep=0.5] at (11,-1) {\small $5$};
				\node[draw,shape=circle, inner sep=0.5] at (13,-1) {\small $6$};
				\node[draw,shape=circle, inner sep=0.5] at (15,-1) {\small $7$};
				\node[draw,shape=circle, inner sep=0.5] at (8,-4) {\small $8$};
				\node[draw,shape=circle, inner sep=0.5] at (10,-4) {\small $9$};
				\node[draw,shape=circle, inner sep=0.5] at (12,-4) {\footnotesize $10$};
				\node[draw,shape=circle, inner sep=0.5] at (14,-4) {\footnotesize $11$};
				\node[draw,shape=circle, inner sep=0.5] at (16,-4) {\footnotesize $12$};
				\node[draw,shape=circle, inner sep=0.5] at (4.7,-1.8) {\footnotesize $13$};
				\node[draw,shape=circle, inner sep=0.5] at (4.7,-2.8) {\footnotesize $14$};
				\node[draw,shape=circle, inner sep=0.5] at (6,-4) {\footnotesize $15$};
				\node[draw,shape=circle, inner sep=0.5] at (2.1,-2.35) {\footnotesize $16$};
				\node[draw,shape=circle, inner sep=0.5] at (2.3,-3.1) {\footnotesize $17$};
				\node[draw,shape=circle, inner sep=0.5] at (2,-4) {\footnotesize $18$};

			\end{tikzpicture}
			\caption{$T(16\aaa\bbb\ddd,2\bbb\ccc^4,2\aaa\ccc^5\ddd)$.} \label{fig 3-2}
		\end{figure}
		
		If $\aaa^2\ccc\ddd^2$ appears  ($\aaa\ccc^5\ddd$, $\ccc^9$ never appear), 
	    it has only two possible AAD since there is no vertex $\bbb^2\cdots$ by $\bbb>1$. In the first picture of Fig.\,\ref{fig 3-3}, $\aaa^2\ccc\ddd^2=\thin^{\bbb}\aaa_1^{\ddd}\thick^{\ddd}\aaa_2^{\bbb}\thin^{\ccc}\ddd_3^{\aaa}\thick^{\aaa}\ddd_4^{\ccc}\thin^{\bbb}\ccc_5^{\ddd}\thin$ determines $T_1,T_2,T_3,T_4,T_5$. Then $\bbb_1\ddd_5\cdots=\aaa_6\bbb_1\ddd_5$ determines $T_6$. We have $\bbb_6\ccc_1\cdots=\bbb_6\ccc_1\ccc_7\ccc^2$. By $\ccc_7$, we have $\bbb_7\ddd_1\ddd_2\cdots$ or $\ddd_1\ddd_2\ddd_7\cdots$, contradicting the AVC.	
		In the second picture of Fig.\,\ref{fig 3-3}, $\aaa^2\ccc\ddd^2=\thin^{\bbb}\aaa_1^{\ddd}\thick^{\aaa}\ddd_2^{\ccc}\thin^{\bbb}\aaa_3^{\ddd}\thick^{\aaa}\ddd_4^{\ccc}\thin^{\bbb}\ccc_5^{\ddd}\thin$ determines $T_1,T_2,T_3,T_4,T_5$. Then $\bbb_1\ddd_5\cdots=\aaa_6\bbb_1\ddd_5$ determines $T_6$. Then $\bbb_5\ccc_4\cdots=\bbb_5\ccc_4\ccc_7\ccc_8\ccc_9$. By $\ccc_7$, we get $\bbb_4\cdots=\bbb_4\ddd_7\cdots$ which determines $T_7$. Similarly, we can determine $T_{8},T_{9}$. Then $\bbb_4\ddd_{7}\cdots=\aaa_{10}\bbb_4\ddd_{7},\aaa_4\ddd_3\cdots=\aaa_4\bbb_{10}\ddd_3$ determine $T_{10}$. Similarly, we can determine $T_{11}$, $T_{12}$. After repeating the process one more time, we can determine $T_{13}$, $T_{14}$, $\dots$, $T_{18}$.
		This tiling can also be obtained by applying the second flip modification in Fig.\,\ref{flip1} two times.	\label{discrete-14}	
		\begin{figure}[htp]
			\centering
			\begin{tikzpicture}[>=latex,scale=0.44*1.45] 
							
				\draw (0,0.5/1.45)--(2/1.45,-2/1.45)--(2/1.45,-4/1.45)--(3/1.45,-6/1.45)--(1/1.45,-6/1.45)--(0,-8/1.45)
				(2/1.45,-2/1.45)--(0,-2/1.45)--(0,-4/1.45)--(1/1.45,-6/1.45)
				(0,0.5/1.45)--(-2/1.45,-2/1.45)--(-2/1.45,-4/1.45)--(0,-4/1.45)
				(-2/1.45,-4/1.45)--(-3/1.45,-6/1.45)--(-1/1.45,-6/1.45)--(0,-8/1.45);
				\draw[line width=1.5](-2/1.45,-2/1.45)--(0,-2/1.45)
				(-1/1.45,-6/1.45)--(0,-4/1.45)--(2/1.45,-4/1.45);
				\draw[dotted] (2/1.45,-2/1.45)--(4/1.45,-2/1.45)--(4/1.45,-4/1.45)--(2/1.45,-4/1.45);
				
				\node at (0,-0.2/1.45){\small $\ccc$};
				\node at (-1.2/1.45,-1.6/1.45){\small $\ddd$};
				\node at (0,-1.7/1.45){\small $\aaa$};
				\node at (1/1.45,-1.6/1.45){\small $\bbb$};
				\node at (-0.4/1.45,-2.5/1.45){\small $\ddd$};
				\node at (0.4/1.45,-2.5/1.45){\small $\bbb$};
				\node at (-1.6/1.45,-2.5/1.45){\small $\aaa$};
				\node at (1.6/1.45,-2.5/1.45){\small $\ccc$};
				\node at (-1.6/1.45,-3.5/1.45){\small $\bbb$};
				\node at (1.6/1.45,-3.5/1.45){\small $\ddd$};
				\node at (-0.4/1.45,-3.5/1.45){\small $\ccc$};
				\node at (0.4/1.45,-3.5/1.45){\small $\aaa$};
				\node at (0.6/1.45,-4.4/1.45){\small $\aaa$};
				\node at (-0.6/1.45,-4.4/1.45){\small $\ddd$};
				\node at (0,-4.7/1.45){\small $\ddd$};
				\node at (1.8/1.45,-4.4/1.45){\small $\ddd$};
				\node at (-1.8/1.45,-4.4/1.45){\small $\ccc$};
				\node at (1.2/1.45,-5.6/1.45){\small $\bbb$};
				\node at (-1.2/1.45,-5.6/1.45){\small $\aaa$};
				\node at (2.4/1.45,-5.6/1.45){\small $\ccc$};
				\node at (-2.4/1.45,-5.6/1.45){\small $\bbb$};
				\node at (0.6/1.45,-6/1.45){\small $\ccc$};
				\node at (-0.6/1.45,-6/1.45){\small $\aaa$};
				\node at (0,-7.2/1.45){\small $\bbb$};
				\node at (2.4/1.45,-2.5/1.45){\small $\ccc$};

				\node[draw,shape=circle, inner sep=0.5] at (1/1.45,-3/1.45) {\small $1$};
				\node[draw,shape=circle, inner sep=0.5] at (1.3/1.45,-4.8/1.45) {\small $2$};
				\node[draw,shape=circle, inner sep=0.5] at (0,-5.6/1.45) {\small $3$};
				\node[draw,shape=circle, inner sep=0.5] at (-1.3/1.45,-4.8/1.45) {\small $4$};
				\node[draw,shape=circle, inner sep=0.5] at (-1/1.45,-3/1.45) {\small $5$};
				\node[draw,shape=circle, inner sep=0.5] at (0,-1/1.45) {\small $6$};
				\node[draw,shape=circle, inner sep=0.5] at (3/1.45,-3/1.45) {\small $7$};
		
			\begin{scope}[xshift=4 cm] 
				
				\draw [line width=4pt, gray!50](0,0)--(0,-2)--(1,-3)--(1,-5)
				(6,0)--(6,-2)--(7,-3)--(7,-5)
				(12,0)--(12,-2)--(13,-3)--(13,-5);
				\foreach \a in {0,1}
				{
					\begin{scope}[xshift=6*\a cm] 
						\draw (0,0)--(0,-2)--(1,-3)--(1,-5)
						(1,-3)--(2,-1.5)--(4,-1)--(6,-2)--(6,0)
						(1,-3)--(3,-2.5)--(4,-2.5)--(6,-2)
						(1,-3)--(3,-4)--(5,-3.5)--(6,-2)--(7,-3)--(7,-5);
						\draw[line width=1.5] (2,0)--(2,-1.5)
						(4,-1)--(3,-2.5)
						(3,-4)--(4,-2.5)
						(5,-3.5)--(5,-5);
						\node at (1,0){\small $\aaa$};
						\node at (1.7,-1.5){\small $\ddd$};
						\node at (0.3,-1.8){\small $\bbb$};
						\node at (1,-2.6){\small $\ccc$};
						\node at (4,0){\small $\ddd$};
						\node at (2.3,-1.2){\small $\aaa$};
						\node at (2.2,-1.7){\small $\bbb$};
						\node at (1.5,-2.65){\small $\ccc$};
						\node at (1.6,-3.1){\small $\ccc$};
						\node at (1.2,-3.4){\small $\ccc$};
						\node at (4,-0.65){\small $\bbb$};
						\node at (3.5,-1.35){\small $\aaa$};
						\node at (4,-1.35){\small $\ddd$};
						\node at (2.8,-2.3){\small $\ddd$};
						\node at (3.4,-2.3){\small $\aaa$};
						\node at (4,-2.25){\small $\bbb$};
						\node at (3,-2.8){\small $\bbb$};
						\node at (3.6,-2.7){\small $\aaa$};
						\node at (4.2,-2.7){\small $\ddd$};
						\node at (2.9,-3.7){\small $\ddd$};
						\node at (3.5,-3.7){\small $\aaa$};
						\node at (3,-4.3){\small $\bbb$};
						\node at (4.9,-3.3){\small $\bbb$};
						\node at (4.7,-3.8){\small $\aaa$};
						\node at (5.3,-3.6){\small $\ddd$};
						\node at (3,-5){\small $\ddd$};
						\node at (6.75,-3.3){\small $\bbb$};
						\node at (6,-5){\small $\aaa$};
						\node at (5.8,-1.7){\small $\ccc$};
						\node at (5.4,-1.95){\small $\ccc$};
						\node at (5.5,-2.4){\small $\ccc$};
						\node at (6,-2.4){\small $\ccc$};
					\end{scope}
				}
				\draw (6+6,0)--(6+6,-2)--(7+6,-3)--(7+6,-5)
				(8+6,0)--(8+6,-2)--(9+6,-3)--(9+6,-5);
				\draw[line width=1.5] (7+6,-3)--(8+6,-2);
				
				\node at (7+6,0){\small $\ccc$};
				\node at (6.3+6,-1.9){\small $\bbb$};
				\node at (7.7+6,-1.9){\small $\ddd$};
				\node at (7+6,-2.6){\small $\aaa$};
				\node at (8+6,-2.4){\small $\aaa$};
				\node at (8+6,-5){\small $\ccc$};
				\node at (7.3+6,-3.1){\small $\ddd$};
				\node at (8.7+6,-3.1){\small $\bbb$};

				\node[draw,shape=circle, inner sep=0.5] at (1,-1) {\small $1$};
				\node[draw,shape=circle, inner sep=0.5] at (5,-0.5) {\small $2$};
				\node[draw,shape=circle, inner sep=0.5] at (7,-1) {\small $3$};
				\node[draw,shape=circle, inner sep=0.5] at (11,-0.5) {\small $4$};
				\node[draw,shape=circle, inner sep=0.5] at (13,-1) {\small $5$};
				\node[draw,shape=circle, inner sep=0.5] at (14,-4) {\small $6$};
				\node[draw,shape=circle, inner sep=0.5] at (10.7,-1.8) {\small $7$};
				\node[draw,shape=circle, inner sep=0.5] at (10.7,-2.8) {\small $8$};
				\node[draw,shape=circle, inner sep=0.5] at (12,-4) {\small $9$};
				\node[draw,shape=circle, inner sep=0.5] at (8.2,-2.3) {\footnotesize $10$};
				\node[draw,shape=circle, inner sep=0.5] at (8.3,-3.1) {\footnotesize $11$};
				\node[draw,shape=circle, inner sep=0.5] at (8,-4.3) {\footnotesize $12$};
				\node[draw,shape=circle, inner sep=0.5] at (4.7,-1.85) {\footnotesize $13$};
				\node[draw,shape=circle, inner sep=0.5] at (4.8,-2.7) {\footnotesize $14$};
				\node[draw,shape=circle, inner sep=0.5] at (6,-4) {\footnotesize $15$};
				\node[draw,shape=circle, inner sep=0.5] at (2.2,-2.3) {\footnotesize $16$};
				\node[draw,shape=circle, inner sep=0.5] at (2.3,-3.1) {\footnotesize $17$};
				\node[draw,shape=circle, inner sep=0.5] at (2,-4.3) {\footnotesize $18$};
			\end{scope}
			\end{tikzpicture}
			\caption{ $T(14\aaa\bbb\ddd,2\aaa^2\ccc\ddd^2,4\bbb\ccc^4)$.} \label{fig 3-3}
		\end{figure}
		
		If $\text{AVC}\sub\{\aaa\bbb\ddd,\bbb\ccc^4\}$, there is no solution satisfying Balance Lemma \ref{balance}.   
		
		\subsubsection*{Subcase $\aaa=\frac{\ccc}{4},\bbb=1+\frac{\ccc}{2},\ddd=\frac{1}{3}+\frac{\ccc}{4},\frac{1}{3}<\ccc\le\frac23$} \label{discrete-2''}
		
		By $R(\bbb\ddd\cdots)<5\aaa,\bbb,\ddd$, Parity Lemma and Lemma \ref{proposition-3}, we deduce that $\aaa\bbb\ddd$ or $\aaa^3\bbb\ddd$ is a vertex.  If $\aaa\bbb\ddd$ is a vertex, then $\aaa=\frac16,\bbb=\frac43,\ccc=\frac23,\ddd=\frac12$. There is only one solution satisfying Balance Lemma \ref{balance}: $\{6\aaa\bbb\ddd,2\ccc^3\}$, and it gives a $2$-layer earth map tiling by Lemma \ref{proposition-7}. This also gives Case $(1,8,4,3)/6$ in Table \ref{Tab-1.1}.	
		If $\aaa^3\ccc\ddd$ is a vertex, then $\aaa=\frac19,\bbb=\frac{11}{9},\ccc=\frac49,\ddd=\frac{4}{9}$, which does not admit any any degree $3$ vertex, a contradiction.

		\subsection*{Case $4$: $\{ x_1,x_2,x_3,x_4\}$ are in Table  \ref{Tab-3.1}.}		
		There are $8\times5\times15=600$ subcases to consider, but most are ruled out by violating $0<\aaa,\ccc,\ddd<1,1<\bbb<2$, $f$ being even integer or $\bbb+\ddd<2$. Only ten subcases are left in Table \ref{Tab-3.5}. But six of them are ruled out by not admitting any degree $3$ vertex. Four remaining subcases are boxed.
		
		\begin{table}[htp]                        
			\centering            
			\begin{minipage}{0.49\textwidth}
				\scriptsize{\begin{tabular}{ccc}				
				$(\aaa,\bbb,\ccc,\ddd)$&$f$& \\
				\hline 
				(5,32,18,11)/30&20&no degree $3$ vertex\\
				\hline
				\multicolumn{1}{|c}{(5,32,6,23)/30}&20&\multicolumn{1}{c|}{}\\
				\multicolumn{1}{|c}{(1,42,4,17)/30}&30&\multicolumn{1}{c|}{}\\
				\hline
				(1,17,9,4)/15&60&no degree $3$ vertex\\
				\hline
				\multicolumn{1}{|c}{(1,21,5,8)/15}&12&\multicolumn{1}{c|}{}\\
				\hline
			\end{tabular}}
		\end{minipage}	
	     \raisebox{0em}{\begin{minipage}{0.49\textwidth}
	     		\scriptsize{\begin{tabular}{ccc}
	     		$(\aaa,\bbb,\ccc,\ddd)$&$f$& \\
				\hline
				(13,66,32,29)/60&12&no degree $3$ vertex\\
				\hline
				\multicolumn{1}{|c}{(5,32,14,13)/30}&30&\multicolumn{1}{c|}{$\bbb\cdots=\bbb\ccc^2$ or $\aaa^3\bbb\ddd$}\\
				\hline
				(3,32,22,7)/30&30&no degree $3$ vertex\\
				(1,19,3,8)/15&60&no degree $3$ vertex\\
				(7,66,8,49)/60&24&no degree $3$ vertex
				
			\end{tabular}}
		\end{minipage}}
	    \caption{Ten solutions induced from Table \ref{Tab-3.1}.}\label{Tab-3.5}  
	   \end{table}

		\subsubsection*{Subcase $(5,32,6,23)/30$}  \label{discrete-9}
		
		By the angle values and Parity Lemma, we get
		$\bbb\cdots=\aaa\bbb\ddd$ or $\aaa^2\bbb\ccc^3$. By $\#\bbb=\#\aaa$, we get $\bbb\cdots=\aaa\bbb\ddd$, which determines a $2$-layer earth map tiling $T(20\aaa\bbb\ddd,2\ccc^{10})$ in Table \ref{Tab-1.1} by Lemma \ref{proposition-7}. 
		
		\subsubsection*{Subcase  $(1,42,4,17)/30$} \label{discrete-11}
		
		By the angle values and Parity Lemma, we get
		$\bbb\cdots=\aaa\bbb\ddd$, $\aaa^2\bbb\ccc^4$, $\aaa^6\bbb\ccc^3$, $\aaa^{10}\bbb\ccc^2$, $\aaa^{14}\bbb\ccc$ or  $\aaa^{18}\bbb$. By $\#\bbb=\#\aaa$, we get $\bbb\cdots=\aaa\bbb\ddd$, which determines a $2$-layer earth map tiling $T(30\aaa\bbb\ddd,2\ccc^{15})$ in Table \ref{Tab-1.1} by Lemma \ref{proposition-7}.
		
		\subsubsection*{Subcase  $(1,21,5,8)/15$} \label{discrete-5}
		
		By the angle values and Parity Lemma, we get
		$\bbb\cdots=\aaa\bbb\ddd$ or $\aaa^4\bbb\ccc$. By $\#\bbb=\#\aaa$, we get $\bbb\cdots=\aaa\bbb\ddd$, which determines a $2$-layer earth map tiling $T(12\aaa\bbb\ddd,2\ccc^{6})$ in Table \ref{Tab-1.1} by Lemma \ref{proposition-7}.
		
		\subsubsection*{Subcase  $(5,32,14,13)/30$}
		
		By the angle values and Parity Lemma, we get
		$\bbb\cdots=\bbb\ccc^2$ or $\aaa^3\bbb\ddd$. There is no solution satisfying Balance Lemma \ref{balance}.

	   

	\section{Degenerate case $\bbb=1$}
	\label{sec-degenerate-b}	
	
    If $\bbb=1$, the quadrilateral degenerates to an isosceles triangle in Fig.\,\ref{fig 7-pi}.
	
	\begin{figure}[htp]
		\centering
		\begin{tikzpicture}[>=latex,scale=0.7]

			\draw (0,0)--(4,0)
			(0,0)--(1.5,2);
			\draw[line width=1.5] (1.5,2)--(4,0);
			\node at (0.6,0.2){\small $\ccc$};
			\node at (1.5,1.6){\small $\ddd$};
			\node at (3.3,0.2){\small $\aaa$};
			\fill (2,0) circle (0.05);
			\node at (0.7,1.3){\small $a$};	  	\node at (2.8,1.3){\small $b$};		\node at (2,-0.3){\small $2a$};  		
		
	   \begin{scope}[xshift=3 cm]		
			\draw[line width=1.5] (5.26,1.8)--(6.8,0.2); 
			\node at (6,2.6){\small $\ccc$}; \node at (6,-0.8){\small $\ccc$};\node at (6.6,0){\small $\ddd$};
			\node at (6.7,0.8){\small $\aaa$}; \node at (5.6,1.9){\small $\ddd$}; \node at (5.4,1.2){\small $\aaa$};
			\fill (6.73,1.8) circle (0.05); \fill (5.21,0.2) circle (0.05);
			\draw (6,3) arc (45:-45:3); \draw (6,3) arc (135:225:3);
     \end{scope}		
		\end{tikzpicture}
		\caption{ $\bbb=1$ and the subcase $\aaa+\ddd=1$.} \label{fig 7-pi}
	\end{figure}
	
	By $\bbb=1$, we have $\aaa+\ccc+\ddd=(1+\frac{4}{f})$.
	By Lemma $\ref{geometry1}'$ and Lemma $\ref{geometry4}'$, we get $\ddd>\aaa$ and $\ccc+2\ddd>1$. By $a+b>2a$, we get $b>a$. This implies $\ccc>\aaa$. If  $\aaa\ge\frac12$, then $R(\bbb\cdots)=1\le2\aaa<2\ccc,2\ddd$. So $\bbb\cdots=\aaa^2\bbb$, contradicting Balance Lemma. We conclude that $\aaa<\frac12$ and $\aaa^2\bbb$ is never a vertex.
	
	If $\aaa\bbb\ddd$ is a vertex, we have $\aaa+\ddd=1,\bbb=1,\ccc=\frac 4f$, as shown in the second picture of Fig.\,\ref{fig 7-pi}. So $a=\frac13$, and we get $\aaa=\arctan( 2\tan\frac{2\pi}{f} )$ by the cosine law. This is equivalent to $\cos(\frac{\pi}{2}-\aaa-\frac{2\pi}{f})-3\cos(\frac{\pi}{2}-\aaa+\frac{2\pi}{f})=0$ by the product to sum formula. Then Theorem $6$ of Conway-Jones \cite{cj} implies that $\aaa$ is irrational for any even integer $f\ge6$. Thus this belongs to the irrational angle case in \cite{lpwx2}. Such quadrilaterals always admit $2$-layer earth map tilings  
	for any even integer $f\ge6$, together with their flip modifications when $f=4k$ as shown in Fig.\,\ref{fig 7-3}.

		\begin{figure}[htp]
			\centering
			\begin{tikzpicture}[>=latex,scale=0.48] 
				
				\draw (0,0)--(0,-8)
				(6,0)--(6,-8)
				(0,-6)--(2,-2)--(6,-2)--(4,-6)--(0,-6);
				\draw[line width=1.5] (2,0)--(2,-2)
				(4,-6)--(4,-8);
				\node at (1,0){\small $\aaa$};
				\node at (4,0){\small $\ddd$};
				\node at (0.3,-2){\small $\bbb$};
				\node at (1.7,-2){\small $\ddd$};
				\node at (2.3,-1.7){\small $\aaa$};
				\node at (4,-1.6){\small $\bbb$};
				\node at (5.7,-1.7){\small $\ccc$};
				\node at (0.25,-5.1){\small $\ccc$};
				\node at (5.75,-3.1){\small $\ccc$};
				\node at (0.3,-6.3){\small $\ccc$};
				\node at (2,-6.4){\small $\bbb$};
				\node at (3.7,-6.3){\small $\aaa$};
				\node at (4.3,-6){\small $\ddd$};
				\node at (5.7,-6){\small $\bbb$};
				\node at (2,-8){\small $\ddd$};
				\node at (5,-8){\small $\aaa$};
				
				\fill (4,-2) circle (0.04);
				\fill (2,-6) circle (0.04);

				\fill (1,-5.8) circle (0.05);
				\fill (0.8,-5.6) circle (0.05);
				\fill (0.6,-5.4) circle (0.05);
				
				\fill (1+4.4,-5.8+3.1) circle (0.05);
				\fill (0.8+4.4,-5.6+3.1) circle (0.05);
				\fill (0.6+4.4,-5.4+3.1) circle (0.05);
				
				\fill (0,-2) circle (0.04);

				\foreach \a in {0,1}
				{
					\begin{scope}[xshift=2*\a cm] 
						\draw (8,0)--(8,-8);
						\draw[line width=1.5] (6,-6)--(8,-2);
						\node at (7,0){\small $\ccc$};
						\node at (6.3,-2){\small $\bbb$};
						\node at (7.7,-2){\small $\ddd$};
						\node at (6.3,-4.8){\small $\aaa$};
						\node at (7.7,-3.2){\small $\aaa$};
						\node at (6.3,-6){\small $\ddd$};
						\node at (7.7,-6){\small $\bbb$};
						\node at (7,-8){\small $\ccc$};
						
					\end{scope}
				}

				\fill (10.8,-4) circle (0.05);
				\fill (11.1,-4) circle (0.05);
				\fill (11.4,-4) circle (0.05);
				
				\fill (10,-6) circle (0.04);
				\node at (5.5,-9.5){\small $T((4k-2)\aaa\bbb\ddd,2\aaa\ccc^{k}\ddd,2\bbb\ccc^{k})$};

		\begin{scope}[xshift=15 cm]
			\foreach \a in {0,1}
			{
				\begin{scope}[xshift=6*\a cm] 
					\draw (0,0)--(0,-8)
					(6,0)--(6,-8)
					(0,-6)--(2,-2)--(6,-2)--(4,-6)--(0,-6);
					\draw[line width=1.5] (2,0)--(2,-2)
					(4,-6)--(4,-8);
					\node at (1,0){\small $\aaa$};
					\node at (4,0){\small $\ddd$};
					\node at (0.3,-2){\small $\bbb$};
					\node at (1.7,-2){\small $\ddd$};
					\node at (2.3,-1.7){\small $\aaa$};
					\node at (4,-1.6){\small $\bbb$};
					\node at (5.7,-1.7){\small $\ccc$};
					\node at (0.25,-5.1){\small $\ccc$};
					\node at (5.75,-3.1){\small $\ccc$};
					\node at (0.3,-6.3){\small $\ccc$};
					\node at (2,-6.4){\small $\bbb$};
					\node at (3.7,-6.3){\small $\aaa$};
					\node at (4.3,-6){\small $\ddd$};
					\node at (5.7,-6){\small $\bbb$};
					\node at (2,-8){\small $\ddd$};
					\node at (5,-8){\small $\aaa$};
					
					\fill (4,-2) circle (0.04);
					\fill (2,-6) circle (0.04);

					\fill (1,-5.8) circle (0.05);
					\fill (0.8,-5.6) circle (0.05);
					\fill (0.6,-5.4) circle (0.05);
					
					\fill (1+4.4,-5.8+3.1) circle (0.05);
					\fill (0.8+4.4,-5.6+3.1) circle (0.05);
					\fill (0.6+4.4,-5.4+3.1) circle (0.05);
					
				\end{scope}
			}
			
			\fill (0,-2) circle (0.04);
			\fill (12,-6) circle (0.04);
			\node at (6,-9.5){\small $T((4k-4)\aaa\bbb\ddd,2\aaa^2\ddd^2,4\bbb\ccc^{k})$};
			
         \end{scope}
		\end{tikzpicture}
		\caption{  Two flips of $T(4k\aaa\bbb\ddd,2\ccc^{2k})$ if $\bbb=1$ and $\aaa+\ddd=1$.} \label{fig 7-3}
	\end{figure}
	
	If $\aaa\bbb\ddd$ is not a vertex, then we will find all tilings by discussing all possible $\bbb$-vertices. 
	If $\aaa^x\bbb(x\ge4)$ is a vertex, then its unique AAD $\thick\aaa^{\bbb}\thin^{\bbb}\aaa\thick\cdots$ at $\aaa^x\bbb(x\ge4)$  gives a vertex $\bbb^2\cdots$, contradicting $\bbb=1$.  Similarly, $\aaa^y\bbb\ddd(y\ge5)$, $\aaa^w\bbb\ccc(w\ge4)$, $\aaa^p\ccc^q(p\ge2)$ and $\aaa^z\ddd^2(z\ge4)$ are not vertices.	
			
	\subsubsection*{Subcase $\ccc>\ddd$}
	
	Then we have $\bbb>\ccc>\ddd>\aaa$. By $\ccc+2\ddd>1$, we get $\ccc>\frac13$. By the angle values and Parity Lemma, we get $\bbb\cdots=\bbb\ccc^2,\aaa^2\bbb\ccc,\bbb\ddd^x(x\ge2)$ or $\aaa^p\bbb\ddd^q(p\ge1,p+q\ge4)$. 

	If $\bbb\ccc^2$ is a vertex, then $\ccc=\frac12,\aaa+\ddd>\frac12,\frac14<\ddd<\frac12$. So we have $\bbb\cdots=\bbb\ccc^2$, $\aaa^2\bbb\ccc$ or $\aaa^3\bbb\ddd$. They all satisfy   $\#\aaa+\#\ccc\ge2\#\bbb$.
	If $\aaa^2\bbb\ccc$ or $\aaa^3\bbb\ddd$ is a vertex, then $\#\aaa+\#\ccc>2\#\bbb$, contradicting Balance Lemma \ref{balance}. If $\bbb\cdots=\bbb\ccc^2$, then $\#\ccc>\#\bbb$,  again a contradiction.

	Therefore, we have $\bbb\cdots=\aaa^2\bbb\ccc,\bbb\ddd^x$ or $\aaa^p\bbb\ddd^q$. They all satisfy   $\#\aaa+\#\ddd\ge2\#\bbb=2f$. There is only one solution satisfying Balance Lemma \ref{balance}: $\{\frac{f}{2}\bbb\ddd^2,\frac{f}{2}\aaa^2\bbb\ccc,2\ccc^{\frac f 4}\}$. This implies $\aaa=\frac12-\frac4f$, $\ccc=\frac8f$, $\ddd=\frac12$. By $1>\ccc>\ddd$, we get $8<f<16$, which do not satisfy \eqref{4-7} in Lemma \ref{geometry6}. We conclude that there is no tiling in this case.

	\subsubsection*{Subcase $\ccc<\ddd$}	
	
	Then we have $\bbb>\ddd>\ccc>\aaa$. By $\ccc+2\ddd>1$, we get $\ddd>\frac13$. By the angle values and Parity Lemma, we get $\bbb\cdots=\bbb\ddd^2,\aaa^3\bbb\ddd,\bbb\ccc^x(x\ge2)$ or $\aaa^p\bbb\ccc^q(p\ge2,q\ge1)$. 
	
	If $\bbb\ddd^2$ is a vertex, then $\ccc<\ddd=\frac12$, $\aaa+\ccc>\frac12$. 
	So  $\bbb\cdots=\bbb\ddd^2,\aaa^2\bbb\ccc,\aaa^3\bbb\ddd$ or $\bbb\ccc^y(y\ge3)$. 
	If $\aaa^3\bbb\ddd$ is a vertex, then $\bbb\cdots=\bbb\ddd^2$ or $\aaa^3\bbb\ddd$. So $\#\aaa+\#\ddd>2\#\bbb$, contradicting Balance Lemma \ref{balance}.  Similarly, $\bbb\ccc^y$ is not a vertex. So  $\bbb\cdots=\bbb\ddd^2$ or $\aaa^2\bbb\ccc$.  There is only one solution satisfying Balance Lemma \ref{balance}: $\{\frac{f}{2}\bbb\ddd^2,\frac{f}{2}\aaa^2\bbb\ccc,2\ccc^{\frac f 4}\}$.  We get $\aaa=\frac12-\frac{4}{f},\bbb=1,\ccc=\frac{8}{f},\ddd=\frac12$. By $\ccc<\ddd$, we get $f>16$. By \eqref{4-7} in Lemma \ref{geometry6}, we get $f=16$, a contradiction. 
		
	If $\bbb\ddd^2$ is not a vertex, we have $\bbb\cdots=\aaa^3\bbb\ddd,\bbb\ccc^x$ or $\aaa^p\bbb\ccc^q$.
	They all satisfy   $\#\aaa+\#\ccc\ge2\#\bbb$.
	If $\aaa^3\bbb\ddd$ or $\aaa^p\bbb\ccc^q$ is a vertex, then $\#\aaa+\#\ccc>2\#\bbb$, contradicting Balance Lemma \ref{balance}. If $\bbb\cdots=\bbb\ccc^x$, then $\#\ccc>\#\bbb$,  again a contradiction.
    
	\subsubsection*{Subcase $\ccc=\ddd$} 
	
	   By $\ccc+2\ddd>1$, we get $\ccc=\ddd>\frac13$. By the angle values and Parity Lemma, we get $\bbb\cdots=\bbb\ccc^2,\bbb\ddd^2,\aaa^2\bbb\ccc$ or $\aaa^3\bbb\ddd$. 
	   If $\bbb\ccc^2$ and $\bbb\ddd^2$ are not vertices, we have $\bbb\cdots=\aaa^2\bbb\ccc$ or $\aaa^3\bbb\ddd$, contradicting Balance Lemma \ref{balance}. Therefore, $\bbb\ccc^2$ or $\bbb\ddd^2$ is a vertex. So we get $\ccc=\ddd=\frac12$. Then we get $\aaa=\frac{4}{f}$. By $\ccc=\ddd=\frac12$, we get $b=2a=\frac12$. By the sine law $
	   \frac{\sin \aaa}{\sin a}=\frac{\sin \ccc}{\sin 2a}$, we have $\aaa=\frac14$. This implies $f=16$. By the angle values and Parity Lemma, we get  \[\text{AVC}\sub\{\bbb\ccc^2,\bbb\ddd^2,\aaa^2\bbb\ccc,\ccc^4,\ccc^2\ddd^2,\ddd^4,\aaa^2\ccc\ddd^2\}.\]	   
	   If $\aaa^2\ccc\ddd^2$ is a vertex, it has only two possible AAD. In the left of Fig.\,\ref{fig 2ac2d 1}, $\aaa^2\ccc\ddd^2=\thin^{\bbb}\aaa_1^{\ddd}\thick^{\ddd}\aaa_2^{\bbb}\thin^{\ccc}\ddd_3^{\aaa}\thick^{\aaa}\ddd_4^{\ccc}\thin^{\bbb}\ccc_5^{\ddd}\thin$	determines $T_1,T_2$, $T_3,T_4,T_5$. We have $\bbb_2\ccc_3\cdots=\bbb\ccc^2$ or $\aaa^2\bbb\ccc$. If $\bbb_2\ccc_3\cdots=\aaa^2\bbb\ccc$, we get the AAD $\thick\aaa^{\bbb}\thin^{\bbb}\ccc\thin\cdots$ at $\aaa^2\bbb\ccc$. This gives a vertex $\bbb^2\cdots$, contradicting the AVC. Therefore, $\bbb_2\ccc_3\cdots=\thin^{\ccc}\bbb_2^{\aaa}\thin^{\ddd}\ccc_3^{\bbb}\thin^{\ddd}\ccc_6^{\bbb}\thin$ determines $T_6$. Then $\bbb_3\ddd_6\cdots=\bbb_3\ddd_6\ddd_7$ determines $T_7$. Similarly, we can determine $T_8,T_9$. Then we get $\aaa_3\aaa_4\ccc_7\ccc_8\cdots$, contradicting the AVC.
	   
	   In the right of Fig.\,\ref{fig 2ac2d 1}, $\aaa^2\ccc\ddd^2=\thin^{\bbb}\aaa_1^{\ddd}\thick^{\aaa}\ddd_2^{\ccc}\thin^{\bbb}\aaa_3^{\ddd}\thick\cdots$ determines $T_1,T_2,T_3$. We have $\bbb_3\ccc_2\cdots=\bbb\ccc^2$ or $\aaa^2\bbb\ccc$. If $\bbb_3\ccc_2\cdots=\aaa^2\bbb\ccc$, we get the AAD $\thick\aaa^{\bbb}\thin^{\bbb}\ccc\thin\cdots$ at $\aaa^2\bbb\ccc$. This gives a vertex $\bbb^2\cdots$, contradicting the AVC. Therefore,  $\bbb_3\ccc_2\cdots=\thin^{\ccc}\bbb_3^{\aaa}\thin^{\ddd}\ccc_2^{\bbb}\thin^{\ddd}\ccc_4^{\bbb}\thin$ determines $T_4$. Then $\bbb_2\ddd_4\cdots=\bbb_2\ddd_4\ddd_5$ determines $T_5$; $\thin\ccc_5\thin\aaa_2\thick\ddd_1\thin\cdots=\thin\ccc_5\thin\aaa_2\thick\ddd_1\thin\aaa_6\thick\ddd_7\thin$ determines $T_6,T_7$. Similarly, we can determine $T_8$. Then $\aaa_4\aaa_5\bbb_8\cdots=\thin^{\bbb}\aaa_4^{\ddd}\thick^{\ddd}\aaa_5^{\bbb}\thin^{\ccc}\bbb_8^{\aaa}\thin^{\bbb}\ccc_9^{\ddd}\thin$ determines $T_9$. We get $\bbb_4\ccc_3\ddd_9\cdots$, contradicting the AVC. 
	   
	   Therefore, $\aaa^2\ccc\ddd^2$ is not a vertex. This implies $\aaa\cdots=\aaa^2\bbb\ccc$.

	   \begin{figure}[htp]
	   	\centering
	   	\begin{tikzpicture}[>=latex,scale=0.5]

	   				\draw (0,0)--(2.5,0)--(4,2)--(5.5,0)--(4,-2)--(6,-2)--(8,0)--(8,-4)--(6,-7)--(4,-4)--(2,-7)--(0,-4)--(0,0)
	   				(0,0)--(2,-2)--(6,-2)--(6,-4)--(2,-4)--(2,-2)
	   				(11,1)--(11,-1)--(23,-1)
	   				(19,1)--(19,-1)
	   				(13,-1)--(11,-4)--(19,-4)--(19,-7)--(21,-4)
	   				(15,-1)--(17,-4)
	   				(19,-4)--(21,-1);
	   				
	   				\draw[line width=1.5]	
	   				(0,-4)--(2,-4)
	   				(2.5,0)--(4,-2)--(4,-4)
	   				(5.5,0)--(8,0)
	   				(6,-4)--(8,-4)
	   				(15,1)--(15,-1)--(13,-4)
	   				(15,-4)--(19,-7)
	   				(17,-1)--(19,-4)
	   				(23,1)--(23,-1)
	   				(21,-1)--(21,-4);

	   				\node at (0.65,-0.35) {\small $\ccc$}; \node at (2.3,-0.35) {\small $\ddd$};\node at (2.9,0) {\small $\ddd$};\node at (5,0) {\small $\bbb$};\node at (4,1.4) {\small $\ccc$};\node at (4,-1.4) {\small $\aaa$};\node at (5.5,-0.35) {\small $\ddd$};\node at (7.2,-0.35) {\small $\aaa$};\node at (0.35,-0.9) {\small $\bbb$};\node at (7.65,-0.9) {\small $\bbb$};
	   				\node at (1.65,-2.2) {\small $\ccc$};\node at (2.2,-1.65) {\small $\bbb$};\node at (3.3,-1.7) {\small $\aaa$};\node at (4.7,-1.7) {\small $\ccc$};\node at (5.9,-1.65) {\small $\bbb$};\node at (6.3,-2.2) {\small $\ccc$};\node at (2.3,-2.35) {\small $\ccc$};\node at (3.7,-2.35) {\small $\ddd$};\node at (4.3,-2.35) {\small $\ddd$};\node at (5.7,-2.35) {\small $\ccc$};
	   				\node at (0.3,-3.65) {\small $\aaa$};\node at (1.7,-3.65) {\small $\ddd$};\node at (2.3,-3.65) {\small $\bbb$};\node at (3.7,-3.65) {\small $\aaa$};\node at (4.3,-3.65) {\small $\aaa$};\node at (5.7,-3.65) {\small $\bbb$};\node at (6.3,-3.65) {\small $\ddd$};\node at (7.7,-3.65) {\small $\aaa$};
	   				\node at (0.6,-4.35) {\small $\aaa$};\node at (2,-4.35) {\small $\ddd$};\node at (3.4,-4.35) {\small $\ccc$};\node at (4.6,-4.35) {\small $\ccc$};\node at (6,-4.35) {\small $\ddd$};\node at (7.4,-4.35) {\small $\aaa$};\node at (2,-6.3) {\small $\bbb$};\node at (6,-6.3) {\small $\bbb$};

	   	        	\node[draw,shape=circle, inner sep=0.5] at (4,0) {\small $1$};
	   	        	\node[draw,shape=circle, inner sep=0.5] at (1.7,-0.8) {\small $2$};
	   	        	\node[draw,shape=circle, inner sep=0.5] at (3,-3) {\small $3$};
	   	        	\node[draw,shape=circle, inner sep=0.5] at (5,-3) {\small $4$};
	   	        	\node[draw,shape=circle, inner sep=0.5] at (6.3,-0.8) {\small $5$};
	   	        	\node[draw,shape=circle, inner sep=0.5] at (1,-3) {\small $6$};
	   	        	\node[draw,shape=circle, inner sep=0.5] at (2,-5.3) {\small $7$};
	   	        	\node[draw,shape=circle, inner sep=0.5] at (7,-3) {\small $8$};
	   	        	\node[draw,shape=circle, inner sep=0.5] at (6,-5.3) {\small $9$};

	   		        \node at (13,1) {\small $\aaa$};\node at (17,1) {\small $\ddd$};\node at (21,1) {\small $\aaa$};
	   		        \node at (11.3,-0.6) {\small $\bbb$};\node at (13,-0.6) {\small $\ccc$};\node at (14.7,-0.6) {\small $\ddd$};\node at (15.3,-0.6) {\small $\aaa$};\node at (17,-0.6) {\small $\bbb$};\node at (18.7,-0.6) {\small $\ccc$};\node at (19.3,-0.6) {\small $\bbb$};\node at (21,-0.6) {\small $\ccc$};\node at (22.7,-0.6) {\small $\ddd$};
	   		        \node at (13.1,-1.4) {\small $\bbb$};\node at (14.3,-1.4) {\small $\aaa$};\node at (15,-1.6) {\small $\ddd$};\node at (15.6,-1.4) {\small $\ccc$};\node at (16.8,-1.4) {\small $\ddd$};\node at (17.7,-1.4) {\small $\ddd$};\node at (19,-1.4) {\small $\ccc$};\node at (20.4,-1.4) {\small $\bbb$};\node at (20.7,-1.9) {\small $\ddd$};
	   		        \node at (11.7,-3.6) {\small $\ccc$};\node at (12.8,-3.6) {\small $\ddd$};\node at (13.7,-3.6) {\small $\aaa$};\node at (15,-3.6) {\small $\bbb$};\node at (16.3,-3.6) {\small $\ccc$};\node at (17.2,-3.6) {\small $\bbb$};\node at (18.3,-3.6) {\small $\aaa$};\node at (19,-3.4) {\small $\aaa$};\node at (19.3,-4.2) {\small $\ccc$};\node at (20.7,-3.9) {\small $\aaa$};\node at (16.1,-4.4) {\small $\ddd$};\node at (17,-4.4) {\small $\ccc$};\node at (18.6,-4.4) {\small $\bbb$};
	   		        \node at (18.6,-6.2) {\small $\aaa$};\node at (19.4,-6) {\small $\bbb$};
	   		
	   		        \node[draw,shape=circle, inner sep=0.5] at (13,0.2) {\small $1$};
	   		        \node[draw,shape=circle, inner sep=0.5] at (17,0.2) {\small $2$};
	   		        \node[draw,shape=circle, inner sep=0.5] at (21,0.2) {\small $3$};
	   		        \node[draw,shape=circle, inner sep=0.5] at (19,-2.4) {\small $4$};
	   		        \node[draw,shape=circle, inner sep=0.5] at (17,-2.4) {\small $5$};
	   		        \node[draw,shape=circle, inner sep=0.5] at (13,-2.4) {\small $6$};
	   		        \node[draw,shape=circle, inner sep=0.5] at (15,-2.4) {\small $7$};
	   		        \node[draw,shape=circle, inner sep=0.5] at (18,-5.2) {\small $8$};
	   		        \node[draw,shape=circle, inner sep=0.5] at (20,-4.5) {\small $9$};
	   		        
	   		        \fill (24,0) circle (0.05);\fill (24.3,0) circle (0.05);\fill (24.6,0) circle (0.05);
	   		        
	   	\end{tikzpicture}	
	   	\caption{Two possible $\aaa^2\ccc\ddd^2$ and their AAD.} \label{fig 2ac2d 1}	
	   \end{figure}
	   
	   If $\bbb\ccc^2$ is a vertex, it has only two possible AAD. In the first picture of Fig.\,\ref{fig 2ac2d 2}, $\bbb\ccc^2=\thin^{\ddd}\ccc_1^{\bbb}\thin^{\ddd}\ccc_2^{\bbb}\thin^{\aaa}\bbb_3^{\ccc}\thin$ determines $T_1,T_2,T_3$. Then $\thick\aaa_3\thin\bbb_2\thin\cdots=\thick^{\ddd}\aaa_3^{\bbb}\thin^{\ccc}\bbb_2^{\aaa}\thin^{\bbb}\ccc_5^{\ddd}\thin^{\bbb}\aaa_4^{\ddd}\thick$ determines $T_4,T_5$; $\thick\aaa_2\thin\bbb_5\thin\cdots=\thick^{\ddd}\aaa_2^{\bbb}\thin^{\ccc}\bbb_5^{\aaa}\thin^{\bbb}\ccc_6^{\ddd} 
	   \thin^{\bbb}\aaa_7^{\ddd}\thick$, $\bbb_1\ddd_2\ddd_7\cdots\\=\bbb_1\ddd_2\ddd_7$ determines $T_6,T_7$; $\bbb_7\ddd_6\cdots=\bbb_7\ddd_6\ddd_8$ determines $T_8$. We get $\aaa_1\ccc_7\ccc_8$ $\cdots$, contradicting the AVC.
	   
	   In the second picture of Fig.\,\ref{fig 2ac2d 2}, $\bbb\ccc^2=\thin^{\bbb}\ccc_1^{\ddd}\thin^{\ddd}\ccc_2^{\bbb}\thin^{\aaa}\bbb_3^{\ccc}\thin$ determines $T_1,T_2,T_3$. Then $\thick\aaa_3\thin\bbb_2\thin\cdots=\thick^{\ddd}\aaa_3^{\bbb}\thin^{\ccc}\bbb_2^{\aaa}\thin^{\bbb}\ccc_5^{\ddd}\thin^{\bbb}\aaa_4^{\ddd}\thick$ determines $T_4,T_5$; $\thick\aaa_2\thin\bbb_5\thin\cdots=\thick^{\ddd}\aaa_2^{\bbb}\thin^{\ccc}\bbb_5^{\aaa}\\\thin^{\bbb}\ccc_6^{\ddd}  
	   \thin^{\bbb}\aaa_7^{\ddd}\thick$  determines $T_6,T_7$; $\ddd_1\ddd_2\ddd_7\cdots=\ddd_1\ddd_2\ddd_7\ddd_8$ determines $T_8$; $\bbb_7\ddd_6\cdots=\bbb_7\ddd_6\ddd_9$ determines $T_9$; $\ccc_7\ccc_8\ccc_9\cdots=\ccc_7\ccc_8\ccc_9\ccc_{10}$. We get $\bbb_8\bbb_{10}\cdots$ or $\bbb_9\bbb_{10}\cdots$, contradicting the AVC. 
	   
	   Therefore, $\bbb\ccc^2$ is not a vertex.

       \begin{figure}[htp]
       	\centering
       	\begin{tikzpicture}[>=latex,scale=0.48]

       		\draw(0,0)--(0,-2)
       		(4,0)--(4,-2)--(2,-2)--(2,-6)--(4,-6)--(6,-4)--(6,-2)--(4,-4)--(2,-2)
       		(6,-2)--(8,-2)--(8,-4)--(10,-4)--(10,-2)--(12,-2)--(12,0)
       		(8,0)--(8,-2)
       		(15,0)--(15,-2)--(17,-2)--(17,-4)--(21,-4)--(21,-2)--(23,-4)--(21,-6)--(19,-6)--(19,0)
       		(21,-2)--(23,-2)--(23,0)
       		(23,-2)--(25,-4)--(27,-4)--(25,-2)--(27,-2)--(27,0);
       		
       		\draw[line width=1.5] (0,-2)--(2,-2)
       		(4,-2)--(6,-2)
       		(4,-4)--(4,-6)
       		(6,-4)--(8,-4)
       		(8,-2)--(10,-2)
       		(17,-2)--(21,-2)
       		(21,-4)--(21,-6)
       		(23,-2)--(25,-2)
       		(23,-4)--(25,-4);
       		
       		\draw[dotted] (17,-4)--(17,-6)--(19,-6);

       		\node at (2,0) {\small $\ccc$};\node at (6,0) {\small $\ccc$};\node at (10,0) {\small $\bbb$};
       		\node at (0.3,-1.6) {\small $\ddd$};\node at (2,-1.6) {\small $\aaa$};\node at (3.7,-1.6) {\small $\bbb$};\node at (4.3,-1.6) {\small $\ddd$};\node at (6,-1.6) {\small $\aaa$};\node at (7.7,-1.6) {\small $\bbb$};\node at (8.3,-1.6) {\small $\aaa$};\node at (10,-1.6) {\small $\ddd$};\node at (11.7,-1.6) {\small $\ccc$};
       		\node at (2.3,-2.9) {\small $\ccc$};\node at (2.7,-2.35) {\small $\ccc$};\node at (4,-2.35) {\small $\ddd$};\node at (5.3,-2.35) {\small $\aaa$};\node at (5.7,-2.9) {\small $\ccc$};\node at (4,-3.5) {\small $\bbb$};\node at (6.3,-2.4) {\small $\bbb$};\node at (7.7,-2.4) {\small $\ccc$};\node at (8.3,-2.4) {\small $\aaa$};\node at (9.7,-2.4) {\small $\ddd$};
       		\node at (2.3,-5.6) {\small $\bbb$};\node at (3.7,-5.6) {\small $\aaa$};\node at (3.7,-4.3) {\small $\ddd$};\node at (4.3,-4.3) {\small $\ddd$};\node at (4.3,-5.2) {\small $\aaa$};\node at (5.7,-3.9) {\small $\bbb$};\node at (6.3,-3.6) {\small $\aaa$};\node at (7.7,-3.6) {\small $\ddd$};\node at (8.3,-3.6) {\small $\bbb$};\node at (9.7,-3.6) {\small $\ccc$};

       		\node[draw,shape=circle, inner sep=0.2] at (2,-0.8) {\small $1$};
       		\node[draw,shape=circle, inner sep=0.2] at (6,-0.8) {\small $2$};
       		\node[draw,shape=circle, inner sep=0.2] at (10,-0.8) {\small $3$};
       		\node[draw,shape=circle, inner sep=0.2] at (9,-3) {\small $4$};
       		\node[draw,shape=circle, inner sep=0.2] at (7,-3) {\small $5$};           
            \node[draw,shape=circle, inner sep=0.2] at (5,-4) {\small $6$};
            \node[draw,shape=circle, inner sep=0.2] at (4.5,-2.8) {\small $7$};
            \node[draw,shape=circle, inner sep=0.2] at (2.7,-4.3) {\small $8$};

       		\node at (17,0) {\small $\ccc$};\node at (21,0) {\small $\ccc$};\node at (25,0) {\small $\bbb$};
       		\node at (15.3,-1.6) {\small $\bbb$};\node at (17,-1.6) {\small $\aaa$};\node at (18.7,-1.6) {\small $\ddd$};\node at (19.3,-1.6) {\small $\ddd$};\node at (21,-1.6) {\small $\aaa$};\node at (22.7,-1.6) {\small $\bbb$};\node at (23.3,-1.6) {\small $\aaa$};\node at (25,-1.6) {\small $\ddd$};\node at (26.7,-1.6) {\small $\ccc$};
       		\node at (17.3,-2.4) {\small $\aaa$};\node at (18.7,-2.4) {\small $\ddd$};\node at (19.3,-2.4) {\small $\ddd$};\node at (20.7,-2.4) {\small $\aaa$};\node at (21.3,-2.8) {\small $\ccc$};\node at (22,-2.4) {\small $\bbb$};\node at (22.9,-2.4) {\small $\ccc$};\node at (23.8,-2.4) {\small $\aaa$};\node at (24.9,-2.4) {\small $\ddd$};
       		\node at (17.3,-3.6) {\small $\bbb$};\node at (18.7,-3.6) {\small $\ccc$};\node at (19.3,-3.6) {\small $\ccc$};\node at (20.7,-3.6) {\small $\bbb$};\node at (21.3,-4) {\small $\ddd$};\node at (22.5,-4.1) {\small $\bbb$};\node at (21.3,-5.2) {\small $\aaa$};\node at (23.1,-3.6) {\small $\aaa$};\node at (24.1,-3.6) {\small $\ddd$};\node at (25.2,-3.6) {\small $\bbb$};\node at (26.1,-3.6) {\small $\ccc$};
       		\node at (18.7,-4.4) {\small $\ccc$};\node at (19.3,-4.4) {\small $\ccc$};\node at (20.7,-4.4) {\small $\ddd$};\node at (19.3,-5.6) {\small $\bbb$};\node at (20.7,-5.6) {\small $\aaa$};
       		
       		\node[draw,shape=circle, inner sep=0.2] at (17,-0.8) {\small $1$};
       		\node[draw,shape=circle, inner sep=0.2] at (21,-0.8) {\small $2$};
       		\node[draw,shape=circle, inner sep=0.2] at (25,-0.8) {\small $3$};
       		\node[draw,shape=circle, inner sep=0.2] at (24.6,-3) {\small $4$};
       		\node[draw,shape=circle, inner sep=0.2] at (22.6,-3) {\small $5$};
       		\node[draw,shape=circle, inner sep=0.2] at (21.9,-4) {\small $6$};
       		\node[draw,shape=circle, inner sep=0.2] at (20,-3) {\small $7$};
       		\node[draw,shape=circle, inner sep=0.2] at (18,-3) {\small $8$};
       		\node[draw,shape=circle, inner sep=0.2] at (20,-5) {\small $9$};
       		\node[draw,shape=circle, inner sep=0.2] at (18,-5) {\footnotesize $10$};

       	\end{tikzpicture}	
       	\caption{Two possible AAD for $\bbb\ccc^2$.} \label{fig 2ac2d 2}	
       \end{figure}
	   
	   	This implies $\text{AVC}\sub\{\bbb\ddd^2,\aaa^2\bbb\ccc,\ccc^4,\ccc^2\ddd^2,\ddd^4\}$. There is only one solution satisfying Balance Lemma \ref{balance}: $\{8\bbb\ddd^2,8\aaa^2\bbb\ccc,2\ccc^4\}$.	   	   
	    We have the AAD $\ccc^4=\thin^{\bbb}\ccc_1^{\ddd}\thin^{\bbb}\ccc_2^{\ddd}\thin^{\bbb}\ccc_3^{\ddd}\thin^{\bbb}\ccc_4^{\ddd}\thin$ which determines $T_1,T_2,T_3$, $T_4$. Then $\bbb_2\ddd_1\cdots=\bbb_2\ddd_1\ddd_5$ determines $T_5$. Similarly, we can determine $T_6,T_7,T_8$. 	    	    
	    Then $\aaa_2\aaa_6\ccc_5\cdots=\thin^{\bbb}\aaa_2^{\ddd}\thick^{\ddd}\aaa_6^{\bbb}\thin^{\ccc}\bbb_9^{\aaa}\thin^{\bbb}\ccc_5^{\ddd}\thin$. So $\bbb_5\cdots=\aaa_9\bbb_5\cdots$ or $\bbb_5\ccc_9\cdots$, shown in two pictures  of Fig.\,\ref{fig 7-1}  respectively.
	    
	    \begin{figure}[htp]
	    	\centering
	    	\begin{tikzpicture}[>=latex,scale=0.45]
	    		
	    		\foreach \b in {0,1,2,3}
	    		{
	    			\begin{scope}[xshift=3*\b cm]	
	    				\draw (0,0)--(0,-2)--(1.5,-2)--(3,-5)--(4.5,-2)--(3,-2)--(3,0)
	    				(4.5,-2)--(6,-5)--(4.5,-5)--(4.5,-7)
	    				(1.5,-7)--(1.5,-5)--(3,-5);
	    				
	    				\draw[line width=1.5]	
	    				(1.5,-2)--(3,-2)
	    				(3,-5)--(4.5,-5);

	    				\node at (1.5,0.2) {\small $\ccc$};\node at (1.5,-1.6) {\small $\aaa$};
	    				\node at (0.35,-1.6) {\small $\bbb$}; \node at (2.65,-1.6) {\small $\ddd$};
	    				\node at (3,-2.4) {\small $\ddd$}; \node at (2.1,-2.4) {\small $\aaa$};
	    				\node at (3.9,-2.5) {\small $\ccc$}; \node at (4.5,-2.8) {\small $\bbb$};
	    				\node at (3.05,-4.2) {\small $\bbb$}; \node at (4.5,-4.6) {\small $\ddd$};
	    				\node at (3.6,-4.6) {\small $\aaa$}; \node at (5.4,-4.6) {\small $\ccc$};
	    				\node at (1.9,-5.4) {\small $\bbb$};  \node at (3,-5.4) {\small $\aaa$}; 
	    				\node at (4.1,-5.4) {\small $\ddd$};  \node at (3,-7.3) {\small $\ccc$}; 
	    				
	    			\end{scope}
	    		}

	    		\node[draw,shape=circle, inner sep=0.5] at (1.5,-0.75) {\small $1$};
	    		\node[draw,shape=circle, inner sep=0.5] at (4.5,-0.75) {\small $2$};
	    		\node[draw,shape=circle, inner sep=0.5] at (7.5,-0.75) {\small $3$};
	    		\node[draw,shape=circle, inner sep=0.5] at (10.5,-0.75) {\small $4$};
	    		
	    		\node[draw,shape=circle, inner sep=0.5] at (3,-3.2) {\small $5$};
	    		\node[draw,shape=circle, inner sep=0.5] at (3+3,-3.2) {\small $6$};
	    		\node[draw,shape=circle, inner sep=0.5] at (3+6,-3.2) {\small $7$};
	    		\node[draw,shape=circle, inner sep=0.5] at (3+9,-3.2) {\small $8$};
	    		
	    		\node[draw,shape=circle, inner sep=0.5] at (4.5,-3.7) {\small $9$};
	    		\node[draw,shape=circle, inner sep=0.5] at (4.5+3,-3.7) {\footnotesize $10$};
	    		\node[draw,shape=circle, inner sep=0.5] at (4.5+6,-3.7) {\footnotesize $11$};
	    		\node[draw,shape=circle, inner sep=0.5] at (4.5+9,-3.7) {\footnotesize $12$};
	    		
	    		\node[draw,shape=circle, inner sep=0.5] at (3,-6.3) {\footnotesize $13$};
	    		\node[draw,shape=circle, inner sep=0.5] at (3+3,-6.3) {\footnotesize $14$};
	    		\node[draw,shape=circle, inner sep=0.5] at (3+6,-6.3) {\footnotesize $15$};
	    		\node[draw,shape=circle, inner sep=0.5] at (3+9,-6.3) {\footnotesize $16$};

	    	\begin{scope}[xshift=15 cm]
	    		
	    		\foreach \b in {0,1,2,3}
	    		{
	    			\begin{scope}[xshift=3*\b cm]	
	    				\draw (0,0)--(0,-2)--(1.5,-2)--(3,-5)--(4.5,-2)--(3,-2)--(3,0)
	    				(4.5,-2)--(6,-5)--(4.5,-5)--(4.5,-7)
	    				(1.5,-7)--(1.5,-5)--(3,-5)
	    				(3,-5)--(4.5,-5);
	    				
	    				\draw[line width=1.5]	
	    				(1.5,-2)--(3,-2)
	    				(4.5,-5)--(6,-5)
	    				(1.5,-5)--(3,-5);

	    				\node at (1.5,0.2) {\small $\ccc$};\node at (1.5,-1.6) {\small $\aaa$};
	    				\node at (0.35,-1.6) {\small $\bbb$}; \node at (2.65,-1.6) {\small $\ddd$};
	    				\node at (3,-2.4) {\small $\ddd$}; \node at (2.1,-2.4) {\small $\aaa$};
	    				\node at (3.9,-2.5) {\small $\ccc$}; \node at (4.5,-2.8) {\small $\bbb$};
	    				\node at (3.05,-4.2) {\small $\bbb$}; \node at (4.5,-4.6) {\small $\ddd$};
	    				\node at (3.6,-4.6) {\small $\ccc$}; \node at (5.4,-4.6) {\small $\aaa$};
	    				\node at (1.9,-5.4) {\small $\ddd$};  \node at (3,-5.4) {\small $\aaa$}; 
	    				\node at (4.1,-5.4) {\small $\bbb$};  \node at (3,-7.3) {\small $\ccc$}; 
	    				
	    			\end{scope}
	    		}

	    		\node[draw,shape=circle, inner sep=0.5] at (1.5,-0.75) {\small $1$};
	    		\node[draw,shape=circle, inner sep=0.5] at (4.5,-0.75) {\small $2$};
	    		\node[draw,shape=circle, inner sep=0.5] at (7.5,-0.75) {\small $3$};
	    		\node[draw,shape=circle, inner sep=0.5] at (10.5,-0.75) {\small $4$};
	    		
	    		\node[draw,shape=circle, inner sep=0.5] at (3,-3.2) {\small $5$};
	    		\node[draw,shape=circle, inner sep=0.5] at (3+3,-3.2) {\small $6$};
	    		\node[draw,shape=circle, inner sep=0.5] at (3+6,-3.2) {\small $7$};
	    		\node[draw,shape=circle, inner sep=0.5] at (3+9,-3.2) {\small $8$};
	    		
	    		\node[draw,shape=circle, inner sep=0.5] at (4.5,-3.7) {\small $9$};
	    		\node[draw,shape=circle, inner sep=0.5] at (4.5+3,-3.7) {\footnotesize $10$};
	    		\node[draw,shape=circle, inner sep=0.5] at (4.5+6,-3.7) {\footnotesize $11$};
	    		\node[draw,shape=circle, inner sep=0.5] at (4.5+9,-3.7) {\footnotesize $12$};
	    		
	    		\node[draw,shape=circle, inner sep=0.5] at (3,-6.3) {\footnotesize $13$};
	    		\node[draw,shape=circle, inner sep=0.5] at (3+3,-6.3) {\footnotesize $14$};
	    		\node[draw,shape=circle, inner sep=0.5] at (3+6,-6.3) {\footnotesize $15$};
	    		\node[draw,shape=circle, inner sep=0.5] at (3+9,-6.3) {\footnotesize $16$};

	    	\end{scope}
	    	\end{tikzpicture}
	    	\caption{Two tilings for $\{8\bbb\ddd^2,8\aaa^2\bbb\ccc,2\ccc^4\}$.} \label{fig 7-1}	
	    \end{figure}   
    
        \begin{figure}[htp]
        	\centering
        	
        	\begin{tikzpicture}[>=latex,scale=0.65]
        		
        		\foreach \b in {0,1,2,3}
        		{
        			\begin{scope}[rotate=90*\b]	
        				\fill[gray!50]  (0,1)--(0,{1+sqrt(2)}) -- (0,{1+sqrt(2)}) arc ((-45+135:90+135):{sqrt(2)})--(-1,0) arc (135+45:45+45:1);
        				
        				\fill[gray!50]  ({-sqrt(2)/2},{sqrt(2)/2}) arc (135+45:180+45:{sqrt(2)})--({-sqrt(2)/2},{-sqrt(2)/2})--({-sqrt(2)/2},{sqrt(2)/2}) arc (135:225:1);

        				\draw (0,0)--({-sqrt(2)/2},{sqrt(2)/2})
        				(0,1)--(0,3);
        				\draw ({-sqrt(2)/2},{sqrt(2)/2}) arc (135:225:1);
        				\draw[line width=1.5]	({-sqrt(2)/2},{sqrt(2)/2}) arc (135+45:180+45:{sqrt(2)});   			
        				\draw[line width=1.5]	(0,{1+sqrt(2)}) arc ((-45+135:90+135):{sqrt(2)});  	   				
        			\end{scope}
        		}
        		
        	\begin{scope}[xshift=8cm]
        		
        		\foreach \b in {0,1,2,3}
        		{
        			\begin{scope}[rotate=90*\b]	
        				\fill[gray!50]  (0,1)--(0,{1+sqrt(2)}) -- (0,{1+sqrt(2)}) arc ((-45+135:90+135):{sqrt(2)})--(-1,0) arc (135+45:45+45:1);
        				
        				\fill[gray!50]  (0,0)--({-sqrt(2)/2},{sqrt(2)/2})--({-sqrt(2)/2},{sqrt(2)/2}) arc (180+270:135+270:{sqrt(2)})--(0,0);

        				\draw (0,0)--({-sqrt(2)/2},{sqrt(2)/2})
        				(0,1)--(0,3);
        				\draw ({-sqrt(2)/2},{sqrt(2)/2}) arc (135:225:1);
        				\draw[line width=1.5]	({-sqrt(2)/2},{sqrt(2)/2}) arc (180+270:135+270:{sqrt(2)});   			
        				\draw[line width=1.5]	(0,{1+sqrt(2)}) arc ((-45+135:90+135):{sqrt(2)});  	   				
        			\end{scope}
        		}
        	\end{scope}	
        	\end{tikzpicture}	    
        	\caption{Stereo-graphic projection for two tilings of $\{8\bbb\ddd^2,8\aaa^2\bbb\ccc,2\ccc^4\}$. } \label{fig 1-1-4}	
        \end{figure}
    
	    In the left of Fig.\,\ref{fig 7-1}, $\bbb_5\cdots=\aaa_9\bbb_5\cdots$
	    determines $T_9$. Then $\aaa_3\aaa_7\ccc_6\cdots=\aaa_3\aaa_7\bbb_{10}\ccc_6,\bbb_6\ccc_9\cdots=\aaa_{10}\bbb_6\ccc_9\cdots$ determine $T_{10}$. Similarly, we can determine $T_{11},T_{12}$. Then $ \aaa_9\bbb_5\ccc_{12}\cdots=\aaa_9\aaa_{13}\bbb_5\ccc_{12}$ determines $T_{13}$. Similarly, we can determine $T_{14},T_{15},T_{16}$. 
	    
	    In the right of Fig.\,\ref{fig 7-1}, $\bbb_5\cdots=\bbb_5\ccc_9\cdots$ determines $T_9$. Then we get a different tiling by similar deductions. 
	    The $3$D pictures for these two tilings are shown in Fig.\,\ref{fig 1-1-1}.  Their authentic pictures of the stereo-graphic projection are shown in Fig.\,\ref{fig 1-1-4}.
	    This is Case $(1,4,2,2)/4$ in Table \ref{Tab-1.1}. \label{discrete-15}

	   

\newpage	

		\section*{Appendix: Exact \& numerical geometric data}

     \begin{table*}[htp]                    	
     	\centering 
     	\resizebox{\textwidth}{75mm}{\begin{tabular}{c|c}	 
     			
     			angles $ (\aaa,\bbb,\ccc,\ddd)$  & edges \\
     			\hline 			 
     			$(6,3,4,3)/6$&$a=1/2,b=1/6$	 \\
     			\hline
     			$(1,8,4,3)/6$&$a=\arccos(1/3)\approx0.3918,b=1$\\     			
     			\hline
     			\multirow{3}{*}{$(12,4,6,2)/9$}&$a=1-\arcsin \! \left(\frac{\sqrt{2}}{\sqrt{\cos \! \left(\frac{2 \pi}{9}\right)}\, \left(2-2 \cos \! \left(\frac{4 \pi}{9}\right)\right)} \! \right) \approx0.5673$\\
     			&$b= \arccos \! \left(\frac{\sqrt{3} \cot \! \left(\frac{2 \pi}{9}\right) -\cot \! \left(\frac{2 \pi}{9}\right) \sin \! \left(\frac{\pi}{9}\right)}{1+\cos \! \left(\frac{\pi}{9}\right)}\right)\approx0.1741$\\
     			\hline
     			\multirow{3}{*}{$(2,10,3,6)/9$}&$a=\arccos \! \left(\frac{4\sqrt{3} \sin(\frac{2 \pi}{9}) }{3}-1 \! \right)\approx 0.3390 $\\
     			&$b=\arccos \! \left(\frac{8 \cos(\frac{\pi}{9})-4 \sqrt{3}\, \sin(\frac{4 \pi}{9})-1}{3} \! \right )\approx0.5324$\\
     			\hline
     			\multirow{3}{*}{$(1,21,5,8)/15$}&$a=\arccos \! \left(\frac{2 \sin \! \left(\frac{\pi}{15}\right)-\sqrt{3}\, \cos \! \left(\frac{7 \pi}{15}\right)}{\sin \! \left(\frac{7 \pi}{15}\right)}\right) \approx0.4241$\\
     			&$b=\arccos \! \left(\frac{51-90 \sqrt{3}\, \sin \! \left(\frac{2 \pi}{5}\right)-96 \sqrt{3}\, \sin \! \left(\frac{7 \pi}{15}\right)+88 \cos \! \left(\frac{2 \pi}{15}\right)+184 \cos \! \left(\frac{\pi}{15}\right)}{1+6 \cos \! \left(\frac{7 \pi}{15}\right)-2 \cos \! \left(\frac{2 \pi}{15}\right)+6 \cos \! \left(\frac{2 \pi}{5}\right)+2 \cos \! \left(\frac{\pi}{5}\right)}\right)\approx0.7413$\\
     			\hline	
     			\multirow{3}{*}{$(4,9,5,17)/15$}&$a=\arccos \! \left(\frac{2 \sin \! \left(\frac{\pi}{15}\right)-\sqrt{3}\, \cos \! \left(\frac{7 \pi}{15}\right)}{\sin \! \left(\frac{7 \pi}{15}\right)}\right) \approx0.4241$\\
     			&$b=\arccos \! \left( \frac{-3+9 \sqrt{5}-5 \sqrt{3}\, \sqrt{10-2 \sqrt{5}}}{-9-9 \sqrt{5}+\sqrt{3}\, \left(\sqrt{5}+4\right) \sqrt{10-2 \sqrt{5}}} \right)\approx0.1654$\\
     			\hline
     			\multirow{3}{*}{$(9,28,10,23)/30$}&$a= \arccos \! \left(\frac{\cot \! \left(\frac{\pi}{10}\right)-2 \cot \! \left(\frac{\pi}{10}\right) \cos \! \left(\frac{7 \pi}{30}\right) \sin \! \left(\frac{\pi}{5}\right)}{2 \sin \! \left(\frac{\pi}{5}\right) \sin \! \left(\frac{7 \pi}{30}\right)}\right)\approx0.3353$\\
     			&$b=\arccos \! \left( \frac{30+2 \sqrt{5}-\sqrt{3}\, \left(5+\sqrt{5}\right) \sqrt{10-2 \sqrt{5}}}{2-10 \sqrt{5}+3 \sqrt{3}\, \left(\sqrt{5}+1\right) \sqrt{10-2 \sqrt{5}}} \right)\approx0.4159$\\
     			\hline
     			\multirow{3}{*}{$(3,16,10,41)/30$}&$a= \arccos \! \left(\frac{\sqrt{3} \cos \! \left(\frac{11 \pi}{30}\right) -2 \sin \! \left(\frac{\pi}{10}\right)}{\sin \! \left(\frac{11 \pi}{30}\right)}\right)\approx0.4698$\\
     			&$b=\arccos \! \left(\frac{-28+60 \sqrt{3}\, \sin \! \left(\frac{7 \pi}{15}\right)+61 \sqrt{3}\, \sin \! \left(\frac{4 \pi}{15}\right)+61 \sqrt{3}\, \sin \! \left(\frac{\pi}{15}\right)-61 \cos \! \left(\frac{2 \pi}{15}\right)-120 \cos \! \left(\frac{\pi}{15}\right)}{\cos \! \left(\frac{2 \pi}{5}\right)+3 \cos \! \left(\frac{2 \pi}{15}\right)}\right)\approx0.1461$\\
     			\hline
     			\multirow{3}{*}{$(5,32,6,23)/30$}&$a= \arccos \! \left(\frac{\cot \! \left(\frac{\pi}{10}\right)-2 \cot \! \left(\frac{\pi}{10}\right) \cos \! \left(\frac{7 \pi}{30}\right) \sin \! \left(\frac{\pi}{5}\right)}{2 \sin \! \left(\frac{\pi}{5}\right) \sin \! \left(\frac{7 \pi}{30}\right)}\right)\approx0.3353$\\
     			&$b=\arccos \! \left( \frac{30+2 \sqrt{5}-\sqrt{3}\, \left(5+\sqrt{5}\right) \sqrt{10-2 \sqrt{5}}}{2-10 \sqrt{5}+3 \sqrt{3}\, \left(\sqrt{5}+1\right) \sqrt{10-2 \sqrt{5}}} \right)\approx0.4159$\\
     			\hline
     			\multirow{3}{*}{$(1,16,6,43)/30$}&$a= \arccos \! \left(\frac{\sqrt{3} \cos \! \left(\frac{11 \pi}{30}\right) -2 \sin \! \left(\frac{\pi}{10}\right)}{\sin \! \left(\frac{11 \pi}{30}\right)}\right)\approx0.4698$\\
     			&$b=\arccos \! \left( \frac{-7 \sqrt{3}+22 \sqrt{3}\, \cos \! \left(\frac{\pi}{15}\right)-24 \sqrt{3}\, \cos \! \left(\frac{2 \pi}{15}\right)+32 \sin \! \left(\frac{7 \pi}{15}\right)-18 \sin \! \left(\frac{2 \pi}{5}\right)}{21 \sqrt{3}-66 \sqrt{3}\, \cos \! \left(\frac{\pi}{15}\right)+80 \sqrt{3}\, \cos \! \left(\frac{2 \pi}{15}\right)-104 \sin \! \left(\frac{7 \pi}{15}\right)+58 \sin \! \left(\frac{2 \pi}{5}\right)}\right)\approx0.2730$\\
     			\hline

     	\end{tabular}}
     \end{table*}                          
     \begin{table*}[htp]                   	     	
     	\centering          
     	\resizebox{\textwidth}{70mm}{\begin{tabular}{c|c}	      			
     			angles $ (\aaa,\bbb,\ccc,\ddd)$  & edges \\
     			\hline 	
     			\multirow{3}{*}{$(1,42,4,17)/30$}&$a=\arccos \! \left(\frac{2 \sin \! \left(\frac{\pi}{15}\right)-\sqrt{3}\, \cos \! \left(\frac{7 \pi}{15}\right)}{\sin \! \left(\frac{7 \pi}{15}\right)}\right) \approx0.4241$\\
     			&$b=\arccos \! \left(\frac{\sqrt{3}\, \left(9 \sqrt{5}+29\right) \sqrt{10-2 \sqrt{5}}-58 \sqrt{5}-70}{\left(15 \sqrt{5}+27\right) \sqrt{3}\, \sqrt{10-2 \sqrt{5}}-46 \sqrt{5}-146}\right)\approx0.5493$\\
     			\hline
     			\multirow{3}{*}{$(3,20,4,13)/18$}&$a=\arccos \! \left(\frac{4\sqrt{3} \sin(\frac{2 \pi}{9}) }{3}-1 \! \right)\approx 0.3390 $\\
     			&$b=\arccos \! \left(\frac{\cos \! \left(\frac{\pi}{9}\right)-1}{2 \sqrt{3}\, \sin \! \left(\frac{4 \pi}{9}\right)-3 \cos \! \left(\frac{\pi}{9}\right)-1}\right)\approx0.4527$\\
     			\hline
     			$(1,4,2,2)/4$&$a=1/4,b=1/2$\\
     			\hline	     			
     			\multirow{3}{*}{$(5,4,7,3)/9$}&$a= \arccos \! \left(\frac{\sqrt{3} \cot \! \left(\frac{2 \pi}{9}\right) -\cot \! \left(\frac{2 \pi}{9}\right) \sin \! \left(\frac{\pi}{9}\right)}{1+\cos \! \left(\frac{\pi}{9}\right)}\right)\approx0.1741$\\
     			&$b=\arccos \! \left(\frac{68 \sqrt{3}+47\sqrt{3} \cos \! \left(\frac{\pi}{9}\right) +162 \sin \! \left(\frac{2 \pi}{9}\right)+162 \sin \! \left(\frac{\pi}{9}\right)}{99 \sqrt{3}+69\sqrt{3} \cos \! \left(\frac{\pi}{9}\right) +234 \sin \! \left(\frac{2 \pi}{9}\right)+234 \sin \! \left(\frac{\pi}{9}\right)} \right)\approx0.2584$\\
     			\hline
     			\multirow{2}{*}{$(15,6,10,7)/18$}&$a= \arccos \! \left(4 \cos \! \left(\frac{\pi}{9}\right)-3\right)\approx0.2258$\\
     			&$b=\arccos \! \left(28\sqrt{3} \sin \! \left(\frac{4 \pi}{9}\right)-36 \cos \! \left(\frac{\pi}{9}\right)-13 \right)\approx0.1183$\\
     			\hline	
     			\hline
     			\multirow{2}{*}{$(\frac4f,1-\frac4f,\frac4f,1)$}&$a=\arccos \! \left(\frac{\cos \! \left(\frac{4\pi}{f}\right) \left(1-\cos \! \left(\frac{4\pi}{f}\right)\right)}{\sin^2 \! \left(\frac{4\pi}{f}\right)}\right),b=1-2a$\\
     			&$f=10,a\approx0.4241,b\approx0.1517;\,\,\lim_{f\to\infty}a=\lim_{f\to\infty}b=1/3$\\
     			\hline	
     			\multirow{4}{*}{$(\frac2f,\frac{4f-4}{3f},\frac4f,\frac{2f-2}{3f})$}&$a= \arccos \! \left(\frac{\sqrt{3}\, \sin \! \left(\frac{8 \pi}{3 f}\right)-\sqrt{3}\, \sin \! \left(\frac{4 \pi}{3 f}\right)-\cos \! \left(\frac{4 \pi}{3 f}\right)-\cos \! \left(\frac{8 \pi}{3 f}\right)+2}{\sqrt{3}\, \sin \! \left(\frac{8 \pi}{3 f}\right)+\sqrt{3}\, \sin \! \left(\frac{4 \pi}{3 f}\right)+\cos \! \left(\frac{4 \pi}{3 f}\right)-\cos \! \left(\frac{8 \pi}{3 f}\right)}\right)$\\
     			&$b=\arccos \! \left(\frac{\sqrt{3}\, \sin \! \left(\frac{2 \pi}{3 f}\right)+4 \cos \! \left(\frac{2 \pi}{f}\right)-\cos \! \left(\frac{2 \pi}{3 f}\right)}{\sqrt{3}\, \sin \! \left(\frac{2 \pi}{3 f}\right)+3 \cos \! \left(\frac{2 \pi}{3 f}\right)}\right)+ 
     			\arccos \! \left(\frac{\sqrt{3}\, \left(\cos \! \left(\frac{2 \pi}{f}\right)-\cos \! \left(\frac{2 \pi}{3 f}\right)+\sqrt{3}\, \sin \! \left(\frac{2 \pi}{3 f}\right)\right)}{3 \sin \! \left(\frac{2 \pi}{f}\right)}\right)$\\
     			&$f=6,a\approx0.3390,b\approx0.8065;\,\,\lim_{f\to\infty}a=\lim_{f\to\infty}b=\arccos(1/3)$\\
     			\hline
     			\multirow{6}{*}{$(\frac2f,\frac{2f-4}{3f},\frac4f,\frac{4f-2}{3f})$}&$a= \arccos \! \left(\frac{\sqrt{3}\,\sin \! \left(\frac{2 \pi}{3 f}\right)  \cos \! \left(\frac{2 \pi}{f}\right)+\cos \! \left(\frac{2 \pi}{3 f}\right) \cos \! \left(\frac{2 \pi}{f}\right)-1}{\sin \! \left(\frac{2 \pi}{f}\right) \left(\sqrt{3}\cos \! \left(\frac{2 \pi}{3 f}\right) -\sin \! \left(\frac{2 \pi}{3 f}\right)\right)}\right)$\\
     			&$\phi =\arccos \! \left(\frac{\sin \! \left(\frac{2 \pi}{f}\right)-\sin \! \left(\frac{\left(f +4\right) \pi}{6 f}\right) \sin \! \left(\frac{4 \pi}{f}\right)}{\sqrt{-2 \sin \! \left(\frac{4 \pi}{f}\right) \sin \! \left(\frac{2 \pi}{f}\right) \sin \! \left(\frac{\left(f +4\right) \pi}{6 f}\right)-\cos \! \left(\frac{2 \pi}{f}\right)^{2}-\cos \! \left(\frac{4 \pi}{f}\right)^{2}+2}}\right)$\\
     			&$b=\arcsin \! \left(\frac{\sin \! a  \sin \! \left(-\frac{2\pi}{3}+\frac{4\pi}{3 f}+\phi \right)}{\sin \! \left(-\frac{4\pi}{3}+\frac{2\pi}{3 f}+\phi \right)}\right)$\\
     			&$f=10,a\approx0.4698,b\approx0.0898;\,\,\lim_{f\to\infty}a=\lim_{f\to\infty}b=\arccos(1/3)$\\
     			\hline
     	\end{tabular}}
     \end{table*}

\end{document}